\documentclass[a4paper,10pt]{amsart}
\usepackage[foot]{amsaddr}
\usepackage{geometry}
\newtheorem{theorem}{Theorem}[section]

\theoremstyle{remark}
\newtheorem{remark}[theorem]{Remark}
\usepackage{graphicx}
\usepackage{lineno}
\usepackage{bm}
\usepackage{amsfonts} 
\usepackage{amsmath}
\usepackage{amssymb}
\usepackage{framed} 
\usepackage{multicol} 
\usepackage{tabularx}
\usepackage{tabulary}
\usepackage{siunitx}
\usepackage{caption}
\usepackage{float}
\usepackage{url} 
\usepackage{multicol}
\usepackage{breakurl} 
\usepackage[breaklinks]{hyperref} 

\usepackage{subcaption}
\usepackage{xpatch}
\usepackage{nomencl} 
\makenomenclature
\makenomenclature 
\RequirePackage{ifthen}
\let\oldref\ref
\renewcommand{\ref}[1]{(\oldref{#1})}

\usepackage{amssymb}
\providecommand{\norm}[1]{\lVert#1\rVert}

\newcommand{\ltwonorm}[1]{\langle #1 \rangle_{L_2(\Omega)}}

\DeclareMathOperator{\spn}{span}
\DeclareMathAlphabet\mathbfcal{OMS}{cmsy}{b}{n}

\xpatchcmd{\thenomenclature}{%
  \section*{\nomname}
}{
}{\typeout{Success}}{\typeout{Failure}}

\usepackage{ifthen}
\renewcommand{\nomgroup}[1]{%
  \ifthenelse{\equal{#1}{A}}{\item[\textbf{Abbreviations}]}{%
    \ifthenelse{\equal{#1}{G}}{\item[\textbf{Symbols}]}{%
      \ifthenelse{\equal{#1}{C}}{\item[\textbf{Abbreviations}]}{%
        \ifthenelse{\equal{#1}{S}}{\item[\textbf{Subscripts}]}{%
          \ifthenelse{\equal{#1}{Z}}{\item[\textbf{Mathematical Symbols}]}{}
        }
      }
    }
  }
}

\begin{document}
\nomenclature{$\text{ROM}$}{Reduced Order Model}
\nomenclature{$\text{SUP-ROM}$}{ROM with pressure supremizer stabilisation}
\nomenclature{$\text{PPE-ROM}$}{ROM with pressure Poisson equation stabilisation}
\nomenclature{$\text{HF}$}{High Fidelity}
\nomenclature{$\text{POD}$}{Proper Orthogonal Decomposition}

\nomenclature[G]{$\bm{u}$}{velocity field}
\nomenclature[G]{$\bm{f}$}{Dirichlet boundary condition for velocity}
\nomenclature[G]{$\bm{k}$}{initial condition for velocity}
\nomenclature[G]{${p}$}{pressure field}
\nomenclature[G]{${\nu}$}{dimensionless kinematic viscosity}
\nomenclature[G]{${N_u^h}$}{number of unknowns for velocity at full-order level}
\nomenclature[G]{${N_p^h}$}{number of unknowns for pressure at full-order level}
\nomenclature[G]{${N_u^r}$}{number of unknowns for velocity at reduced order level}
\nomenclature[G]{${N_p^r}$}{number of unknowns for pressure at reduced order level}
\nomenclature[G]{$N_s$}{number of snapshots}
\nomenclature[G]{$Q$}{space-time domain}
\nomenclature[G]{$T$}{final time}
\nomenclature[G]{$\mathcal{T}$}{tessellation}
\nomenclature[G]{$\mathcal{P}$}{parameter space}
\nomenclature[G]{$\mathcal{K}$}{training set space}
\nomenclature[G]{${\Omega}$}{bounded domain}
\nomenclature[G]{${\Gamma}$}{boundary of $\Omega$}
\nomenclature[G]{$\mathcal{R}$}{residual}
\nomenclature[G]{$\bm{n}$}{outward normal vector}
\nomenclature[G]{$\bm{S_f}$}{vector area}
\nomenclature[G]{$\bm{\varphi_i}$}{i-th POD basis function for velocity}
\nomenclature[G]{$\bm{\eta_i}$}{i-th POD basis function for supremizers}
\nomenclature[G]{${\chi_i}$}{i-th POD basis function for pressure}
\nomenclature[G]{$L_u$}{reduced basis space for velocity}
\nomenclature[G]{$L_p$}{reduced basis space for pressure}
\nomenclature[G]{$L_s$}{reduced basis space for supremizers}
\nomenclature[G]{$\bm{M}$}{full-order model mass matrix}
\nomenclature[G]{$\bm{C(u)}$}{full-order model convection matrix}
\nomenclature[G]{$\bm{A}$}{full-order model diffusion matrix}
\nomenclature[G]{$\bm{M_r}$}{ROM mass matrix} 
\nomenclature[G]{$\bm{C_r(u)}$}{ROM convection matrix}
\nomenclature[G]{$\bm{A_r}$}{ROM diffusion matrix}
\nomenclature[G]{${N_k}$}{number of parameters in the training set $\mathcal{K}$}
\nomenclature[G]{$\bm{{\mathcal{S}_u}}$}{snapshots matrix for the velocity field}
\nomenclature[G]{$\bm{{\mathcal{S}_p}}$}{snapshots matrix for the pressure field}
\nomenclature[G]{$\bm{{\mathcal{S}_s}}$}{snapshots matrix for supremizers}
\nomenclature[G]{$\otimes$}{tensor product}
\nomenclature[G]{$\bm{\nabla}\cdot$}{divergence operator}
\nomenclature[G]{$\bm{\nabla}\times$}{curl operator}
\nomenclature[G]{$\bm{\nabla}$}{gradient operator}
\nomenclature[G]{$\bm{\nabla}^s$}{symmetric gradient operator}
\nomenclature[G]{$\bm{a}$}{reduced vector of unknowns for velocity}
\nomenclature[G]{$\bm{b}$}{reduced vector of unknowns for pressure}
\nomenclature[G]{$\Delta$}{laplacian operator}
\nomenclature[G]{$\Box_e$}{value of a variable defined at the centre of a cell}
\nomenclature[G]{$\Box_f$}{value of a variable defined at the centre of a face}
\nomenclature[G]{$\beta$}{inf-sup stability constant}
\nomenclature[G]{$\left\lVert \cdot\right\rVert$}{norm in $L^2(\Omega)$}
\nomenclature[G]{$\langle \cdot , \cdot \rangle$}{inner product in $L^2(\Omega)$}


\newcommand{\dif}{\mbox{d}}

\title[]{Finite volume POD-Galerkin stabilised reduced order methods for the parametrised incompressible Navier-Stokes equations}

\author{Giovanni Stabile\textsuperscript{1,*}}
\thanks{\textsuperscript{*}Corresponding Author.}
\address{\textsuperscript{1}SISSA, International School for Advanced Studies, Mathematics Area, mathLab Trieste, Italy.}
\email{gstabile@sissa.it}

\author{Gianluigi Rozza\textsuperscript{1}}
\email{grozza@sissa.it}
\subjclass[2010]{78M34, 97N40, 35Q35}

\keywords{proper orthogonal decomposition; finite volume approximation; Poisson equation for pressure; inf-sup approximation; supremizer velocity space enrichment; Navier-Stokes equations.}

\date{}

\dedicatory{}


\begin{abstract} 
In this work a stabilised and reduced Galerkin projection of the incompressible unsteady Navier-Stokes equations for moderate Reynolds number is presented. The full-order model, on which the Galerkin projection is applied, is based on a finite volumes approximation. The reduced basis spaces are constructed with a POD approach. Two different pressure stabilisation strategies are proposed and compared: the former one is based on the supremizer enrichment of the velocity space, and the latter one is based on a pressure Poisson equation approach. 
\end{abstract}

\maketitle

\begin{multicols}{2}

%
%
%
%
\section{Introduction}\label{sec:intro}


During the last decades several progresses have been done in the field of computational fluid dynamics and more in general into the resolution of problems governed by partial differential equations. Nowadays one can find a wide variety of methods and computational libraries for the resolution of computational fluid dynamic problems. However, there are still many cases where the resolution of the governing equations, using standard discretisation techniques (Finite Element Method, Finite Volume Method and Finite Difference Method), become unfeasible. Such situations occur, for example, when a large number of different system configurations are in need of being tested (uncertainty quantification, optimization, ...) or a limited computational cost is required (real-time control). A possible way to overcome this limitation is the use of reduced order modelling (ROM) techniques \cite{hesthaven2015certified,quarteroniRB2016,ChinestaEnc2017,bennerParSys}. 

This technique is based on the assumption that the evolution in time of the dynamics of the system and its response into the parameter space (physical or geometrical) is governed by a reduced number of dominant modes.

In this work, a reduced basis POD-Galerkin method is considered and the interest is posed on parametrized time-dependent partial differential equations that govern fluid dynamics problems. 

In particular, the attention is devoted on RB-ROM generated starting from high dimensional finite volume approximations. This approximation method, that is particularly widespread for the resolution of fluid dynamic problems in many engineering fields (aeronautics engineering, naval engineering, automotive engineering, civil engineering...), is not particularly exploited in the field of Reduced Order Models. 

The FEM methodology is in fact more widespread for the generation of the RB spaces. RB methods starting from FEM full order approximations have been used to treat several problems based on linear elliptic equations \cite{Rozza2008}, linear parabolic equations \cite{grepl2005} and even non-linear problems \cite{Veroy2003,Grepl2007}. Few research works can be found dealing with finite volume schemes \cite{Haasdonk2008,Lorenzi2016,Drohmann2012,Stabile2017,HOR08}. 

First attempts to apply the reduced basis method in the context of viscous flows and Navier-Stokes equation can be found in \cite{ITO1998403,peterson1989}. 
In these pioneering works the projection is performed on divergence-free spaces without considering the pressure term at reduced order level.

It is well known that ROMs techniques for the Stokes and Navier-Stokes equations, obtained with Galerkin projection methods, are prone to several instability problems. In particular two different kind of instabilities have been observed and treated in literature: instabilities of the resulting system of ODEs for what concerns transient problems \cite{Iollo2000,Akhtar2009,Bergmann2009516,Sirisup2005218,taddei2017}, inf-sup pressure instabilities due to spurious pressure modes when the equivalent inf-sup condition for the reduced system is not fulfilled \cite{Caiazzo2014598,Ballarin2015,Gerner2012,Rozza2007,rozzanumerische}. 

The methods proposed in this work aim to deal with the second type of instabilities. The first aim of this work is to investigate and compare different two strategies for pressure stabilisation in the context of POD-Galerkin ROMs obtained from full-order finite volume approximations. The first proposed method is based on the supremizer enrichment of the velocity space to fulfil a reduced and parametrized version of the inf-sup condition, while the second proposed method is based onto the exploitation of a pressure Poisson equation during the online stage. An objective of this work is also to test the efficiency of the two methods for long-time integration. In the numerical examples, in fact, the attention is also paid to the performances of both stabilisation methods to approximate systems with periodic response, under long-time integration conditions. 

To the best of the authors' knowledge, the supremizer stabilisation technique \cite{Rozza2007} is here introduced for the first time in the context of a finite volume approximation. The work is organized as it follows: in \autoref{sec:math_form} we introduce the formulation and the methods used for the full-order approximation of the equations, in \autoref{sec:ROM} the two reduced order methodologies, object of this manuscript, are introduced and discussed in details. In \autoref{sec:num_exp} the two proposed ROMs techniques are tested on two different numerical benchmarks, dealing with the lid driven cavity problem, and the problem of the flow around a circular cylinder for moderate Reynolds numbers. Finally in \autoref{sec:conclusions} conclusions and perspectives are drawn, highlighting the directives for future improvements and developments.  
%
%
%
%
\section{Mathematical formulation and full-order approximation of the Navier-Stokes Equation}\label{sec:math_form}
The mathematical problem on which this work is focused is given by the unsteady incompressible parametrized Navier-Stokes equations. Considering an Eulerian frame on a space-time domain $Q = \Omega \times [0,T] \subset \mathbb{R}^d\times\mathbb{R}^+$ with $d=2,3$ the problem consists in finding the vectorial velocity field $\bm{u}:Q \to \mathbb{R}^d$ and the scalar pressure field $p:Q \to \mathbb{R}$ such that:
\begin{equation}\label{eq:navstokes}
\begin{cases}
\bm{u_t}+ \bm{\nabla} \cdot (\bm{u} \otimes \bm{u})- \bm{\nabla} \cdot 2 \nu \bm{\nabla^s} \bm{u}=-\bm{\nabla}p &\mbox{ in } Q,\\
\bm{\nabla} \cdot \bm{u}=\bm{0} &\mbox{ in } Q,\\
\bm{u} (t,x) = \bm{f}(\bm{x}) &\mbox{ on } \Gamma_{In} \times [0,T],\\
\bm{u} (t,x) = \bm{0} &\mbox{ on } \Gamma_{0} \times [0,T],\\ 
(\nu(\mu)\nabla \bm{u} - p\bm{I})\bm{n} = \bm{0} &\mbox{ on } \Gamma_{Out} \times [0,T],\\ 
\bm{u}(0,\bm{x})=\bm{k}(\bm{x}) &\mbox{ in } T_0,\\            
\end{cases}
\end{equation}
where $\Gamma = \Gamma_{In} \cup \Gamma_{0} \cup \Gamma_{Out}$ is the boundary of $\Omega$ and, is composed by three different parts $\Gamma_{In}$, $\Gamma_{Out}$ and $\Gamma_0$ that indicates, respectively, inlet boundary, outlet boundary and physical walls. The function $\bm{f}(\bm{x})$ represents the boundary conditions for the non-homogeneous boundary and $\bm{k}(\bm{x})$ denotes the initial condition for the velocity at $t=0$. It is also supposed that the boundary condition $\bm{f}$ is not depending on time. The parameter dependency is given by the kinematic viscosity $\nu(\mu)$ whose values are function of a parameter $\mu\in \mathcal{P}$ with $\mathcal{P}$ denoting the parameter space. It is moreover assumed that the kinematic viscosity is constant in the spacial domain. For sake of brevity, the parameter dependency of $\nu$ will be omitted in the formulations. Here, the equations are presented in its general form for an inlet-outlet problem, in \autoref{sec:num_exp}, where the numerical experiments are presented, it will be better specified the particular boundary conditions.
%
\subsection{The finite Volume Approximation}\label{subsec:FV_approx}
The system of equations in \eqref{eq:navstokes} together with its boundary and initial conditions is approximated at full-order level using a finite volume method. Here the finite volume approximation is briefly recalled, for more details the reader may see \cite{ferziger99:CMFD,ohlbergerFV}. Even though the finite volume approximation is normally derived starting directly from the integrated form of the governing equations here, in order to be consistent with the reduced basis methodology introduced in the next sections, the finite volume discretisation is presented as the restriction of the solution space associated with the weak formulation of the governing equations. The problem associated with the weak formulation of the Navier-Stokes equations consists in finding $(\bm{u},p) \in \mathcal{V}\times\mathcal{Q}$ such that:
%
\begin{equation}\label{eq:weak_NS}
\begin{split}
\mathcal{R}(\bm{u},p;\bm{v},q) = &\int_{\Omega} \frac{\partial \bm{u}}{\partial t} \cdot \bm{v} \dif \Omega + \int_{\Omega} \bm{\nabla} \cdot (\bm{u}\otimes\bm{u}) \cdot \bm{v} \dif \Omega \\
 & -\int_{\Omega} \bm{\nabla} \cdot 2 \nu \bm{\nabla^s}  \bm{u}  \cdot \bm{v} \dif \Omega+ \int_{\Omega} \bm{\nabla}p \cdot \bm{v} \dif \Omega \\
 &+\int_{\Omega} \bm{\nabla}\cdot \bm{u} q \dif \Omega = 0 \mbox{ \hspace{1cm} } \forall (\bm{v},q) \in \mathcal{V}\times\mathcal{Q},
\end{split}
\end{equation}
where $\mathcal{V} = \mathcal{H}^1(\Omega)$ and $\mathcal{Q}=\mathcal{L}^2(\Omega)$ are function spaces for velocity and pressure, respectively, and $\mathcal{R}$ is the residual associated with the weak formulation. The domain $\Omega$ is then divided into a tessellation $\mathcal{T} = \{ \Omega_e \}_{e=1}^{N_h}$ composed by a set of convex and non overlapping polygonals (finite volumes) such that $\Omega = \bigcup_{e=1}^{N_{FV}} \Omega_e$ and $\Omega_i \bigcap \Omega_j = \emptyset \mbox{ for } i\neq j$. The solution is then restricted to the finite dimensional space $\mathcal{V}^h\times\mathcal{Q}^h$ given by the space of the finite volume functions that are piecewise constant functions over each element $\Omega_e$. Note that while in the finite element method the solution space $\mathcal{V}^h\times\mathcal{Q}^h$ is given by suitable continuous piecewise polynomial functions, for the finite volume case the trial functions belong to the discontinuous space given by the finite volume functions. The solution is, in fact, sought into the finite dimensional space:
\begin{equation}
\mathcal{V}^h \times \mathcal{Q}^h = \spn\{ \mathcal{I}_k(\bm{x})\} , \forall k \in \mathcal{T}(\Omega),
\end{equation}
where $\mathcal{I}_k(\bm{x})$ is the basis function of each finite volume:
\begin{equation}
\begin{split}
& \mathcal{I}_k(\bm{x}) : \Omega \to \mathbb{R},\\
& \mathcal{I}_k(\bm{x}) = \left \{ \begin{array}{l} 1 \mbox{ if } \bm{x} \in \Omega_k\\ 0 \mbox{ if } \bm{x} \in \Omega \setminus \Omega_k \end{array} \right . .
\end{split}
\end{equation}
The problem consists then in finding $(\bm{u^h},p^h) \in \mathcal{V}^h\times \mathcal{Q}^h$ such that:
\begin{equation}
\mathcal{R}(\bm{u}^h,p^h;\bm{v}^h,q^h) = 0 \mbox{ \hspace{1cm} } \forall (\bm{v^h},p^h) \in \mathcal{V}^h \times \mathcal{Q}^h.
\end{equation}
Within a finite volume discretisation, all the divergence terms are rewritten in term of fluxes over the boundaries of each finite volume, making use of the Gauss's theorem:
\begin{equation}\label{eq:NS_disc}
\begin{aligned}
&\sum_{e=1}^{N_{FV}} \bm{v}_e\left( \int_{\Omega_e} \bm{u_t} \dif \Omega + \int_{\partial \Omega_e} \bm{n} \cdot (\bm{u}\otimes\bm{u})  \dif \Gamma \right. \\
&- \left. \int_{\partial \Omega_e} \bm{n} \cdot 2 \nu \bm{\nabla^s}\bm{u} \dif \Gamma + \int_{\partial \Omega_e}\bm{n} p  \dif \Gamma \right) \\
&+\sum_{e=1}^{N_{FV}} q_e \left( \int_{\partial\Omega_e} \bm{n} \cdot \bm{u} \dif \Gamma\right).
\end{aligned}
\end{equation}
Equation \eqref{eq:NS_disc} represents the semi-discretised version of the momentum and mass conservation. For sake of completeness the methodologies used to approximate the differential operators within a finite volume approximation are briefly recalled. It is in fact important to recall the procedures used to approximate each term inside equation \eqref{eq:NS_disc}. The same approaches are in fact used also during the generation of the reduced order model when the governing equations are projected onto the reduced basis spaces. In the following expressions, the value of the variables at the centre of the cells are indicated with the subscript $\Box_e$ while values at the centre of the faces are indicated with $\Box_f$. The term $\bm{S_f}$ indicates the surface area vector.  
The acceleration term is discretised, for the moment without specifying the method used to compute the time derivative, as:
\begin{equation}
\int_{\Omega_e} \bm{u_t} \dif \Omega = \bm{u_{t_e}} V_e,
\end{equation}
where $V_e$ denotes the volume of each cell. All the coefficients that multiply the acceleration terms can be recast in matrix form giving raise to the matrix $\bm{M}$ of equation~\eqref{eq:NS_matrix_form}. The non-linear convective term is discretised as:
\begin{equation}\label{eq:convect_discr}
\int_{\partial \Omega_e} \bm{n} \cdot (\bm{u}\otimes\bm{u}) \dif \Gamma = \sum_f \bm{S_f} \cdot \bm{u}_f \otimes \bm{u}_f,
\end{equation}
where $\bm{u}_f$ indicates the velocity at the centre of the faces. In this work the non-linear term is linearised with the substitution, of one of the $\bm{u}_f$ terms inside equation \eqref{eq:convect_discr}, with a previously calculated velocity that satisfies the continuity equation, for more details about this issue we refer to \cite{jasak1996error}. 
This discretisation process produces the matrix $\bm{C}$ of equation~\eqref{eq:NS_matrix_form}.
The diffusive term is discretised as:
\begin{equation}
\int_{\partial \Omega_e} \bm{n} \cdot 2 \nu \bm{\nabla^s}\bm{u} \dif \Gamma = \int_{\partial \Omega_e} \bm{n} \cdot \nu \bm{\nabla}\bm{u} \dif \Gamma = \nu \sum_f \bm{S_f} \cdot (\bm{\nabla u})_f,
\end{equation}
where the first equality follows from the incompressibility constraint and the term  $(\bm{\nabla u})_f$ indicates the gradient of the velocity field at the centre of each face. This is calculated, starting from the values at the centre of the neighbouring cells, using a finite difference scheme that includes a correction in the case of non-orthogonal meshes. For more details on this aspect the we refer to \cite{jasak1996error}. The coefficients obtained with such discretisation are used to assemble the matrix $\bm{A}$ of equation~\eqref{eq:NS_matrix_form}.
The term originated from the gradient of pressure, which gives raise to the matrix $\bm{B}$  of equation~\eqref{eq:NS_matrix_form}, is discretised as:
\begin{equation}
\int_{\partial \Omega_e} \bm{n} p  \dif \Gamma = \sum_f \bm{S_f} p_f,
\end{equation}
while the term originated from the divergence of velocity is discretised as:
\begin{equation}\label{eq:NS_con}
\int_{\partial \Omega} \bm{n} \cdot \bm{u} \dif \Gamma = \sum_{f=1}^{N_f} \bm{S_f} \cdot \bm{u_f}. 
\end{equation}
The coefficients of the above discretisation are used to assemble the matrix $\bm{P}$ of equation~\eqref{eq:NS_matrix_form}.
In the equations above, the values $\bm{u}_f$ and $p_f$, which are the values of the unknowns at the centre of each face, must be rewritten, using appropriate interpolation schemes, as functions of their values at the centre of the cells.
Even though a linear interpolation is appropriate for most of the above terms, in order to obtain an overall stable and accurate procedure the non-linear convective term needs particular attention and several schemes have been developed such as the upwind, second order linear upwind or MUSCL \cite{VANLEER1979101}. However, since it is not the objective of this manuscript to discuss the different types of stabilisation techniques for convection dominated problems, for more details we refer to \cite{ferziger99:CMFD}. 
For each finite volume, the interpolation coefficients obtained during the discretisation process are used to form an algebraic system of equations that can be rearranged in matrix form as:
\begin{equation}\label{eq:NS_matrix_form}
\begin{split}
\bm{M}\bm{\dot{u}}+\bm{C}(\bm{u})\bm{u}+\nu\bm{A}\bm{u}+\bm{B}\bm{p} & = \bm{0} \\
\bm{P}\bm{u} & = 0 ,
\end{split}
\end{equation}
The above system of equations can be solved using both a monolithic and partitioned approach, in the present case, at full-order level, a partitioned approach is preferred. In particular a PIMPLE algorithm is used, it consists into the combination of a SIMPLE \cite{Patankar1972} and PISO \cite{Issa1986} procedure. More details regarding the particular numerical schemes employed in the numerical experiments are reported in \autoref{sec:num_exp}.
The full-order simulations have been performed using the open source C\texttt{++} finite volume library OpenFOAM 5.0 \cite{OF} while the reduced order modelling computations are carried out using ITHACA-FV an in-house C\texttt{++} library. 

\section{Reduced order model with a POD-Galerkin method}\label{sec:ROM}
The full-order model illustrated in \autoref{sec:math_form} is solved for each $\mu^k \in \mathcal{K}=\{ \mu^1, \dots, \mu^{N_k}\} \subset \mathcal{P}$ where $\mathcal{K}$ is a finite dimensional training set of parameters chosen inside the parameter space $\mathcal{P}$. The considered problem can be simultaneously parameter and time dependent so, in order to collect snapshots for the generation of the reduced basis spaces one needs to consider both the time and parameter dependency. For this reason also discrete time instants $t^k \in \{t^1,\dots,t^{N_t}\} \subset [0,T]$ belonging to a finite dimensional training set, which is a subset of the simulation time window are considered as parameters. The total number of snapshots is then equal to $N_s = N_r\cdot N_t$.  The snapshots matrices $\bm{\mathcal{S}_u}$ and $\bm{\mathcal{S}_p}$, for velocity and pressure respectively, are then given by $N_s$ full-order snapshots:
\begin{gather}
\bm{\mathcal{S}_u} = [\bm{u}(\mu^1,t^1),\dots,\bm{u}(\mu^{N_r},t^{N_t})] \in \mathbb{R}^{N_u^h\times N_s},\\
\bm{\mathcal{S}_p} = [p(\mu^1,t^1),\dots,p(\mu^{N_r},t^{N_t})] \in \mathbb{R}^{N_p^h\times N_s}.
\end{gather}
The reduced order problem can be efficiently solved for all the set of parameters and time instants. In order to generate the reduced basis spaces, for the projection of the governing equations, one can find in literature several techniques such as the Proper Orthogonal Decomposition (POD), the Proper Generalized Decomposition (PGD) and the Reduced Basis (RB) with a greedy sampling strategy. For more details about the different strategies the reader may see \cite{Rozza2008,ChinestaEnc2017,Kalashnikova_ROMcomprohtua,quarteroniRB2016,Chinesta2011,Dumon20111387}. In this work a POD strategy is exploited and is chosen to apply the POD onto the full snapshots matrices that include both the time and parameter dependency. In case of parametric and time dependent problems also other approaches are available such use the POD-Greedy approach \cite{Haasdonk2008} or the nested POD approach where the POD is applied before in the time domain and later on the parameter space.
Given a general scalar or vectorial function $\bm{u}(t):Q \to \mathbb{R}^d$, with a certain number of realizations $\bm{u}_1,\dots, \bm{u}_{N_s}$, the POD problem consists in finding, for each value of the dimension of POD space $N_{POD} = 1,\dots,N_s$, the scalar coefficients $a_1^1,\dots,a_1^{N_s},\dots,a_{N_s}^1,\dots,a_{N_s}^{N_s}$ and functions $\bm{\varphi}_1,\dots,\bm{\varphi}_{N_s}$ that minimize the quantity:
\begin{gather}\label{eq:pod_energy}
E_{N_{POD}} = \sum_{i=1}^{N_s}||\bm{u}_i-\sum_{k=1}^{N_{POD}}a_i^k \bm{\varphi_k}|| \hspace{0.5cm}\forall\mbox{ } N_{POD} = 1,\dots,N\\
\mbox{ with } \ltwonorm{\bm{\varphi}_i,\bm{\varphi}_j} = \delta_{ij} \mbox{\hspace{0.5cm}} \forall \mbox{ } i,j = 1,\dots,N_s .
\end{gather}
In this case the velocity field $\bm{u}(t)$ is used as example. It can be shown \cite{Kunisch2002492} that the minimisation problem of Equation~\eqref{eq:pod_energy} is equivalent of solving the following eigenvalue problem:
\begin{gather}
\bm{\mathcal{C}^u}\bm{Q}^u = \bm{Q^u}\bm{\lambda^u} ,\\
\mathcal{C}^u_{ij} = \ltwonorm{\bm{u}_i,\bm{u}_j} \mbox{\hspace{0.5cm} for } i,j = 1,\dots,N_s ,
\end{gather}
where $\bm{\mathcal{C}^u}$ is the correlation matrix obtained starting from the snapshots $\bm{\mathcal{S}_u}$, $\bm{Q^u}$ is a square matrix of eigenvectors and $\bm{\lambda^u}$ is a vector of eigenvalues. 
\begin{remark}
Normally in the standard finite element framework, since for velocity the natural functional space belongs to $\mathcal{H}^1(\Omega)$, to compute its correlation matrix $\bm{\mathcal{C}^u}$ the $H^1$ norm is preferred and the $L^2$ norm is used to compute the correlation matrix of pressure. Here, for both velocity and pressure, the $L^2$ norm is preferred because as illustrated in section \ref{sec:math_form}, using a finite volume method, both the velocity and the pressure belong to discontinuous spaces and, in order to compute the gradient necessary for the $H^1$ norm evaluation, one would introduce further discretisation error. Moreover the $L^2$ norm has a direct physical meaning being directly correlated with the kinetic energy of the system. 
\end{remark}
The basis functions can then be obtained with: 
\begin{equation}
\bm{\varphi_i} = \frac{1}{N_s\lambda_i^u}\sum_{j=1}^{N_s} \bm{u}_j Q^u_{ij}.
\end{equation}
The POD spaces are constructed for both velocity and pressure using the aforementioned methodology resulting in the spaces:
\begin{equation}
\begin{split}
&L_u = [\bm{\varphi_1}, \dots , \bm{\varphi_{N_u^r}}] \in \mathbb{R}^{N_{u}^h \times N_u^r},\\
&L_p = [{\chi_1}, \dots , {\chi_{N_p^r}}] \in \mathbb{R}^{N_{p}^h \times N_p^r}.
\end{split}
\end{equation}
where $N_u^r$, $N_p^r < N_s$ are chosen according to the eigenvalue decay of the vectors of eigenvalues $\bm{\lambda}^u$ and $\bm{\lambda}^p$. 

Once the POD functional spaces are set, the reduced velocity and pressure fields can be approximated with:
\begin{equation}\label{eq:aprox_fields}
\bm{u^r} \approx \sum_{i=1}^{N_u^r} a_i(t,\mu) \bm{\varphi_i}(\bm{x}), \mbox{\hspace{0.5 cm}}
p^r\approx \sum_{i=1}^{N_p^r} b_i(t,\mu)\chi_i(\bm{x}).
\end{equation}
Where the coefficients $a_i$ and $b_i$ depend only on the time and parameter spaces and the basis functions $\bm{\varphi}_i$ and $\bm{\chi}_i$ depend only on the physical space. The unknown vectors of coefficients $\bm{a}$ and $\bm{b}$ can be then obtained through a Galerkin projection of the governing equations onto the POD reduced basis spaces and with the resolution of following reduced algebraic system:
\begin{equation}\label{eq:alg_red}
\begin{split}
\bm{M_r} \bm{\dot{a}} -\nu\bm{A_r}\bm{a}+ \bm{C_r}(\bm{a})\bm{a}+\bm{B_r}\bm{b} = \bm{0} \\
\bm{P_r} \bm{a} = 0,
\end{split}
\end{equation} 
where the terms inside equation~\eqref{eq:alg_red} are evaluated with:
\begin{equation}
\begin{split}\label{eq:red_matrices}
& M_{r_{ij}} = \ltwonorm{\bm{\varphi_i},\bm{\varphi_j} } \mbox{ , } A_{r_{ij}} = \ltwonorm{\bm{\varphi_i}, \bm{\nabla} \cdot 2\bm{\nabla^s\varphi_j}} \mbox{ , } \\
& B_{r_{ij}} = \ltwonorm{\bm{\varphi_i},\bm{\nabla}{\chi_j}} \mbox{ , } P_{r_{ij}} = \ltwonorm{\chi_i,\bm{\nabla} \cdot \bm{\varphi_{j}}}.
\end{split}
\end{equation}
Once the reduced basis spaces $L_u$ and $L_p$ are defined through the basis functions $\bm{\varphi_i}$ and $\chi_i$, all the reduced matrices of equation \eqref{eq:alg_red} can be precomputed during an offline stage without difficulties with the exception of the reduced matrix $\bm{C_r}(\bm{a})$, which is originated by the non-linear convective term. The strategy employed here consists into the storage of a third-order tensor $\bm{\mathsf{C}}_r$ \cite{quarteroni2007numerical,Rozza2009} whose entries are given by: 
\begin{equation}\label{eq:third_order}
\mathsf{C}_{r_{ijk}} = \ltwonorm{\bm{\varphi_i}, \bm{\nabla} \cdot (\bm{\varphi_j} \otimes \bm{\varphi_k})}. 
\end{equation}
During the online stage, at each fixed point iteration of the solution procedure, each entry of the contribution to the reduced residual given by the convective term $\mathbfcal{R}^r_c =\bm{C_r}(\bm{a})\bm{a}$, can be computed with:
\begin{equation}\label{eq:third_order_res}
\mathbfcal{R}^r_{c_i} = (\bm{C_r}(\bm{a})\bm{a})_i = \bm{a}^T \bm{\mathsf{C}}_{r_{i \bullet \bullet}} \bm{a}.
\end{equation}
\begin{remark}
The dimension of the $\bm{\mathsf{C}_r}$ tensor is increasing with the cube of the number of basis functions. For this reason, when a large number of basis functions are employed this approach may lead to high storage costs. In the present case a relatively small number of basis functions is considered ($N<20$) but in case of richer reduced spaces other approaches, such as EIM-DEIM \cite{Xiao20141,BARRAULT2004667} or Gappy-POD \cite{Carlberg2013623} could become more affordable. 
\end{remark} 
\subsection{Initial conditions}
The initial conditions for the ROM system of ODEs of equation \eqref{eq:alg_red} are obtained performing a Galerkin projection of the initial full-order condition $\bm{u}(0)$ onto the POD basis spaces. As it will be shown in the next sections, also functional spaces with non-orthogonal basis functions $\bm{\varphi_i}$ are considered. For this reason, the initial coefficients $\bm{a_0}$ have to be obtained solving the following linear system of equations:
\begin{equation}
\bm{M_r} \bm{a_0} = \bm{e},
\end{equation}
where $\bm{M_r}$ is obtained following the expression of equation \eqref{eq:red_matrices} and the components of the $\bm{e}$ vector are obtained with $e_i = \ltwonorm{\bm{\varphi_i},\bm{u}(0)}$.

\subsection{Stability Issues}\label{sec:stab_issues}
The reduced problem, as formulated in \autoref{sec:ROM}, presents stability issues. It is well known in fact that, using a mixed formulation for the approximation of the incompressible Navier-Stokes equations, the approximation spaces need to satisfy the inf-sup (Ladyzhenskaya-Brezzi-Babuska) condition \cite{BREZZI199027,boffi_mixed}. It is required that there should exist a constant $\beta > 0 $, independent to the discretisation parameter $h$, such that:
\begin{equation}
\inf_{q_h \in \mathcal{Q}} \sup_{\bm{v_h} \in \mathcal{V}} \frac{\langle \nabla \cdot \bm{v_h}, q_h \rangle}{\norm{\nabla \bm{v_h}}\norm{q_h}} \ge \beta > 0.
\end{equation}
Dealing with finite element methods, for what concerns the full-order level, this requirement can be met choosing appropriate finite element spaces such as the standard Taylor-Hood ($\mathbb{P}_2$ - $\mathbb{P}_1$). In this case, at full-order level, since a finite volume formulation is used, no attention is paid to this issue but, at reduced-order level, where a mixed formulation based on a projection method is used, one has to ensure that a \emph{reduced} version of the LBB condition is fulfilled; in fact, at reduced order level, where a Galerkin approach is exploited, two different spaces are used to approximate the velocity and the pressure variables. Regardless the full-order discretisation technique, even though the snapshots have been obtained by stable numerical methods, there is no guaranty that the original properties of the full-order system are preserved after the Galerkin projection onto the RB spaces \cite{Rozza2007,Gerner2012,Ballarin2015}. 
To overcome this issue, most of the contributions available in literature do not attempt to recover the pressure field and, at reduced order level, resolve only the momentum equation neglecting the contribution of the gradient of pressure. This choice is justified by the fact that, the projection of the pressure gradient onto the POD spaces is numerically zero for the case of enclosed flows as presented in \cite{deane1991,Ma2002,noack1994}, or in the case of inlet-outlet problems with outlet far from the obstacle \cite{Akhtar2009}. However, as highlighted in \cite{noack2005}, in many applications the pressure term is needed  and cannot be neglected. This work aims at comparing two different strategies for pressure stabilisation during the resolution of the reduced problem. 
In the first proposed approach the velocity space is enriched in order to satisfy a reduced version of the inf-sup condition \cite{Ballarin2015,Rozza2007}, this approach will be henceforth denoted as SUP-ROM. The second approach is based on a Leray-Helmholtz projection by exploiting at reduced order level a Poisson equation for pressure \cite{Stabile2017,Akhtar2009}; it will be henceforth denoted as PPE-ROM. The two methods proposed in this work are just two options among the possible choices to obtain stable ROMs for what concerns both velocity and pressure fields. It is worth mentioning also other possibilities that rely on pressure stabilised Petrov-Galerkin (PSPG) methods during the online procedure \cite{Caiazzo2014598,Baiges2014189}. In other approaches it is assumed that velocity and pressure share the same temporal coefficients and during the online procedure only the momentum equation is exploited \cite{Bergmann2009516,Lorenzi2016}.    
\begin{remark}
It is important to remark that the reduced order model is obtained with a projection method. For this reason, regardless from the approximation procedure used to produce the snapshots matrices for the generation of the snapshots, even though we are dealing with finite volume full-order approximations, which do not require the fulfilment of the inf-sup condition, at reduced order level this condition becomes relevant and needs to be met.
\end{remark}

\subsection{Supremizer enrichment}\label{subsec:sup_enrich}
The first proposed approach relies onto the fulfilment of a reduced and also parametric, in case, version of the inf-sup condition. As mentioned in section \ref{sec:stab_issues}, the problem, formulated using a mixed formulation, in order to be solvable and stable needs to meet the inf-sup condition. Within this approach, the velocity supremizer basis functions $L_s$ are computed and added to the reduced velocity space which is transformed into $\tilde{L}_u$:
\begin{equation}
\begin{split}
&L_s =[\bm{\eta_1}, \dots, \bm{\eta_{N_s^r}}] \in \mathbb{R}^{N_u^h \times N_s^r}, \mbox{ } \\
&\tilde{L}_u = [\bm{\varphi_1},\dots,\bm{\varphi_{N_u^r}},\bm{\eta_1},\dots,\bm{\eta_{N_s^r}}] \in \mathbb{R}^{N_u^h \times (N_u^r+N_s^r)}.
\end{split}
\end{equation}
We remark that, in this case, the space $\tilde{L}_u$ is not any-more formed by only orthogonal basis functions. The POD is in fact applied separately onto the velocity snapshots and onto the supremizer snapshots. These basis functions are chosen solving a supremizer problem which ensures that a reduced version of the inf-sup condition is fulfilled. The supremizer solution $\bm{s_i}$ is the element that, given a certain pressure basis function $p_i$, permits the realization of the inf-sup condition. For each pressure basis function the corresponding supremizer element can be found solving the following problem:
\begin{equation}\label{eq:sup_problem}
\begin{cases}
\Delta \bm{s_i} = - \bm{\nabla} p_i &\mbox{ in } \Omega,\\
\bm{s_i}=\bm{0} &\mbox{ on } \partial\Omega.
\end{cases}
\end{equation}
In this case the supremizer problem, which in a standard finite element setting is solved starting directly from the weak formulation, is expressed in strong form and solved using the full-order finite volume solver. For more details regarding the derivation one may see \cite{Rozza2007,Gerner2012,Ballarin2015}. As presented in \cite{Ballarin2015}, two different strategies can be employed to enrich the velocity space and select the supremizer space $L_s$ such that the inf-sup condition is met: an \emph{exact supremizer enrichment} procedure and an \emph{approximate supremizer enrichment} procedure. 

In the \emph{exact} approach, for each basis of the pressure space $\chi_i$, the problem of equation~\eqref{eq:sup_problem} is solved and the resulting solution is used as additional basis function for the velocity space. Using such an approach it can be proven that the resulting ROM that is obtained by the Galerkin projection onto the RB spaces is inf-sup stable \cite{Ballarin2015}.
In the approximated approach the problem is solved for each pressure snapshot $p(\mu,t)$ and a snapshots matrix of supremizer is assembled:
\begin{equation}\label{eq:sup_snap}
\bm{S_s} = [\bm{s}(\mu^1,t^1),\dots,\bm{s}(\mu^{N_r},t^{N_t})] \in \mathbb{R}^{N_u^h\times N_s}.
\end{equation}
A POD procedure is then applied to the resulting snapshots matrix in order to obtain the supremizer POD basis functions $\bm{\eta_i}$. This procedure permits to strongly reduce the online computational cost. The supremizer basis functions do not depend, in fact, on the particular pressure basis functions but are computed during the offline phase, starting directly from the pressure snapshots. However, with such an approach it is not possible to rigorously show that the inf-sup condition is satisfied and it is only possible to rely on heuristic criteria or to check it during a post-processing stage \cite{Ballarin2015}, such as a computational validation. In general this is true and reliable for non-geometric parametrization. 
\subsection{Pressure Poisson Equation}
The second approach is based on an alternative form of the Navier-Stokes equations where the incompressibility constraint $\nabla \cdot \bm{u} = 0$ is replaced by a Poisson equation for pressure. This alternative form of the Navier-Stokes equations was firstly proposed in the context of projection methods introduced by Chorin \cite{chorin1967} and Teman \cite{temam1968methode}. In these methods an intermediate velocity is first computed and later projected onto the space of divergence-free vector fields through the solution of a Poisson equation. For a thorough review on projection methods the reader may refer to \cite{GUERMOND20066011}. These methods, that can interpreted as a variant of pressure stabilisation methods \cite{Rannacher1992}, can be successfully applied also with functional spaces that do not satisfy the inf-sup condition \cite{Guermond1998,LIU20103428}. The idea of projection methods in the context of POD-Galerkin ROMs is attractive since the POD velocity modes, constructed from divergence-free snapshots, are indeed divergence-free (up to numerical precision). This approach was firstly proposed in \cite{Akhtar2009} and recently re-proposed in a finite volume setting in \cite{Stabile2017}.  
The modified set of equations considered here read:
\begin{equation}\label{eq:NS_PPE}
\begin{cases}
\bm{u_t}+ \bm{\nabla} \cdot (\bm{u} \otimes \bm{u})- \bm{\nabla} \cdot 2 \nu \bm{\nabla^s} \bm{u}=-\bm{\nabla}p &\mbox{in } Q\\
\Delta p = -\bm{\nabla} \cdot (\bm{\nabla} \cdot (\bm{u} \otimes \bm{u})) & \mbox{in } Q,\\
\bm{u}(t,\bm{x}) = \bm{0} & \mbox{on } \Gamma_0 \times [0,T],\\
\bm{u}(t,\bm{x}) = \bm{f}(\bm{x}) & \mbox{on } \Gamma_{In}, \\
\frac{\partial p}{\partial \bm{n}} = -\nu \bm{n} \cdot (\bm{\nabla} \times \bm{\nabla} \times \bm{u} ) - \bm{n} \cdot \bm{f}_t  & \mbox{on } \Gamma.
\end{cases}
\end{equation} 
In equation~\eqref{eq:NS_PPE} the Poisson equation for pressure is obtained taking the divergence of the momentum equation and exploiting the continuity constraint. The above formulation can be derived only under the assumption of sufficient smoothness of the solution $\bm{u}-p$, so that the divergence of the momentum equation makes sense. The last term of the above equation is a Neumann boundary condition for the pressure Poisson equation. This boundary condition is here introduced, to the best of authors knowledge, for the first time in the context of POD-Galerkin methods. 
In \cite{Caiazzo2014598,noack2005} an homogeneous Neumann boundary condition is prescribed while in \cite{giere2006} a different Neumann condition for pressure is employed. This boundary condition is derived starting from the enforcement of the divergence-free constraint on the boundary $\nabla \cdot \bm{u} = 0 |_{\Gamma}$, more details concerning the derivation can be found in \cite{Orszag1986,JOHNSTON2004221} where this condition is proposed in the context of a full order finite element formulation. Alternative ways to enforce a boundary condition for the pressure term are given in \cite{Gresho1987,LIU20103428,GUERMOND20066011}. The system of equations~\eqref{eq:NS_PPE} is used as starting point to derive a weak formulation and to construct the Galerkin system. The reduced system is obtained substituting the velocity and pressure field expansions of equation \eqref{eq:aprox_fields} and projecting the momentum and pressure equations onto the subspaces spanned by the velocity and pressure modes $\bm{\varphi_i}$ and $\chi_i$, respectively:
\begin{subequations}
\begin{align} \label{eq:PPE_weak(a)}
&\ltwonorm{\bm{\varphi_i}, \bm{u_t}+\bm{\nabla} \cdot (\bm{u} \otimes  \bm{u}) + \bm{\nabla} p - \bm{\nabla} \cdot 2 \nu \bm{\nabla^s} \bm{u}} = 0, \\
\begin{split}
&\ltwonorm{\nabla \chi_i,\nabla p} + \ltwonorm{\nabla{\chi_i},\bm{\nabla} \cdot (\bm{u} \otimes \bm{u})}
\\ \label{eq:PPE_weak(b)}
& - \nu \langle \bm{n} \times \bm{\nabla} \chi_i, \bm{\nabla} \times \bm{u} \rangle_{\Gamma} - \langle \chi_i, \bm{n} \cdot \bm{f_t} \rangle_{\Gamma} = 0. 
\end{split}
\end{align}
\end{subequations}
Where the equations \eqref{eq:PPE_weak(b)} has been obtained with integration by part of the laplacian term and exploiting the pressure boundary condition. In equation \eqref{eq:PPE_weak(b)} only first order derivatives appear, for this reason, during the Galerkin projection, the numerical error introduced by the numerical differentiation can be significantly reduced. We remark that, since we are using a finite volume formulation, there are no (theoretical) limitations regarding the achievable order of differentiation; however, as highlighted in \autoref{sec:math_form} derivatives are approximated by numerical methods and therefore an higher order of derivation would introduced a higher numerical error.
Performing a substitution of the velocity and pressure fields with the approximate expansion, it is possible to obtain the Galerkin system, which consists into a system of ODEs and reads:
\begin{subequations}
\begin{align}
&\bm{M_r}\bm{\dot{a}} - \nu\bm{A_r}\bm{a} + \bm{a}^T \bm{\mathsf{C_r}} \bm{a} + \bm{B_r} \bm{b} = 0, \label{eq:PPE_alg_mom} \\
&\bm{D_r}\bm{b}+\bm{a}^T\bm{\mathsf{G}_r}\bm{a} - \nu\bm{N_r}\bm{a} - \bm{F_r} = 0. \label{eq:PPE_alg_PPE}
\end{align}
\end{subequations}
Where the matrices and the tensor inside equation \eqref{eq:PPE_alg_mom} are obtained using the expressions given in \eqref{eq:red_matrices} and \eqref{eq:third_order}, while the matrices and the tensor inside \eqref{eq:PPE_alg_PPE} are given by:
\begin{equation}\label{eq:PPE_red_matrices}
\begin{split}
&D_{r_{ij}} = \ltwonorm{\bm{\nabla} \chi_i,\bm{\nabla} \chi_j} \mbox{ , }  \\ 
&\mathsf{G_r}_{ijk} = \ltwonorm{\bm{\nabla} \chi_i, \bm{\nabla} \cdot (\bm{\varphi_j} \otimes \bm{\varphi_k})},\\
& N_{r_{ij}} = \langle \bm{n} \times \bm{\nabla} \chi_i, \bm{\nabla} \times \bm{\varphi_j} \rangle_{\Gamma}, \\
& F_{r_i} = \langle \chi_i , \bm{n} \cdot \bm{f_t} \rangle_{\Gamma}.
\end{split}
\end{equation}
The residual associated with the non-linear term in the equation \eqref{eq:PPE_alg_PPE} is evaluated using the same strategy proposed in equation \eqref{eq:third_order_res}, i.e. storing the third order tensor $\bm{\mathsf{G}_r}$. In the numerical experiments considered in this work, the boundary condition are not varying in time, for this reason the term $\bm{f}_t$ is identically equal to zero and so is the reduced vector $\bm{F_r}$. The pressure boundary condition of equation~\eqref{eq:NS_PPE}, firstly proposed in \cite{JOHNSTON2004221}, to the best of the authors knowledge, is here introduced for the first time in the context of POD-Galerkin methods. These additional terms are neglected in \cite{Akhtar2009}, while a different boundary condition for pressure is considered in \cite{Caiazzo2014598}.  
Making a comparison with the SUP-ROM, it is possible to notice that the resulting ROM has an additional complexity due to the computation of the terms $\bm{N_r}$ and $\bm{F_r}$ but, due to the absence of the additional supremizer modes, produces a reduced dynamical which has a smaller dimension. It is worth mentioning, as highlighted in \cite{JOHNSTON2004221}, that the additional boundary condition for pressure of equation~\eqref{eq:NS_PPE} is not consistent in the case of steady flows. 
\section{Numerical Experiments}\label{sec:num_exp}
In this section the two different proposed stabilisation methods are tested and compared on two benchmark test cases. The first benchmark consists into the well known and studied lid driven cavity problem \cite{SCHREIBER1983310}. The second benchmark consists into the flow around a circular cylinder for moderate Reynolds number ($100<\mbox{Re}<200$) \cite{Schafer1996}. In the first case, any kind of parametrisation is introduced, while in the second case the kinematic viscosity is parametrised. 
\subsection{Lid driven cavity problem} 
\begin{figure*}
\begin{minipage}{0.48\textwidth}
\begin{minipage}{0.5\textwidth}
\includegraphics[width=\textwidth]{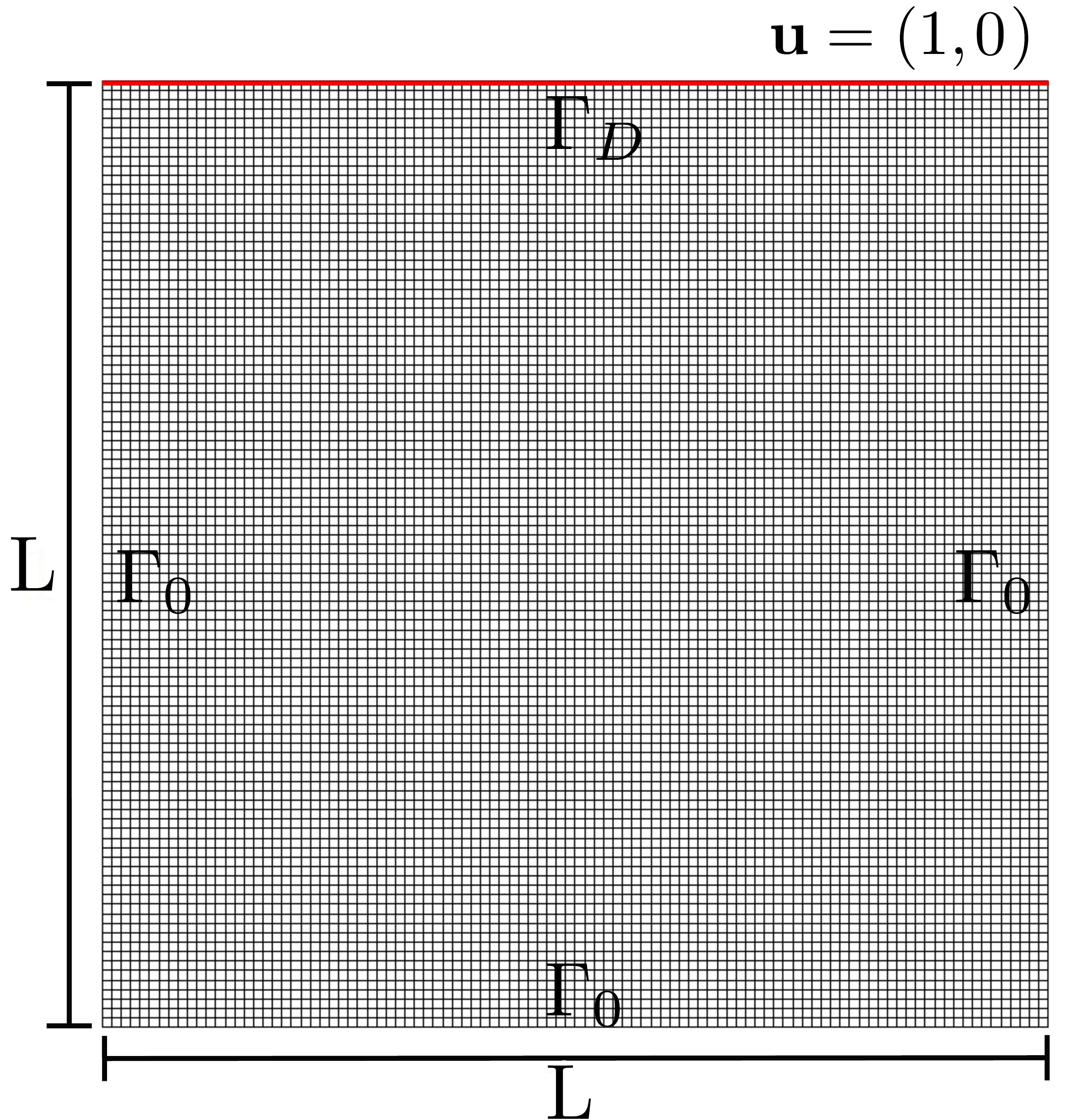}
\end{minipage} 
\begin{minipage}{0.45\textwidth}
\resizebox{\textwidth}{!}{%
\begin{tabular}{ l | c | c }
    & $\Gamma_{D}$ & $\Gamma_{0}$ \\ \hline 
    $\bm{u}$ & $\bm{u} = (1,0)$ & $\bm{u} = (0,0)$\\\hline 
    $p$ & $\nabla p \cdot \bm{n} = 0$ & $\nabla p \cdot \bm{n} = 0$
   \end{tabular}} 
\end{minipage}
\end{minipage}
\caption{Sketch of the mesh for the lid driven cavity problem together with the boundary subdivisions and boundary conditions.}\label{fig:mesh_cavity}
\end{figure*}

\begin{figure*}
\centering
\includegraphics[width=0.40\textwidth]{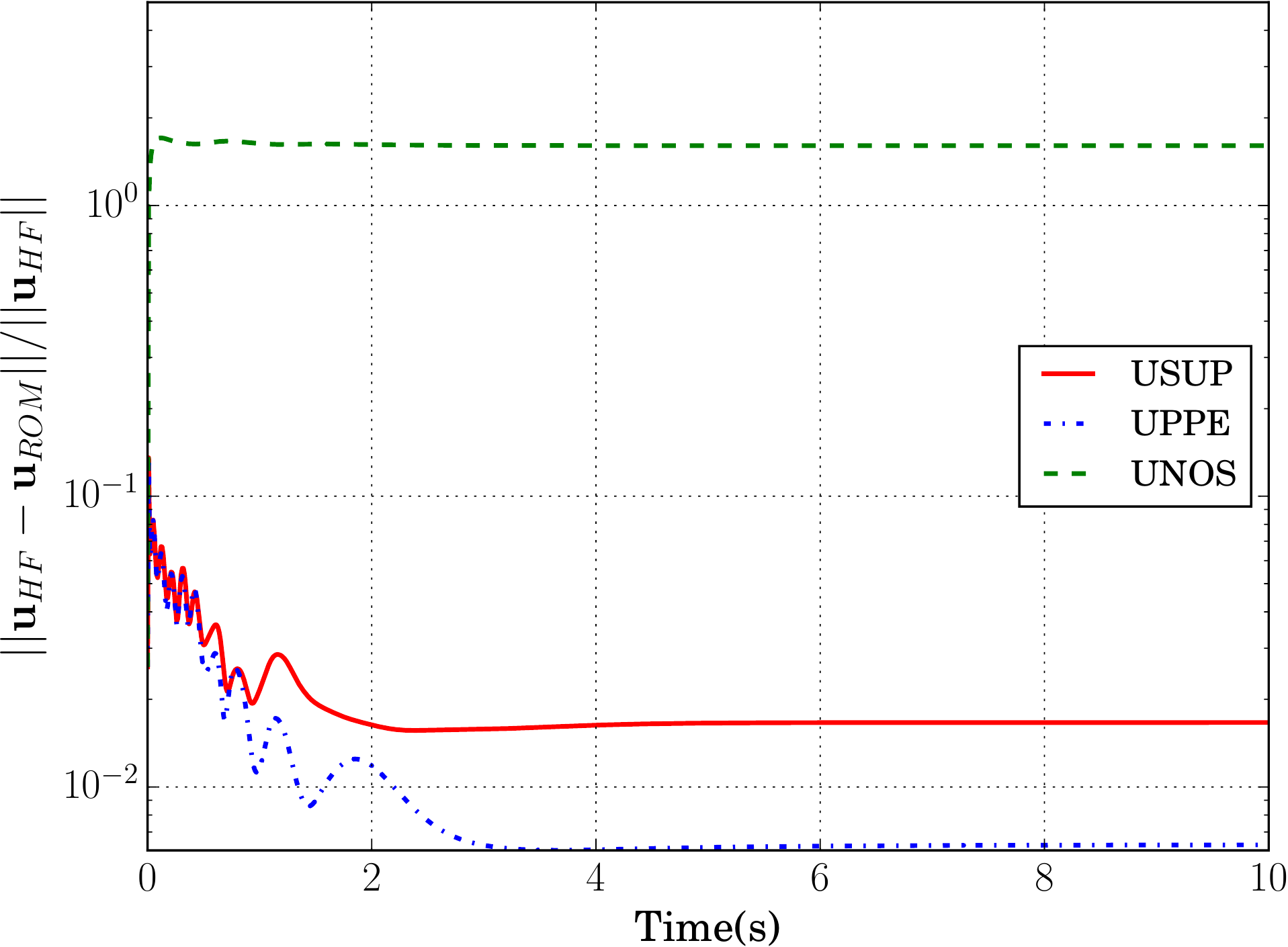}
\includegraphics[width=0.40\textwidth]{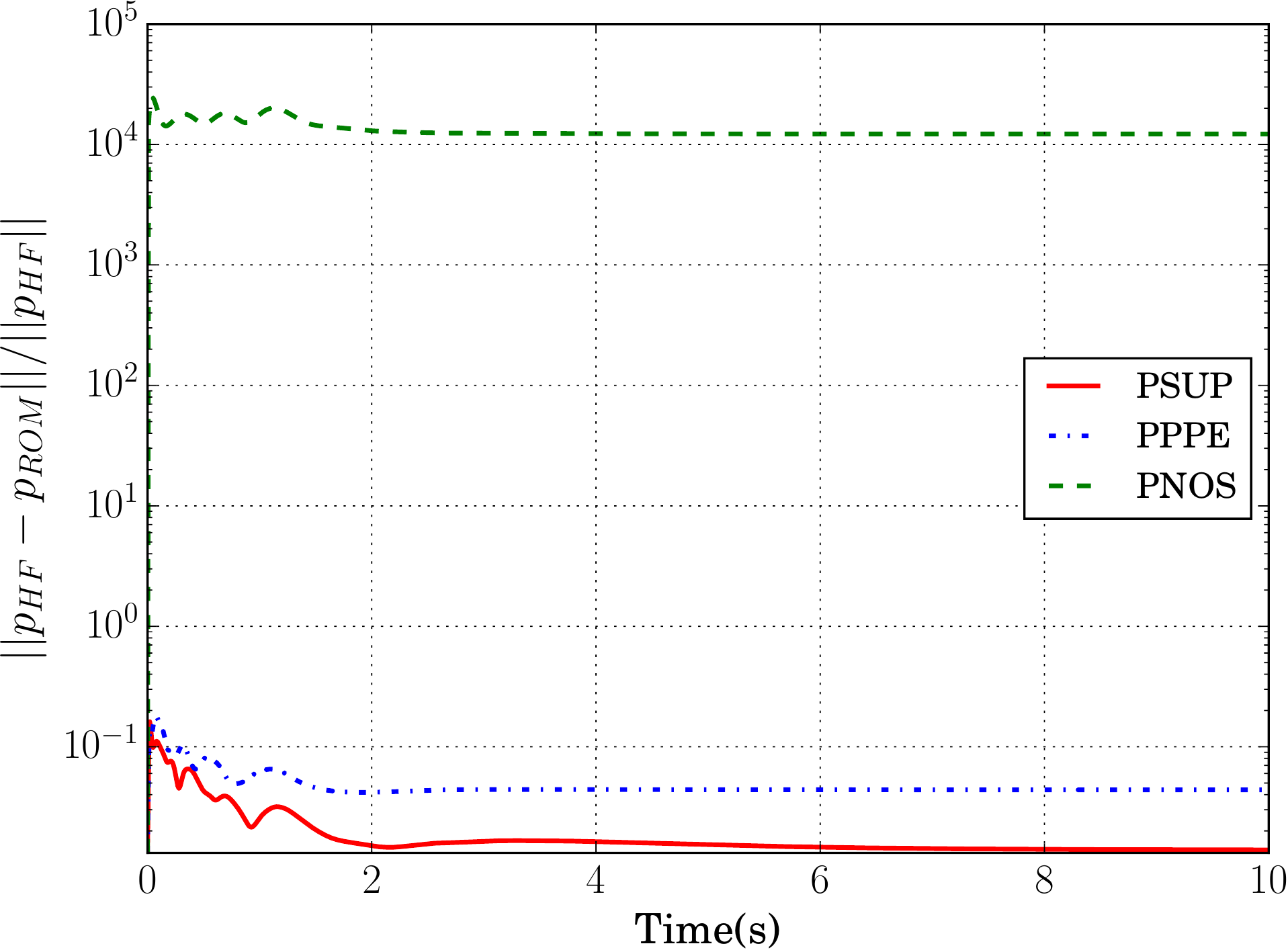}
\caption{Error analysis for the velocity field. The $L^2$ norm of the relative error is plotted over time for three different models: with supremizer stabilisation (USUP - continuous red line), pressure Poisson equation stabilisation (UPPE - dotted blue line), and without stabilisation (PNOS - dashed green line). The ROMs are obtained with 10 modes for velocity, pressure and supremizer.}\label{fig:error_cavity_u}
\end{figure*}

As said, the first proposed benchmark consists into the well known lid driven cavity problem. The simulation is carried on a two-dimensional square domain of length $L = 0.1 \mbox{m}$. The boundary is subdivided into two different parts $  \Gamma = \Gamma_D \cup \Gamma_0$ and the boundary conditions for velocity and pressure are set according to figure~\ref{fig:mesh_cavity}. At the top of the cavity a constant uniform and horizontal velocity equal to $u_x = 1 \si{m/s}$ is prescribed. The mesh is structured and counts $40000$ quadrilateral cells, $200$ on each dimension of the square. The kinematic viscosity is equal to $\nu = \num{1e-4} \si{m^2/s}$ that leads to a Reynolds number of $1000$. For what concerns the full-order simulation, the time discretisation is treated using a second order backward differencing scheme, while the  discretisation in space is performed with a forth order interpolation scheme. The time step is kept constant and equal to $\Delta t = \num{5e-4}$ and the simulation is run till $T = 10 \si{s}$. The snapshots are acquired every $0.01 \si{s}$ giving a total number of snapshots equal to $1000$. For what concerns the reduced order model, the dimension of the reduced spaces for velocity and pressure is set, for both the presented methodologies, equal to $N_u^r = 10$ and $N_p^r = 10$. In the SUP-ROM the reduced space for velocity is enriched with $10$ additional supremizer modes. This selection is done according to table~\ref{tab:cum_eig} where it is possible to observe that, such number of modes, is sufficient to retain more than $99.9 \%$ of the energy for both velocity and pressure. We remark that, in this numerical experiment, any kind of parametrisation is introduced. The ROM is in fact used to simulate the same conditions tested in the full-order setting and the results are compared against the full-order simulation results. The time discretisation, at reduced order level, is treated making use of a first order backward Newton method. 
Figure~\ref{fig:comparison_cavity} depicts a comparison between the HF simulation and the ROM one, for both velocity and pressure fields, at different time instants. As one can see from the figure, both models are capable of reproducing the main flow pattern for both the two fields. Figure~\ref{fig:error_cavity_u}  reports the evolution in time of the $L^2$ relative error for velocity and pressure respectively. The plots report also the error without any type of stabilisation. It is clear that, without stabilisation, even though the ROM is not diverging, both the velocity and pressure fields are completely unreliable. For this particular numerical test, the SUP-ROM produces, with respect to the PPE-ROM, worse results for what concerns the velocity field but better results for what concerns the pressure field. This difference can be justified by the fact that, within a supremizer stabilisation technique, the POD velocity space is enriched by non-necessary (for the correct reproduction of the velocity field) supremizer modes. During the initial transient, both fields present a higher relative error and this fact is due to the relatively low number of snapshots acquired during the initial transient. The snapshots, as highlighted above, are equally distributed in time and, to enhance the performance of the ROM one should concentrate the snapshots in the time span where the system exhibits the most non-linear behaviour. For what concerns the SUP-ROM, according to the indication reported in \cite{Ballarin2015}, the number of supremizer modes is chosen equal to the number of pressure modes. Table~\ref{tab:cum_eig} report also the value of the inf-sup constant $\beta$, obtained keeping constant the number of velocity and pressure modes (10 modes for velocity and 10 modes for pressure) and varying the number of supremizer modes. As one can observe from the table, by increasing the number of supremizer modes, leads to a remarkable increase of the inf-sup constant. 
\begin{figure*}
\centering
\begin{minipage}[c]{0.16\textwidth}
\centering
\includegraphics[width=\textwidth]{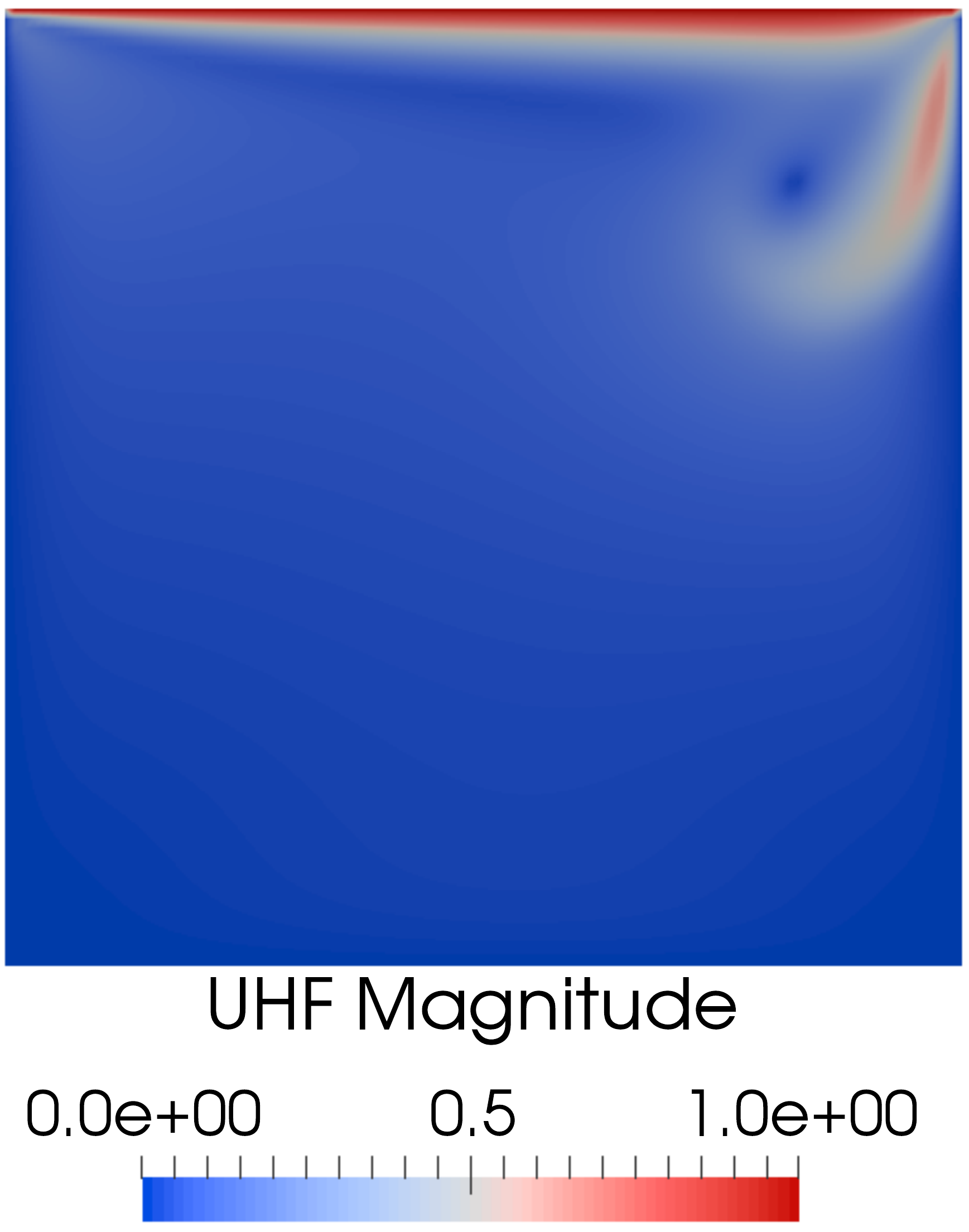}
\end{minipage}
\begin{minipage}[c]{0.16\textwidth}
\centering
\includegraphics[width=\textwidth]{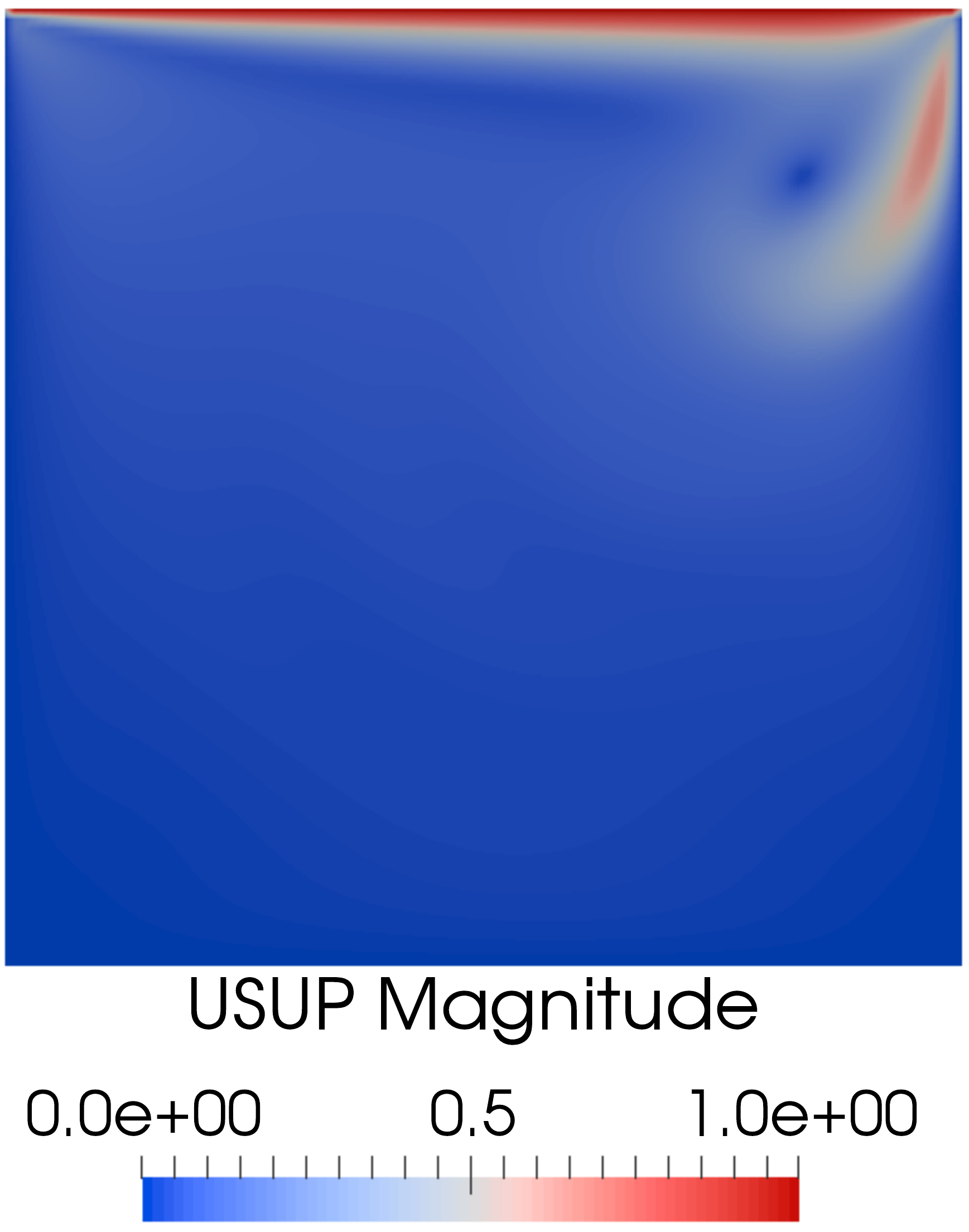}
\end{minipage}
\begin{minipage}[c]{0.16\textwidth}
\centering
\includegraphics[width=\textwidth]{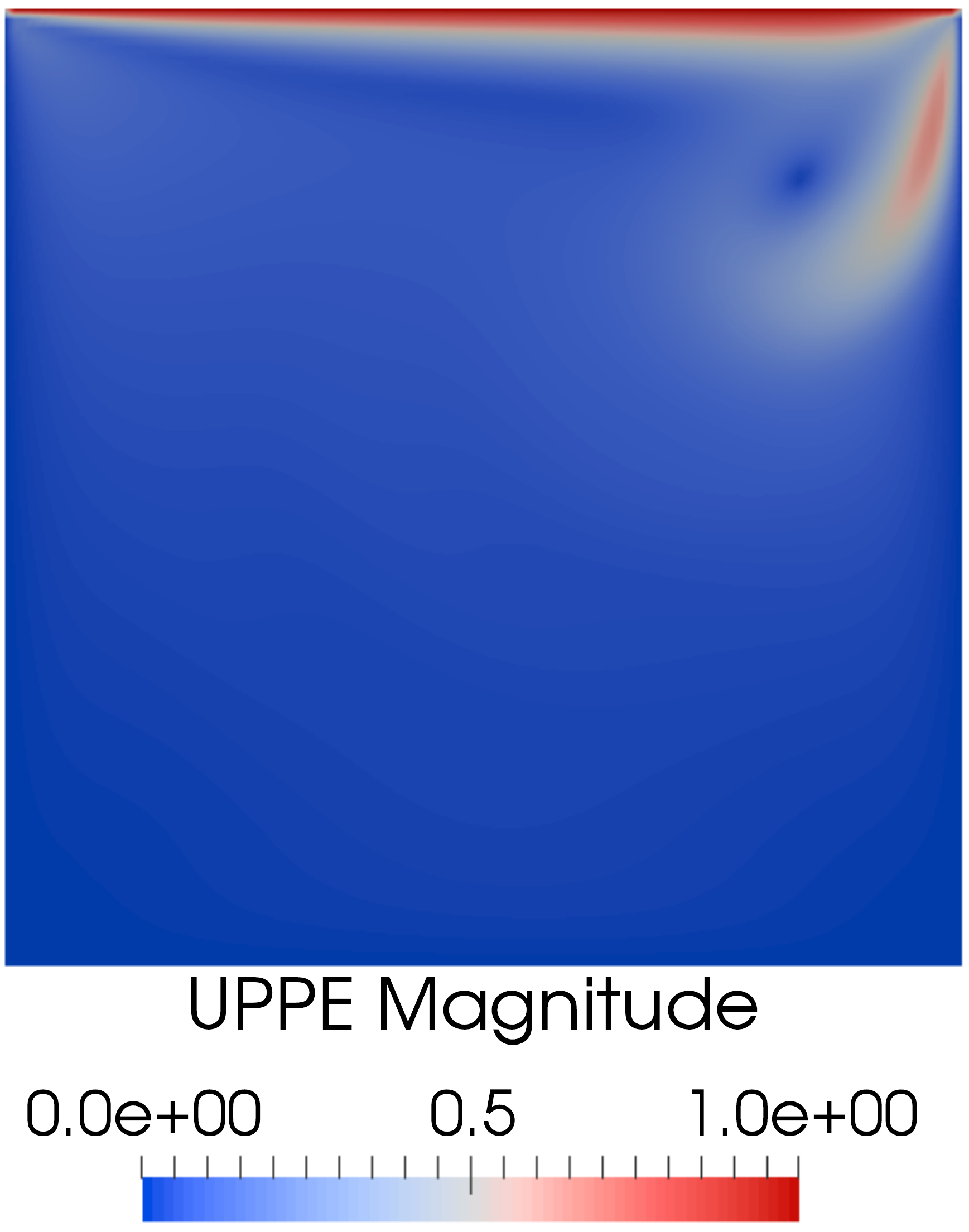}
\end{minipage}
\begin{minipage}[c]{0.16\textwidth}
\centering
\includegraphics[width=\textwidth]{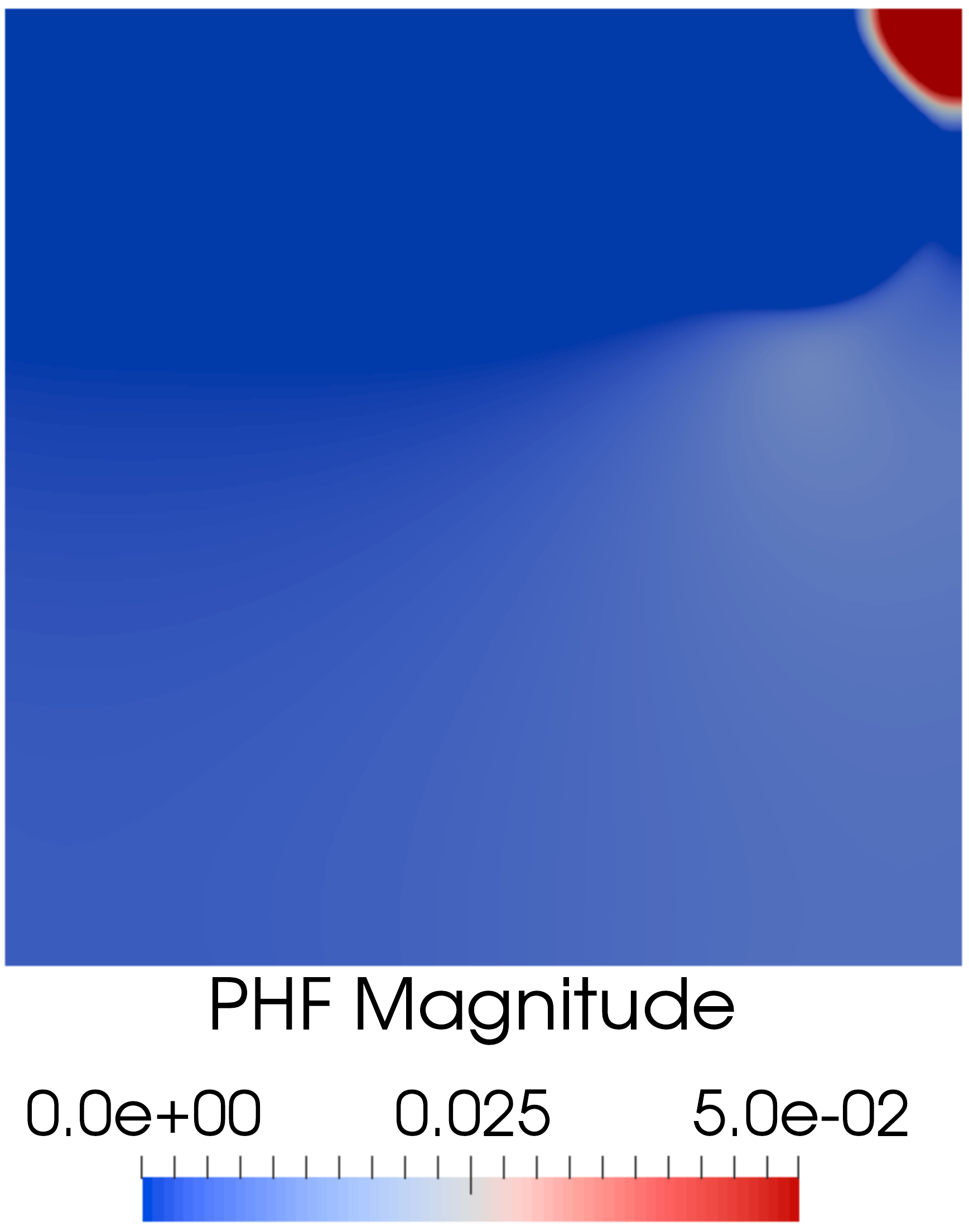}
\end{minipage}
\begin{minipage}[c]{0.16\textwidth}
\centering
\includegraphics[width=\textwidth]{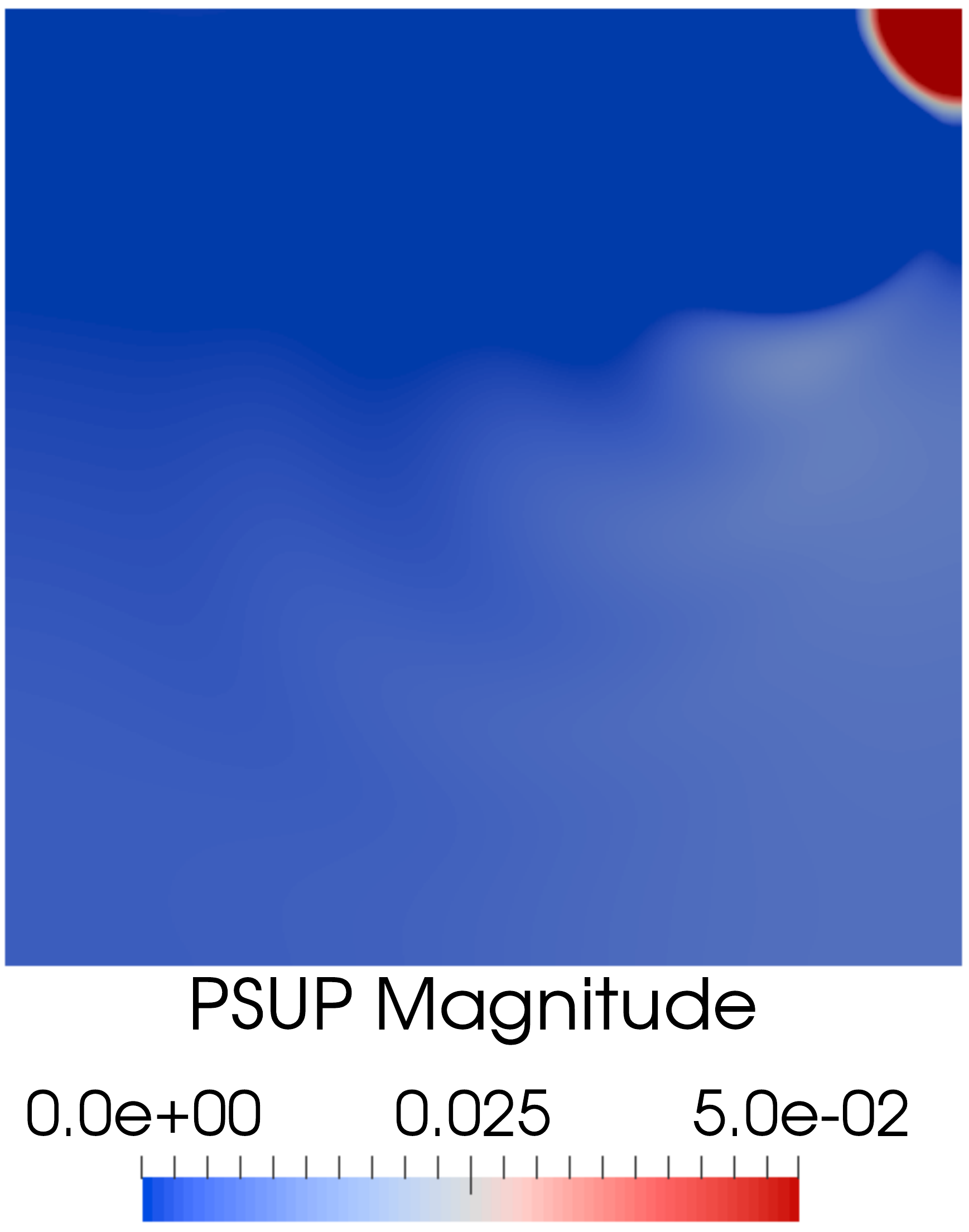}
\end{minipage}
\begin{minipage}[c]{0.16\textwidth}
\centering
\includegraphics[width=\textwidth]{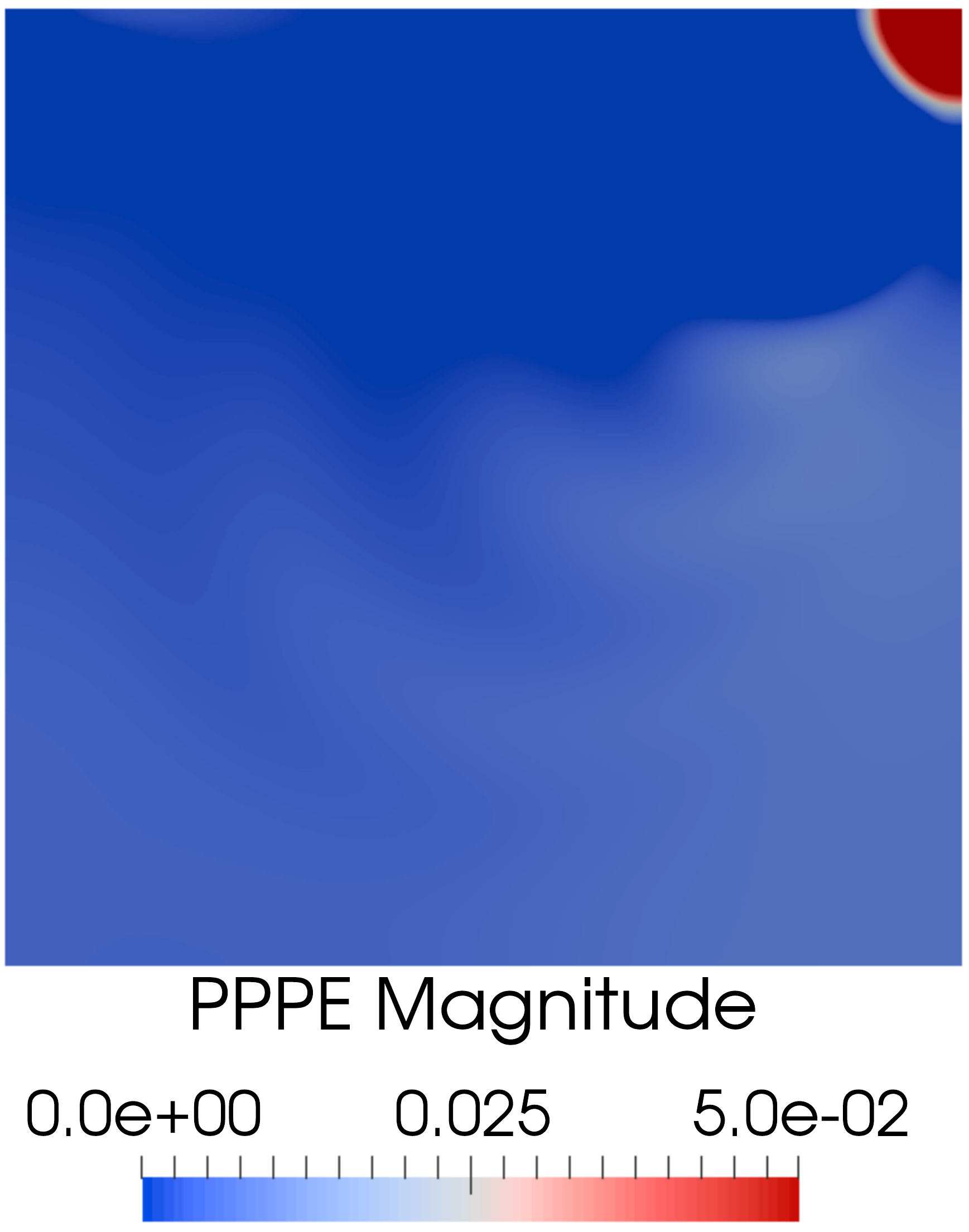}
\end{minipage}
\begin{minipage}[c]{0.16\textwidth}
\centering
\includegraphics[width=\textwidth]{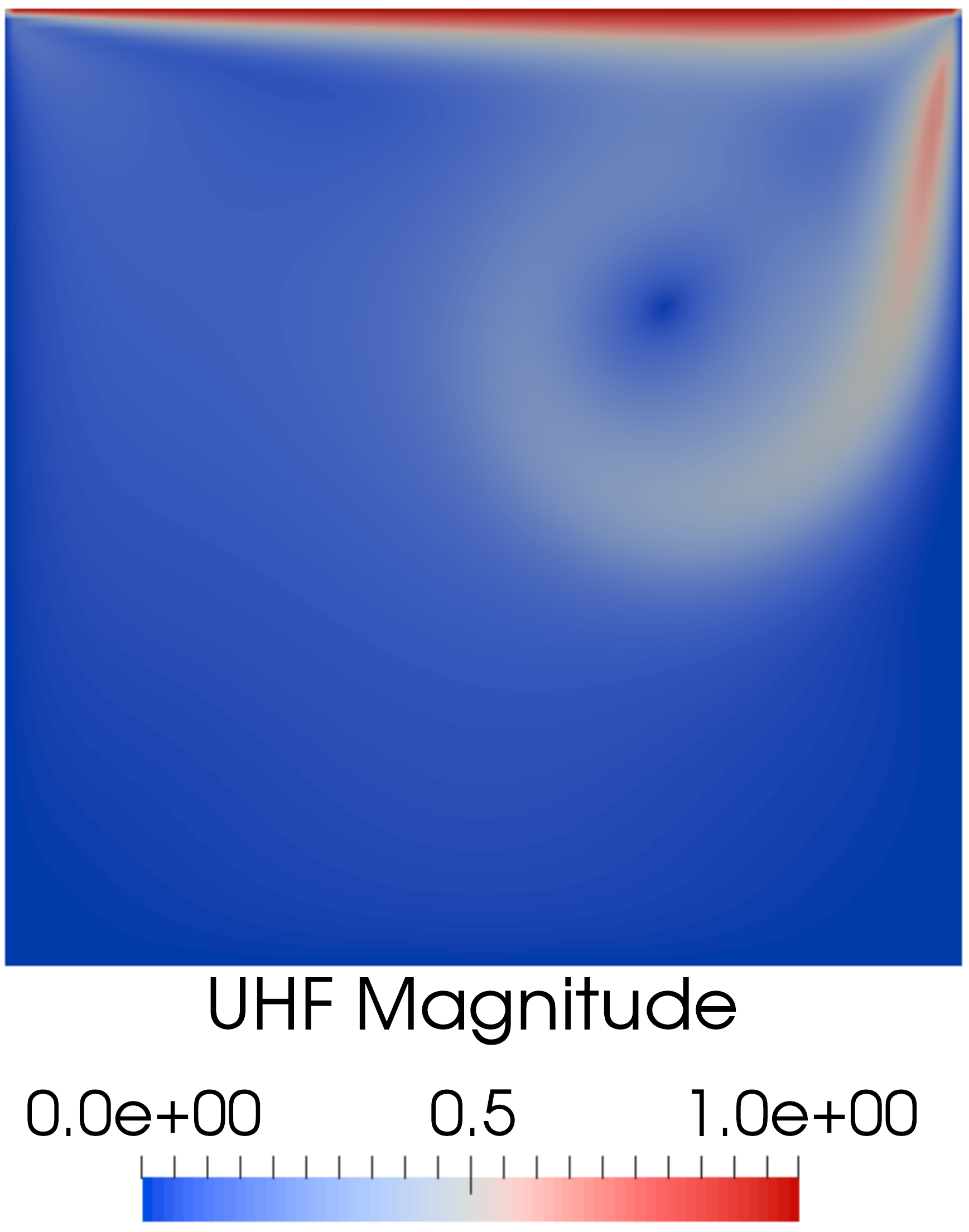}
\end{minipage}
\begin{minipage}[c]{0.16\textwidth}
\centering
\includegraphics[width=\textwidth]{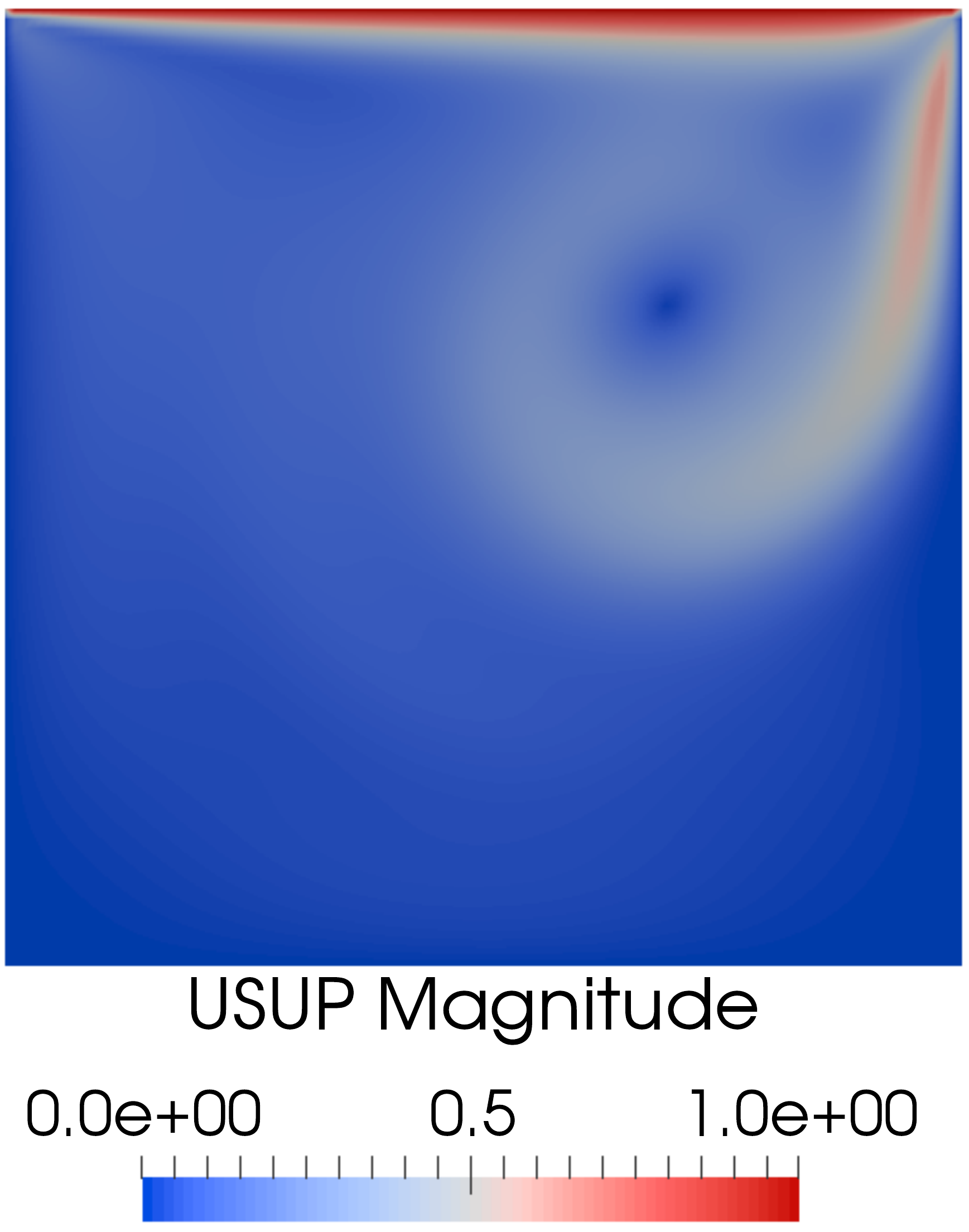}
\end{minipage}
\begin{minipage}[c]{0.16\textwidth}
\centering
\includegraphics[width=\textwidth]{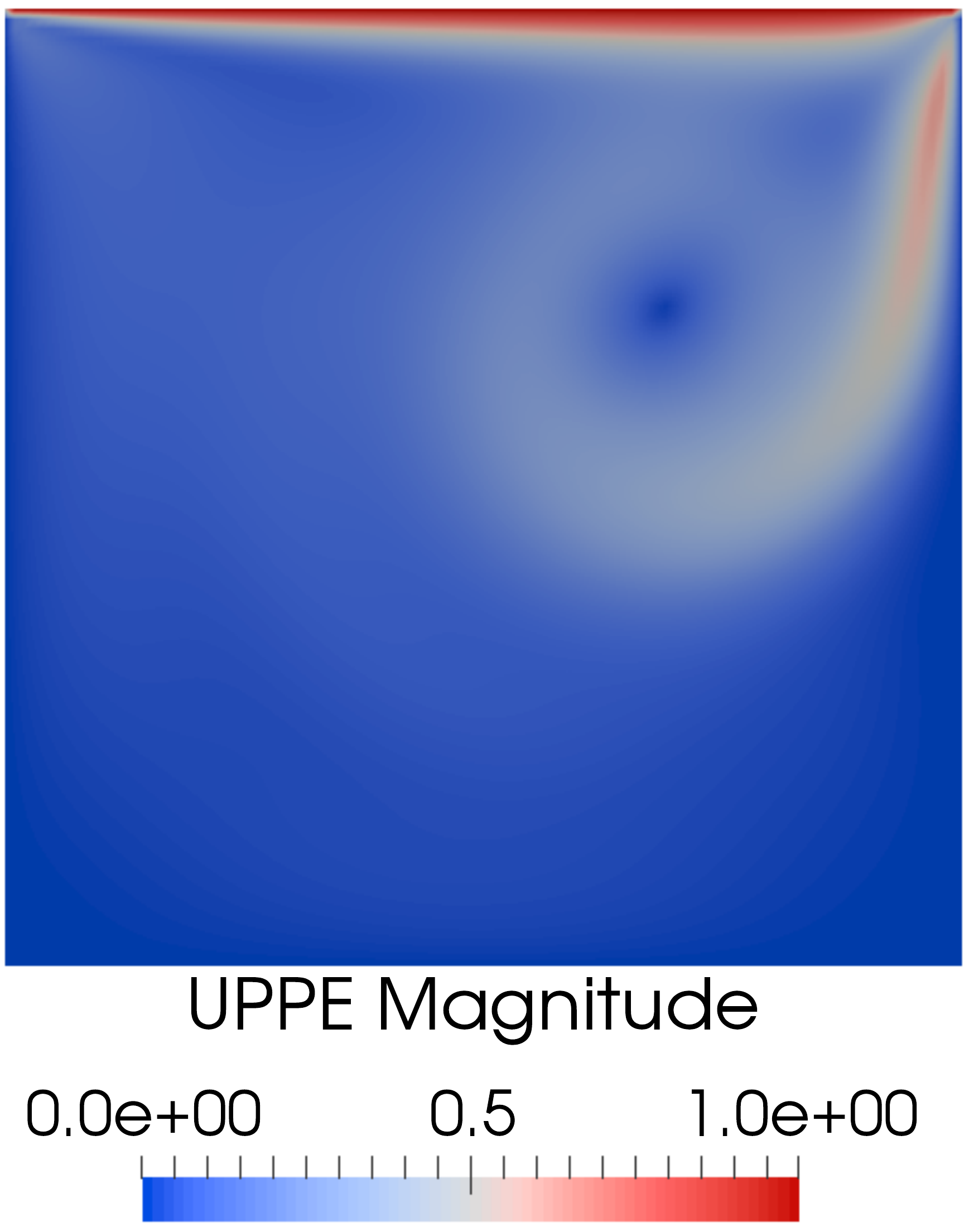}
\end{minipage}
\begin{minipage}[c]{0.16\textwidth}
\centering
\includegraphics[width=\textwidth]{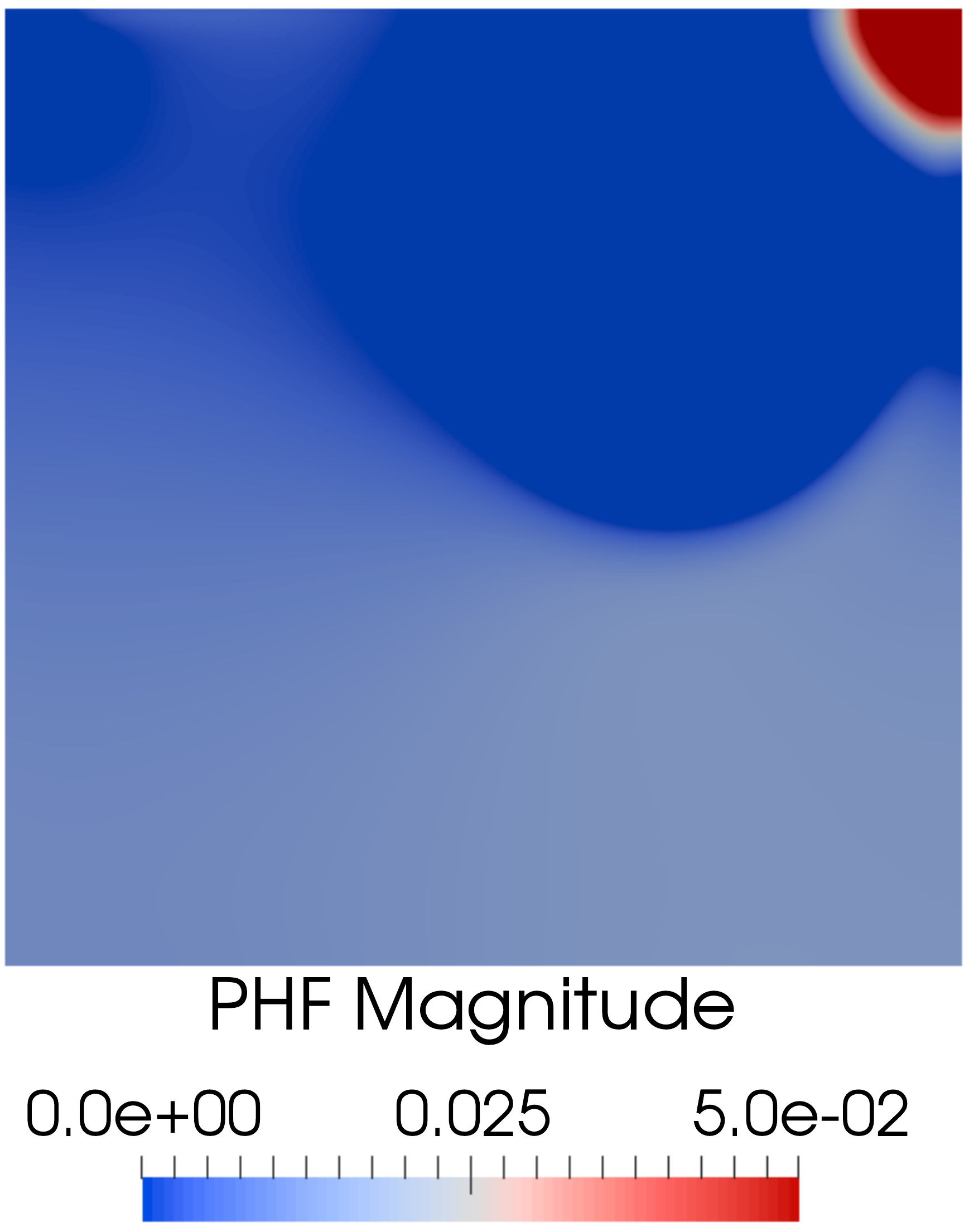}
\end{minipage}
\begin{minipage}[c]{0.16\textwidth}
\centering
\includegraphics[width=\textwidth]{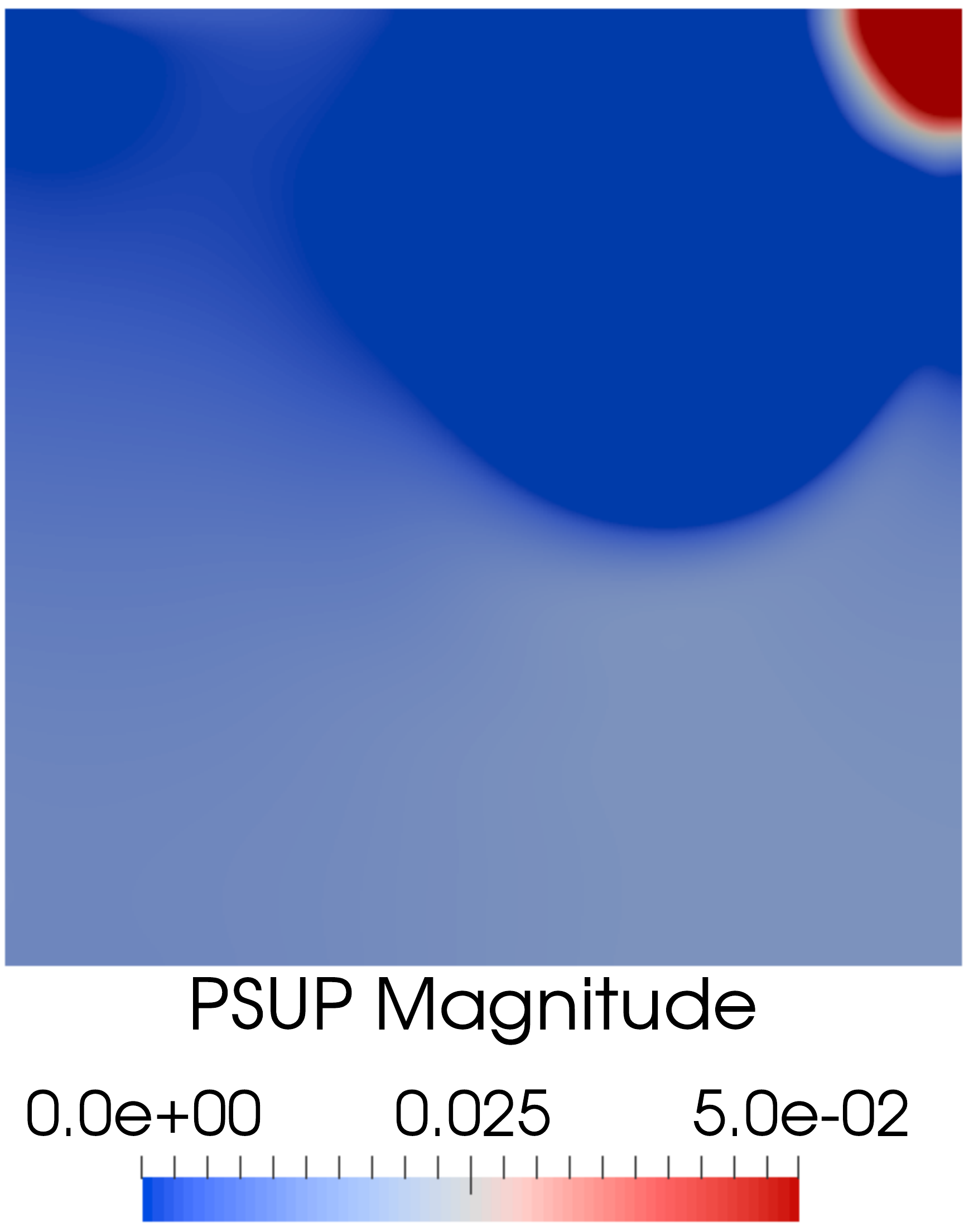}
\end{minipage}
\begin{minipage}[c]{0.16\textwidth}
\centering
\includegraphics[width=\textwidth]{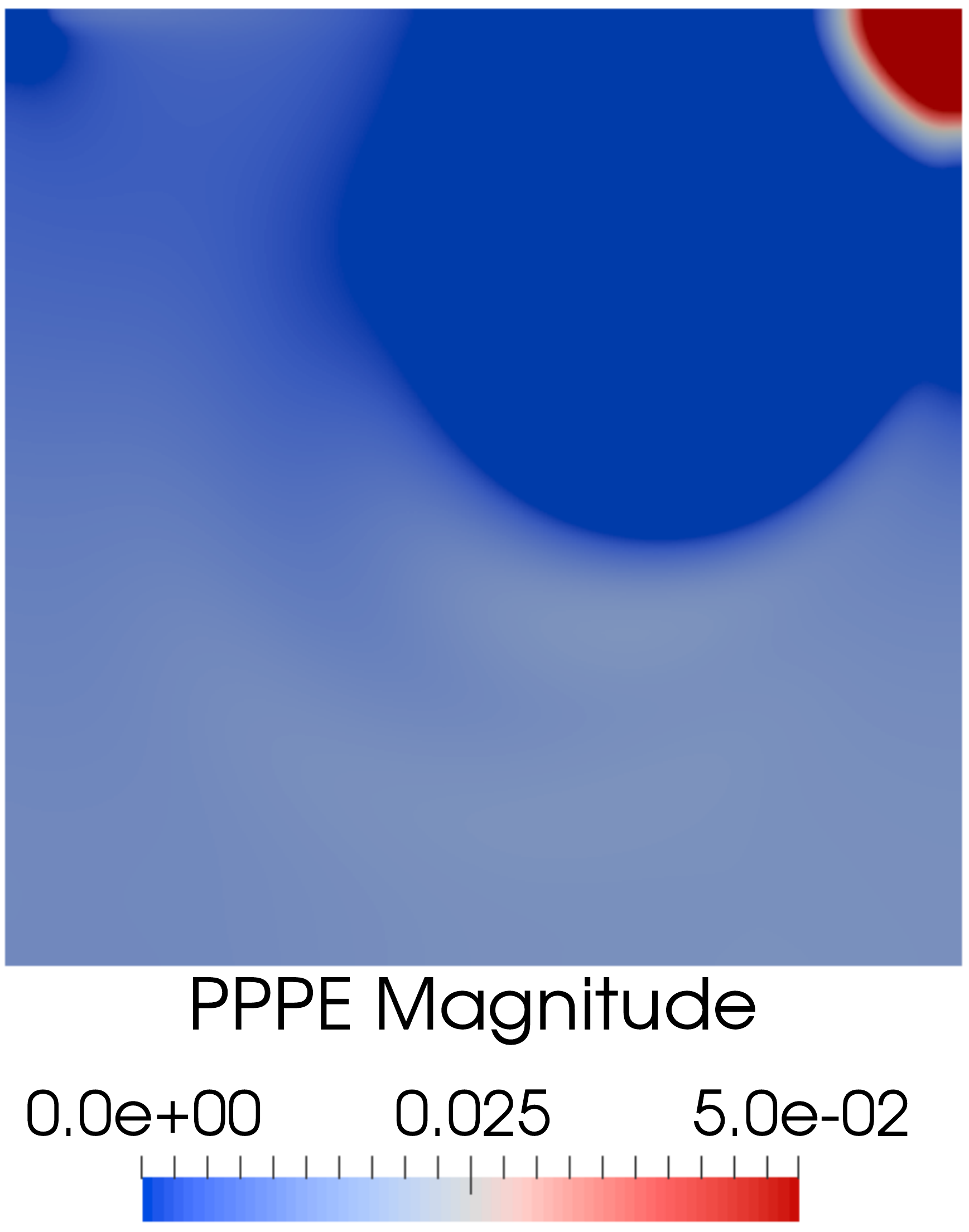}
\end{minipage}
\begin{minipage}[c]{0.16\textwidth}
\centering
\includegraphics[width=\textwidth]{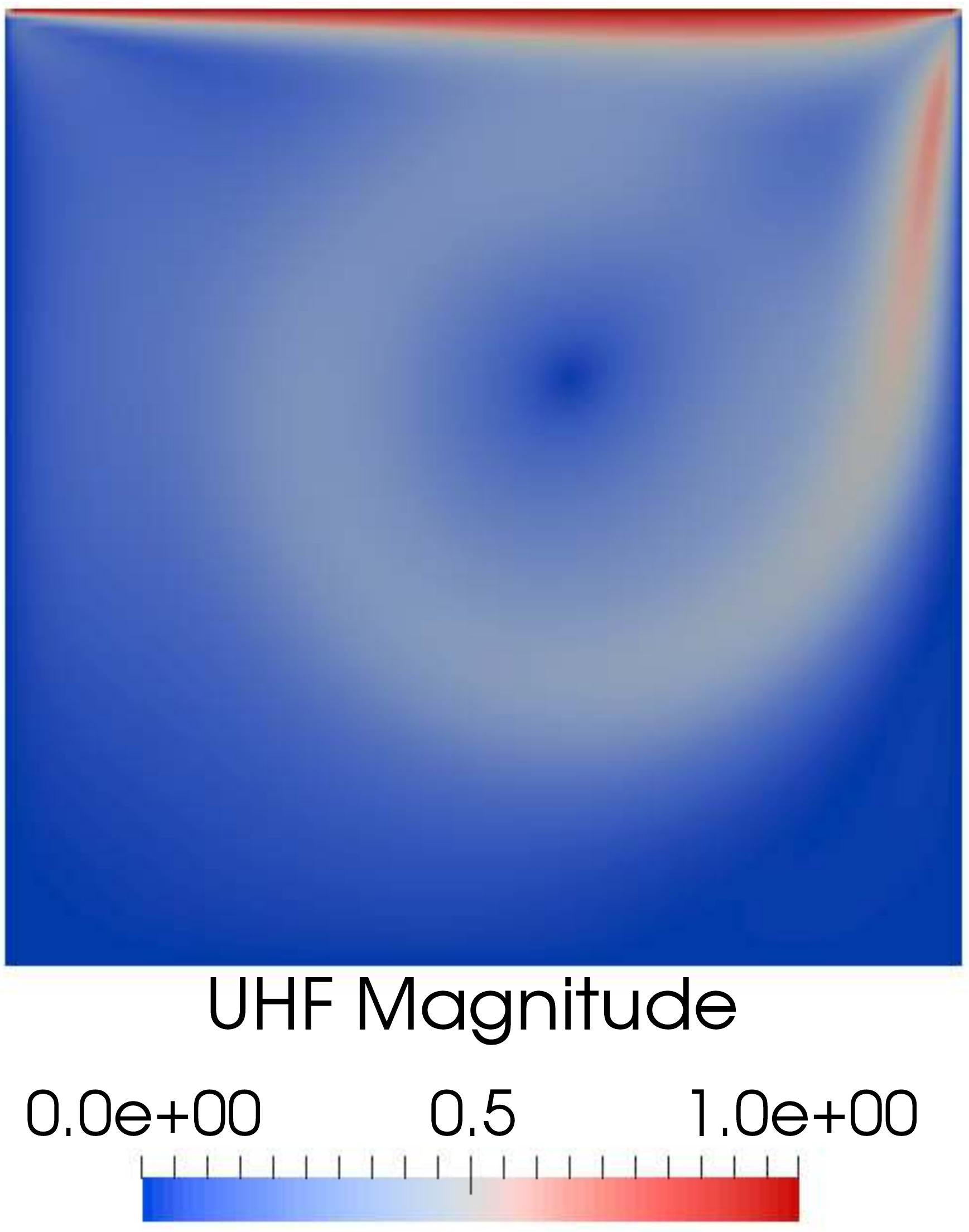}
\end{minipage}
\begin{minipage}[c]{0.16\textwidth}
\centering
\includegraphics[width=\textwidth]{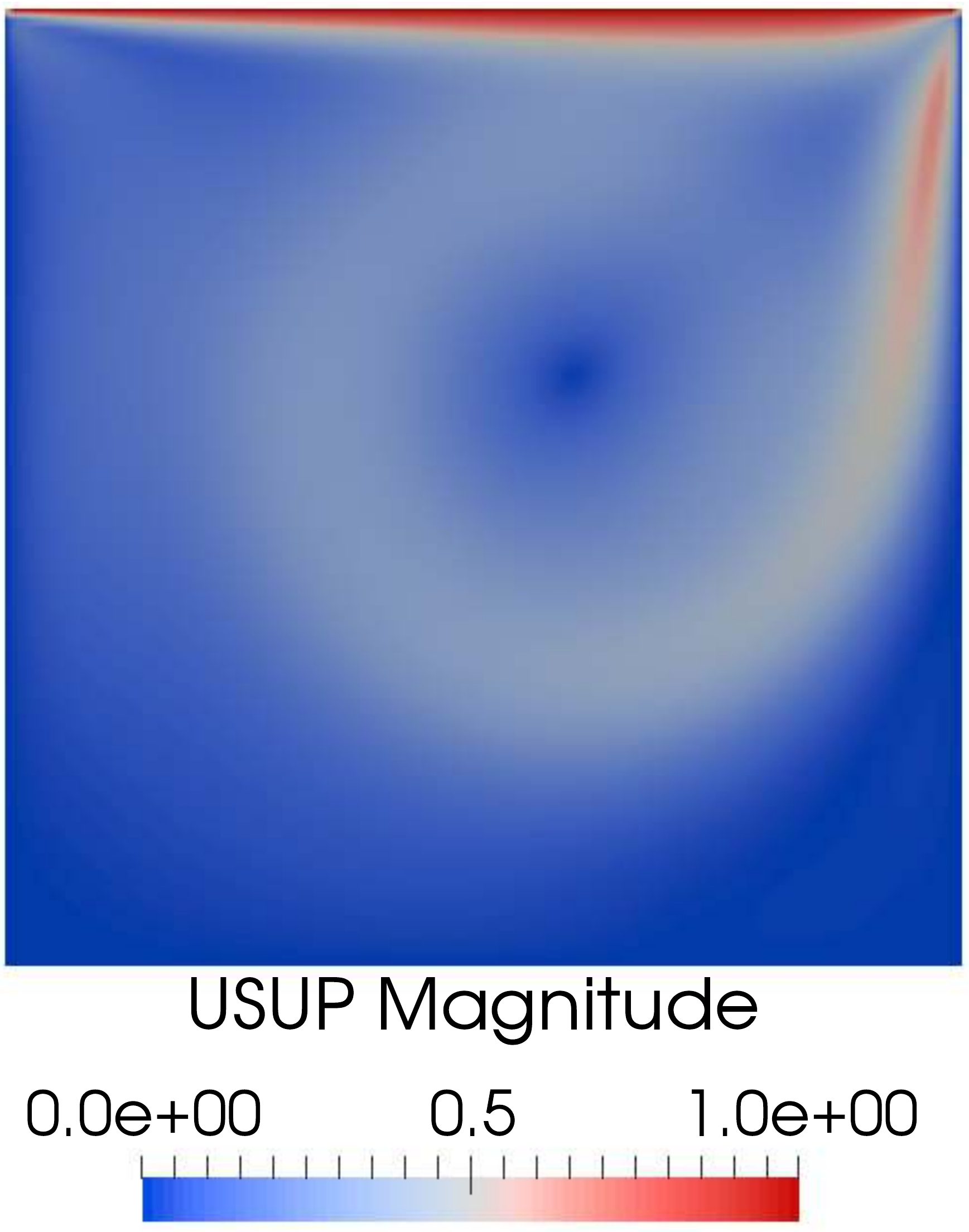}
\end{minipage}
\begin{minipage}[c]{0.16\textwidth}
\centering
\includegraphics[width=\textwidth]{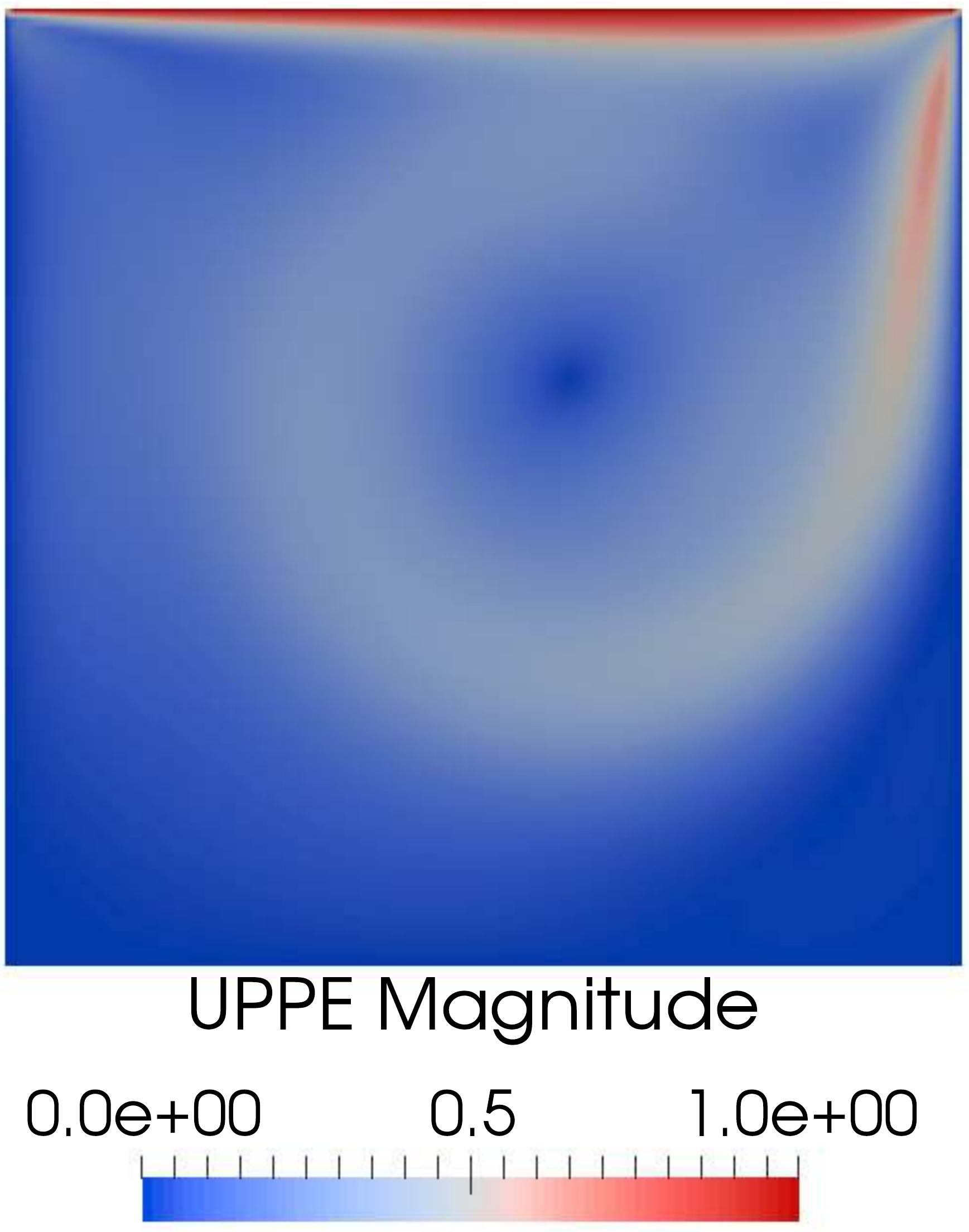}
\end{minipage}
\begin{minipage}[c]{0.16\textwidth}
\centering
\includegraphics[width=\textwidth]{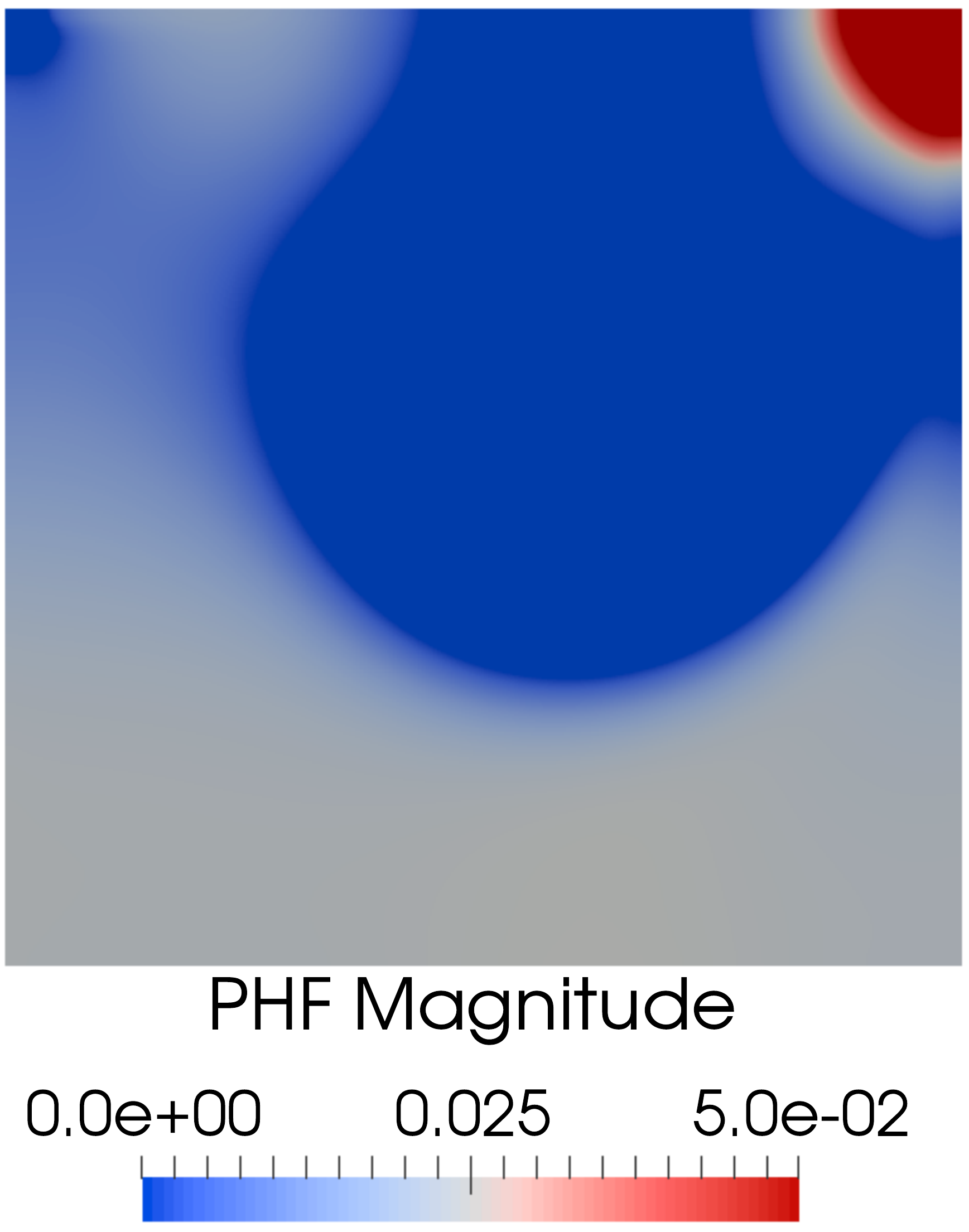}
\end{minipage}
\begin{minipage}[c]{0.16\textwidth}
\centering
\includegraphics[width=\textwidth]{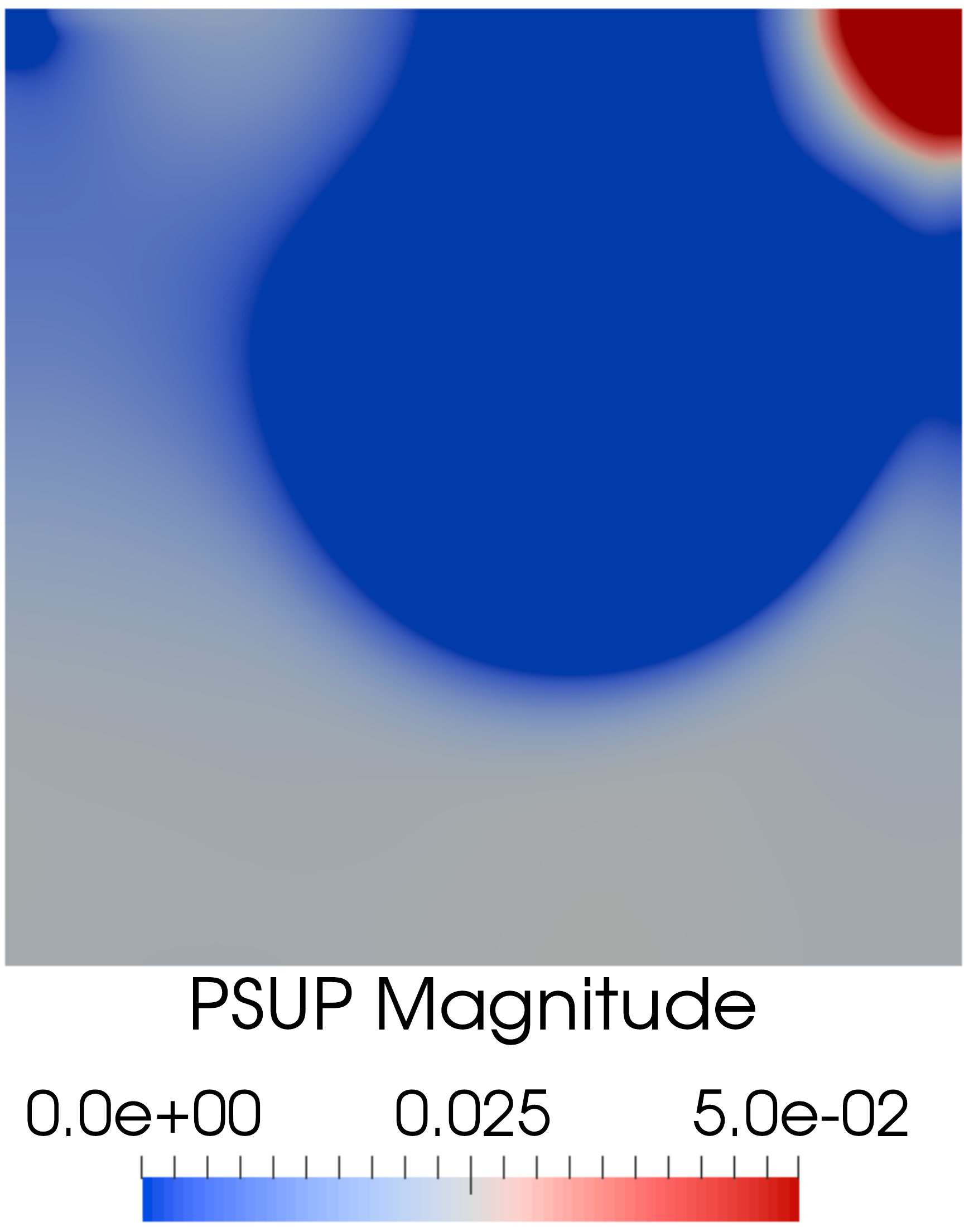}
\end{minipage}
\begin{minipage}[c]{0.16\textwidth}
\centering
\includegraphics[width=\textwidth]{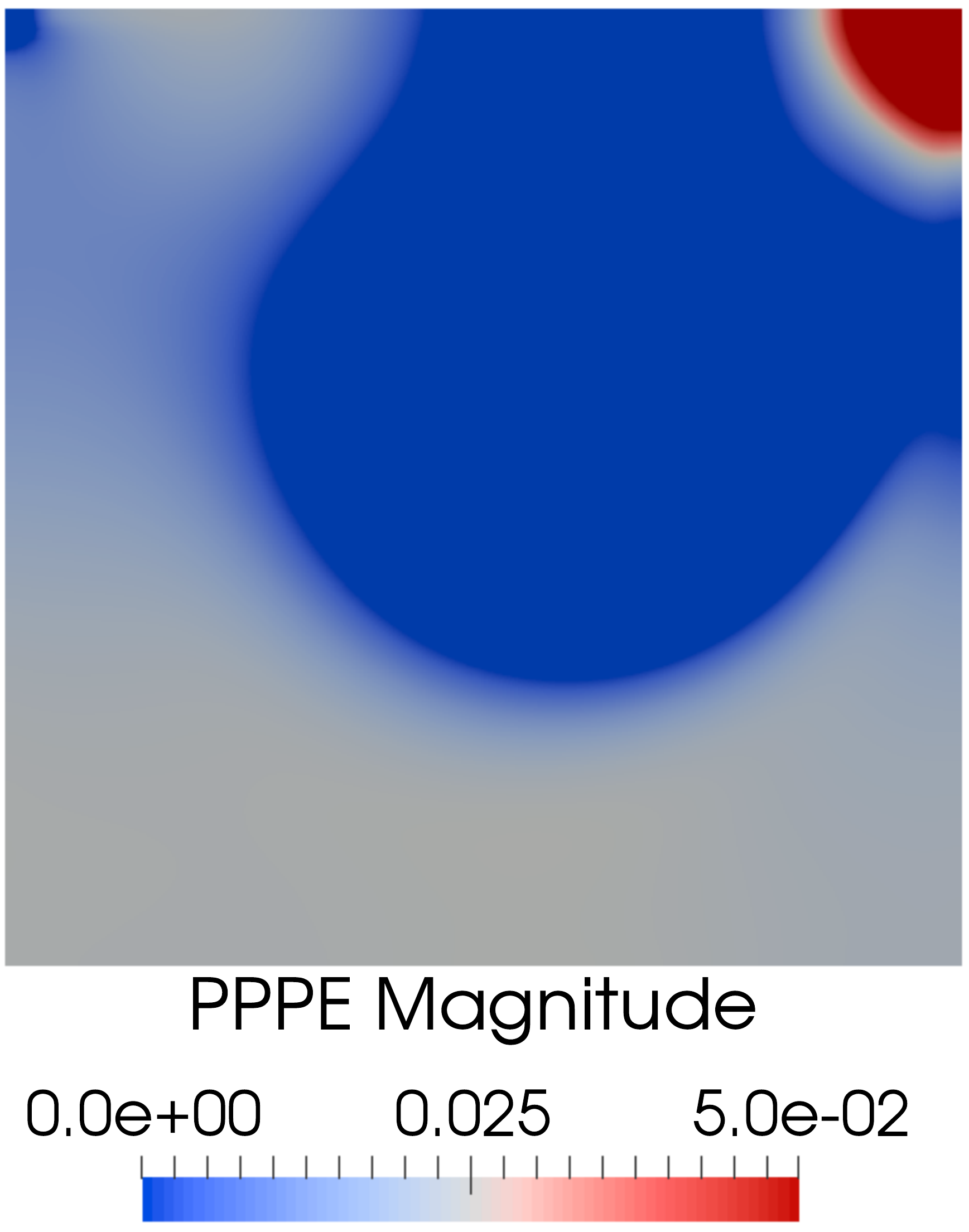}
\end{minipage}
\begin{minipage}[c]{0.16\textwidth}
\centering
\includegraphics[width=\textwidth]{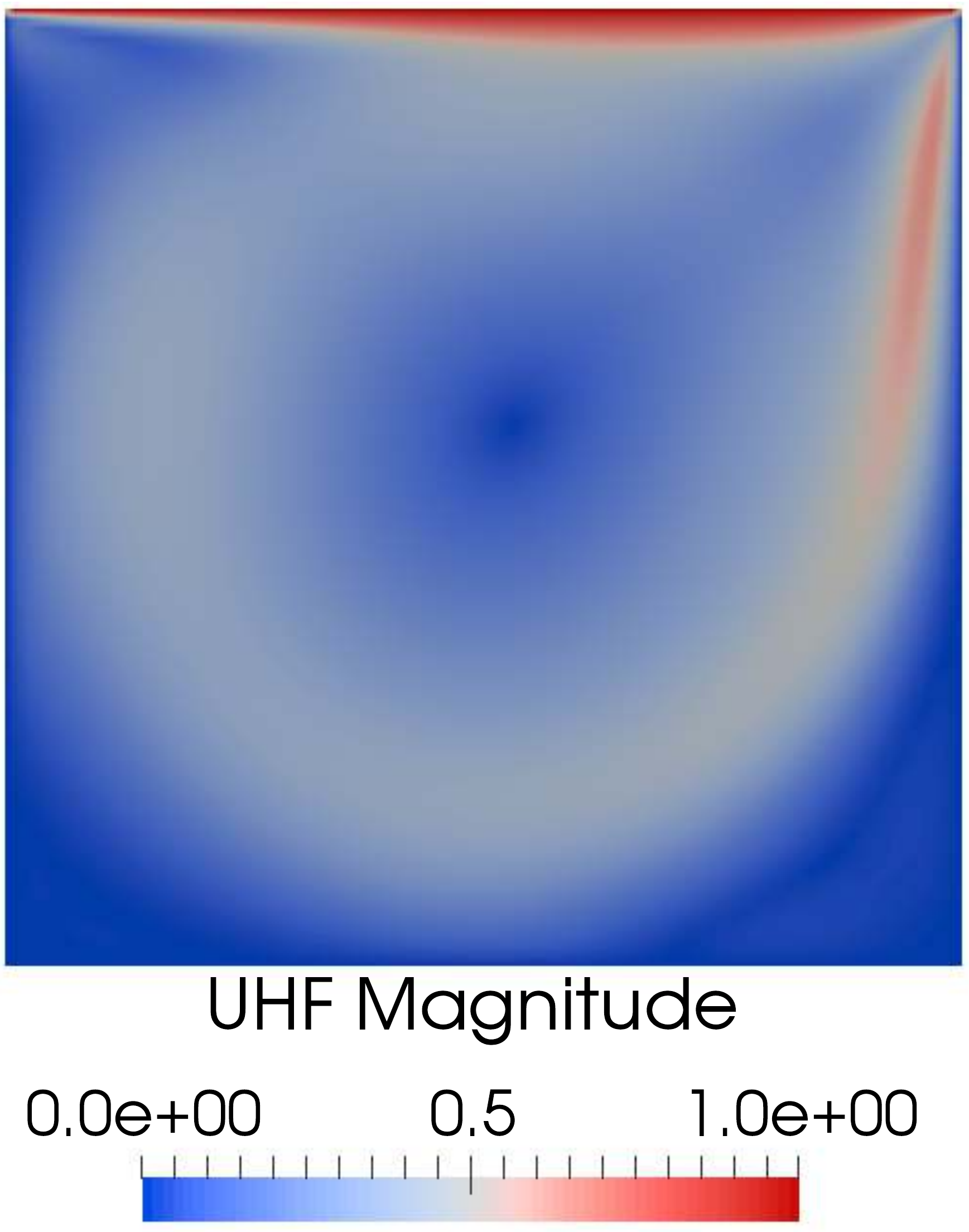}
\end{minipage}
\begin{minipage}[c]{0.16\textwidth}
\centering
\includegraphics[width=\textwidth]{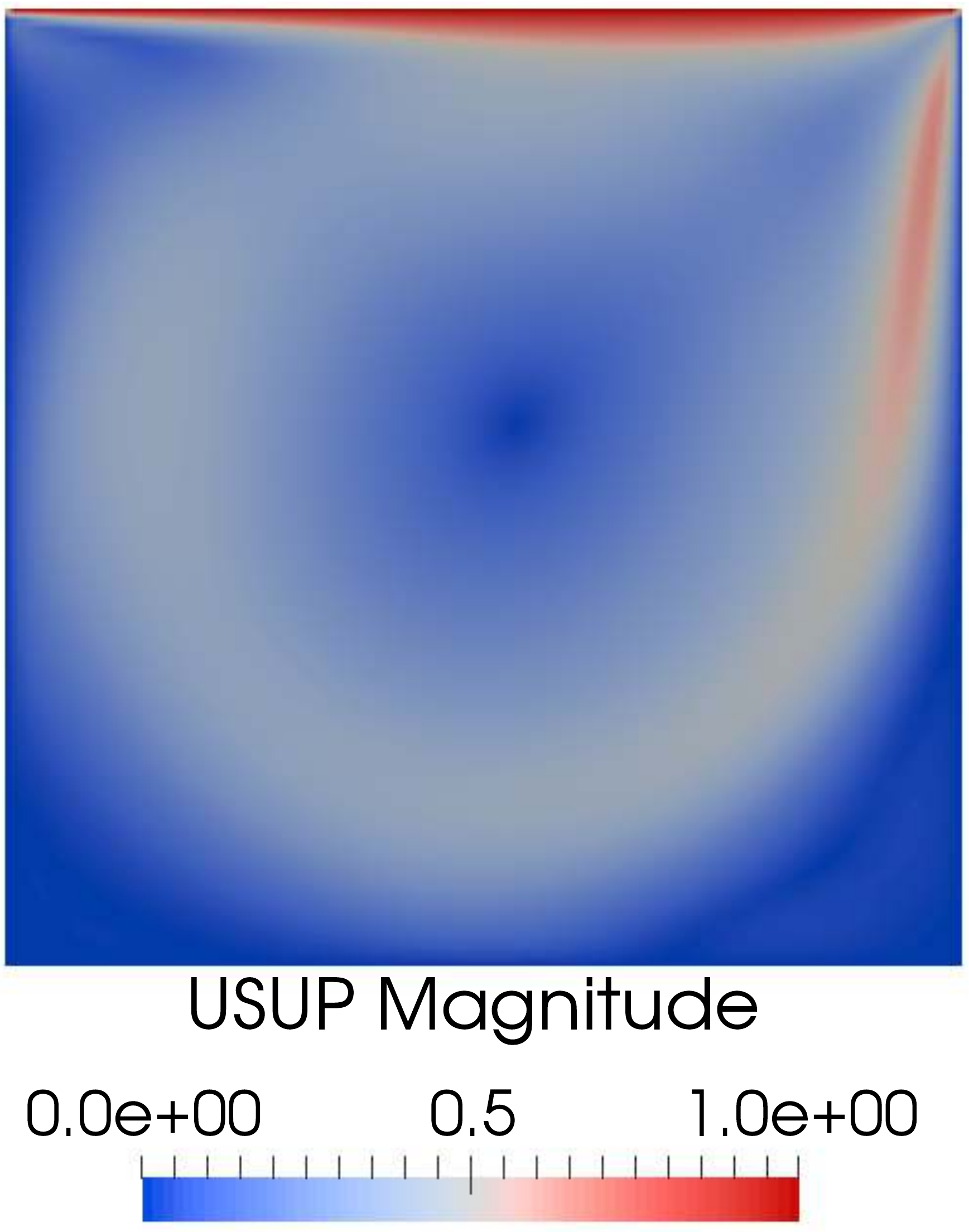}
\end{minipage}
\begin{minipage}[c]{0.16\textwidth}
\centering
\includegraphics[width=\textwidth]{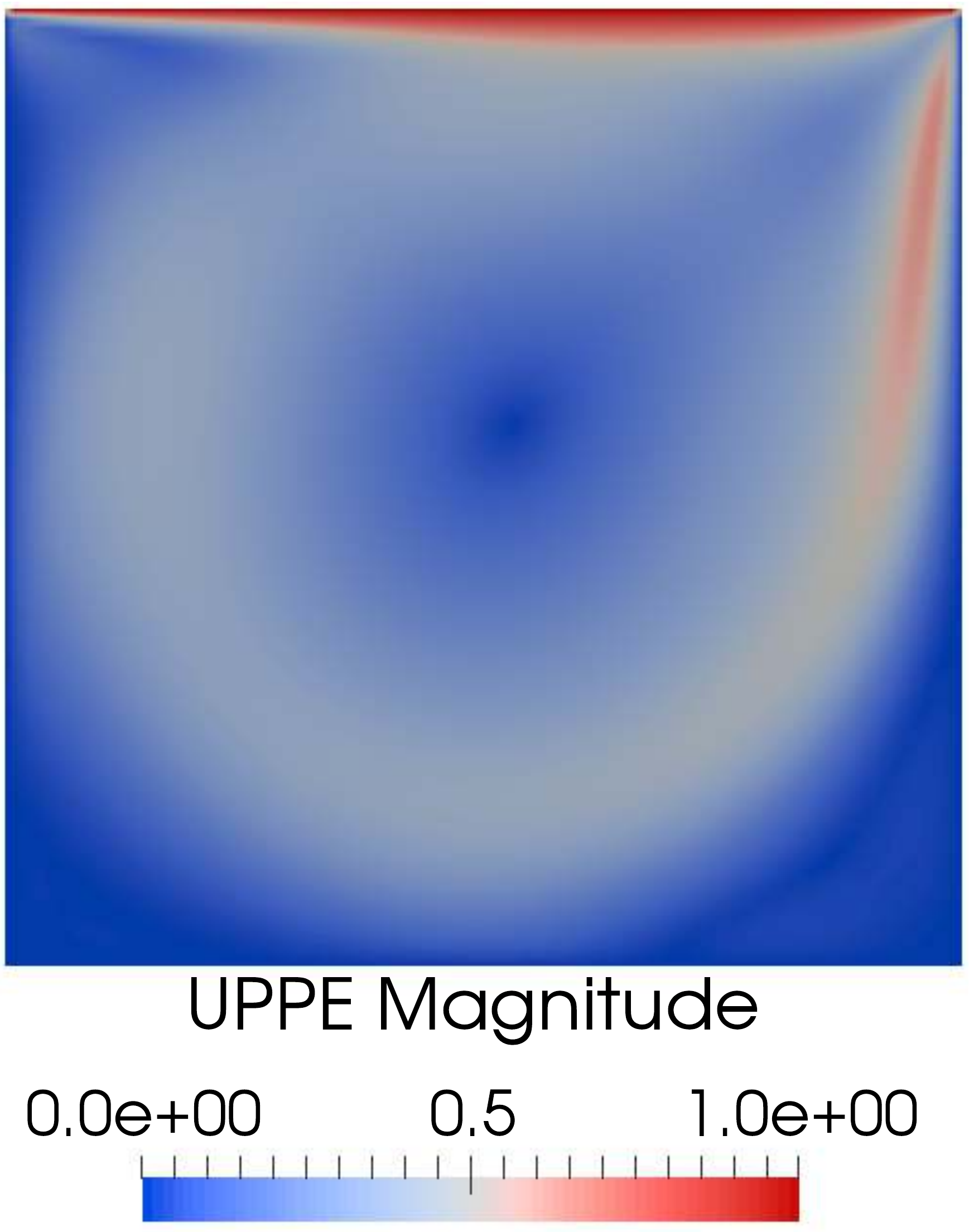}
\end{minipage}
\begin{minipage}[c]{0.16\textwidth}
\centering
\includegraphics[width=\textwidth]{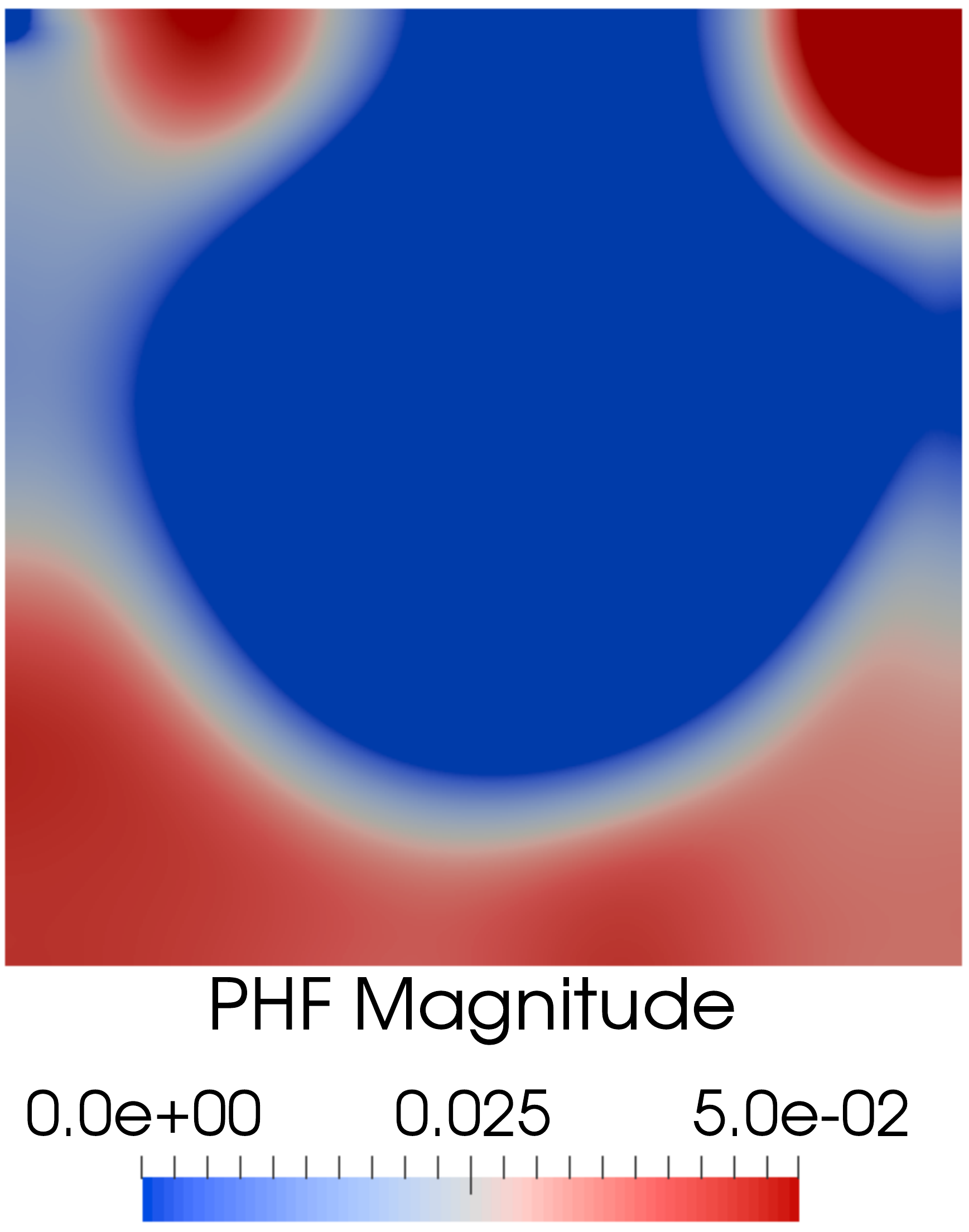}
\end{minipage}
\begin{minipage}[c]{0.16\textwidth}
\centering
\includegraphics[width=\textwidth]{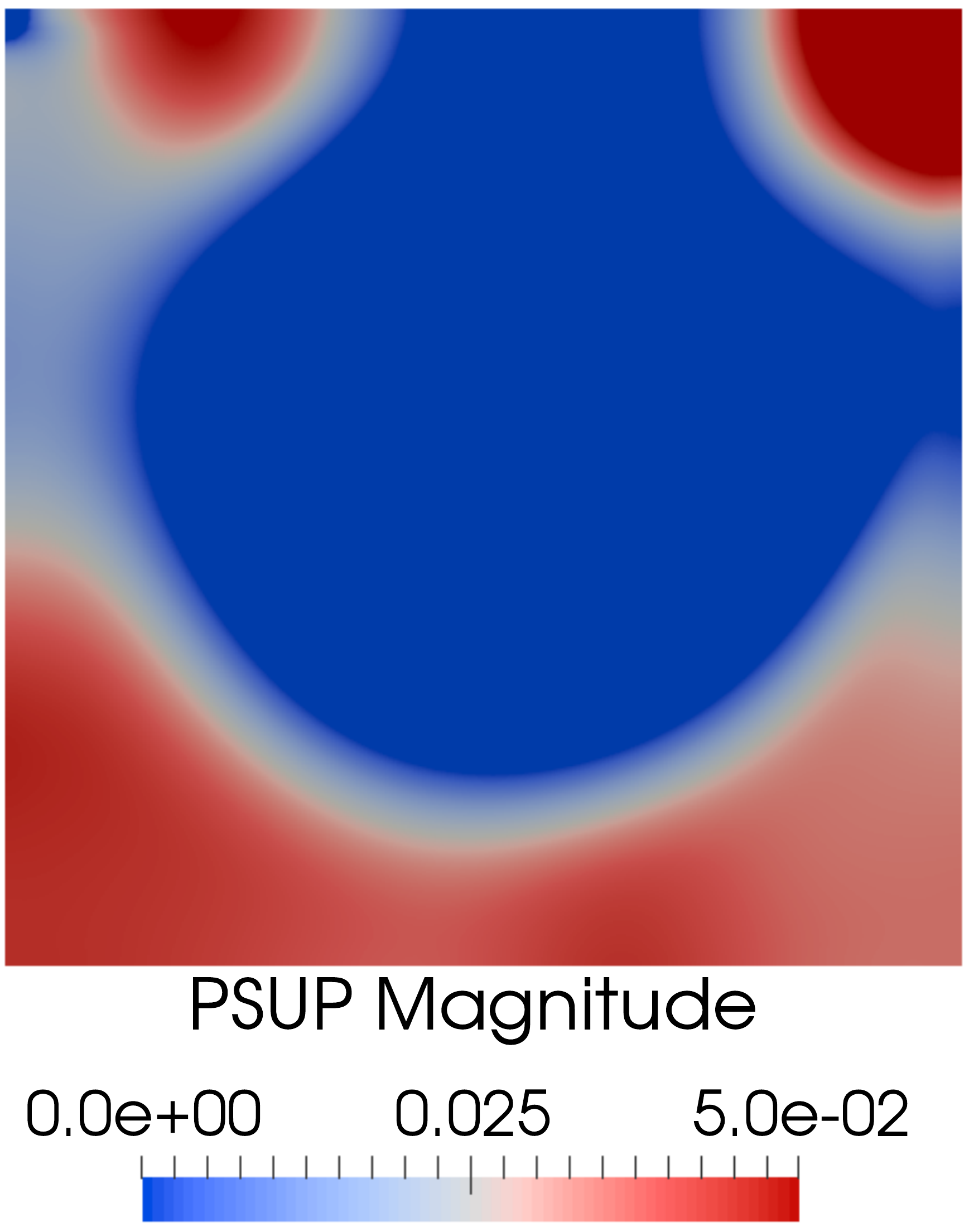}
\end{minipage}
\begin{minipage}[c]{0.16\textwidth}
\centering
\includegraphics[width=\textwidth]{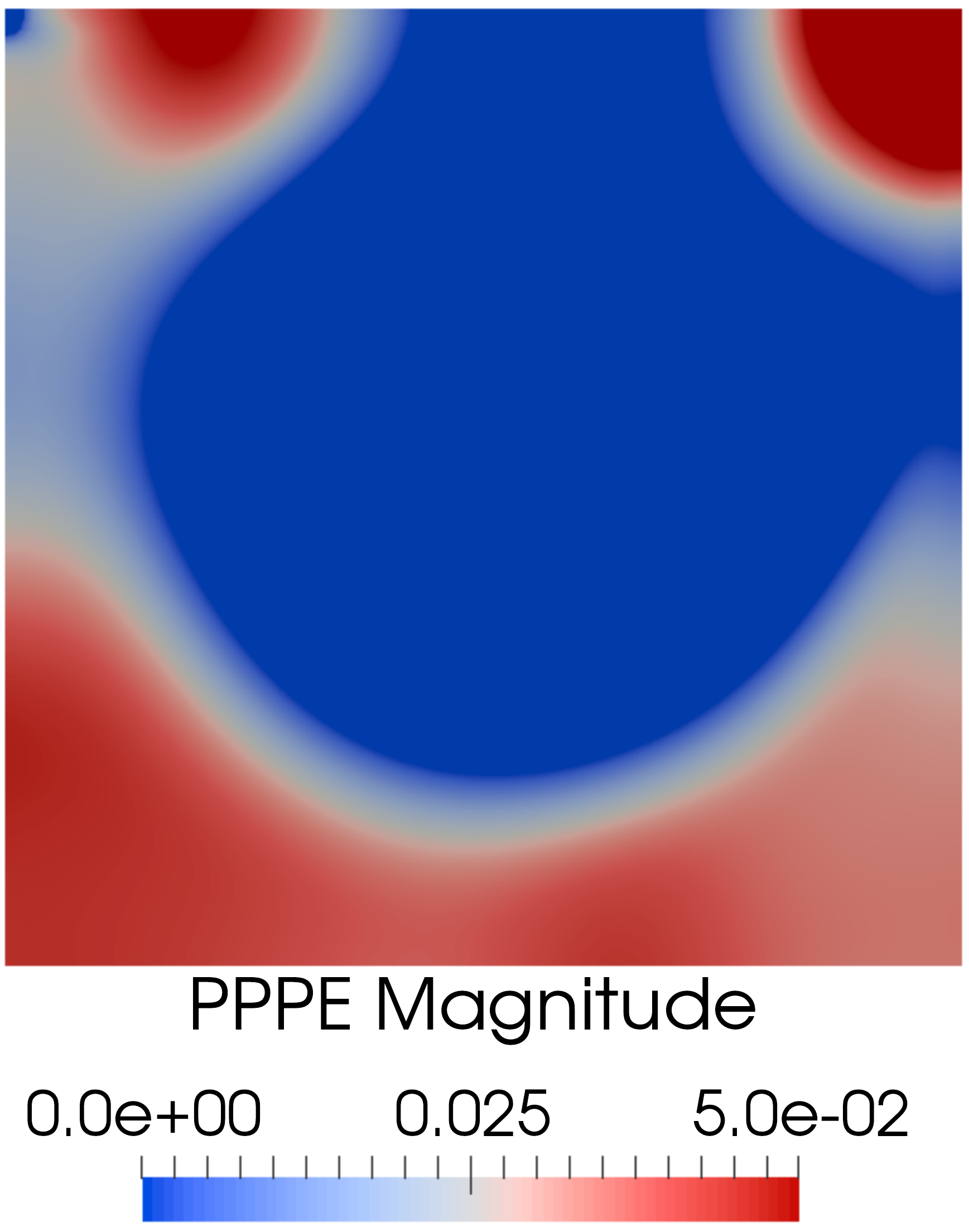}
\end{minipage}
\caption{Comparison of the velocity and pressure fields for high fidelity (UHF - $1$ column, PHF - $4$ column), SUP-ROM (USUP - $2$ column, PSUP - $5$ column) and PPE-ROM (UPPE $3$ column, PPEE $6$ - column). The fields are depicted for different time instant equal to $t=0.2 \si{s},0.5 \si{s},1 \si{s}$ and $5 \si{s}$, respectively, and increasing in the image from top to bottom. The ROM models are obtained with $10$ modes for velocity and pressure and only for the SUP-ROM with 10 additional supremizer modes. The velocity and pressure magnitudes are shown in the image legends.}\label{fig:comparison_cavity}
\end{figure*}
\begin{table*}
{
\centering
{
\begin{tabular}{ c | c | c | c | c}
N Modes & $\bm{u}$ & $p$ & $\bm{s}$ & $\beta$ \\ 
\hline
1&0.978946&0.975406&0.980260  & 9.264e-05  \\ 
2&0.994184&0.991528&0.995232  & 9.264e-05  \\ 
3&0.997737&0.995385&0.997912  & 7.175e-04  \\ 
4&0.998990&0.998116&0.999400  & 7.175e-04  \\ 
5&0.999483&0.999270&0.999844  & 7.175e-04  \\ 
10&0.999971&0.999971&0.999997 & 1.551e-02
\end{tabular}}\caption{The table contains the cumulative eigenvalues for the lid driven cavity test. The first, second and third columns report the cumulative eigenvalues for the velocity, pressure and supremizer fields, respectively. The last column contains the value of the inf-sup constant, in the supremizer stabilisation case, for different different number of supremizer modes and with a fixed number of velocity and pressure modes (10 modes for velocity and 10 modes for pressure)}\label{tab:cum_eig}
}
\end{table*}

\subsection{Flow around a circular cylinder}
\begin{figure*}
\begin{minipage}{0.48\textwidth}
\centering
\includegraphics[width=\textwidth]{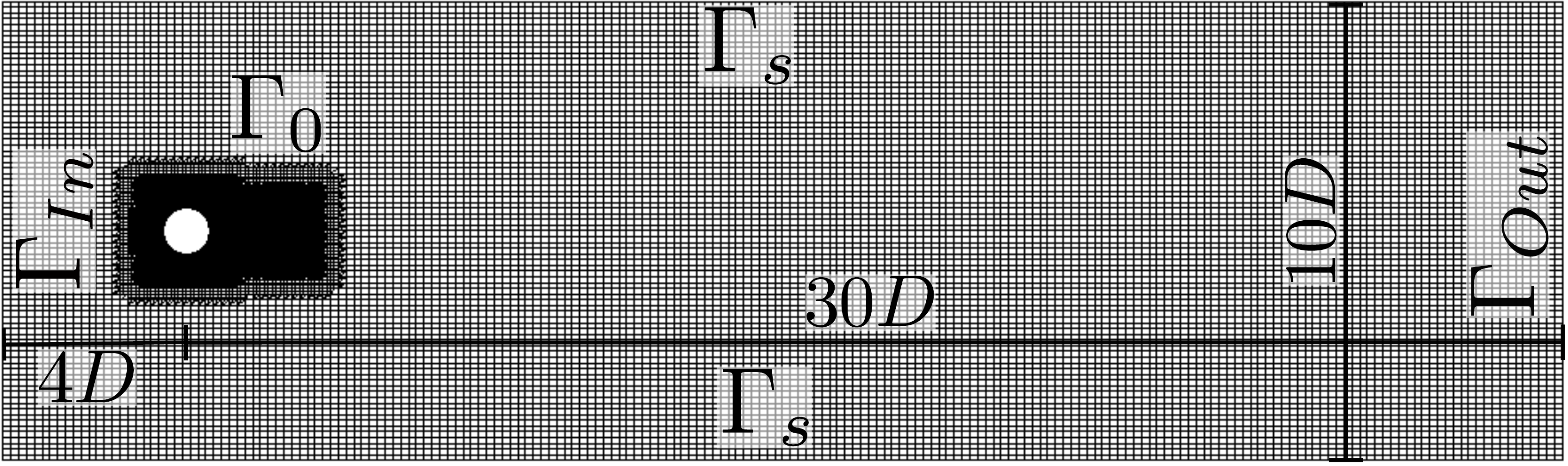}\vspace{2mm}\\
\resizebox{\textwidth}{!}{%
\begin{tabular}{ l | c | c | c | c}
    & $\Gamma_{In}$ & $\Gamma_{0}$ & $\Gamma_{s}$ & $\Gamma_{Out}$\\ \hline 
    $\bm{u}$ & $\bm{u} = (1,0)$ & $\bm{u} = (0,0)$ & $\bm{u}\cdot \bm{n} = 0$ & $\nabla \bm{u} \cdot \bm{n} = 0$ \\\hline 
    $p$ & $\nabla p \cdot \bm{n} = 0$ & $\nabla p \cdot \bm{n} = 0$ & $\nabla p \cdot \bm{n} = 0$ & $p = 0$
\end{tabular}}\vspace{2mm}\\
\includegraphics[width=\textwidth]{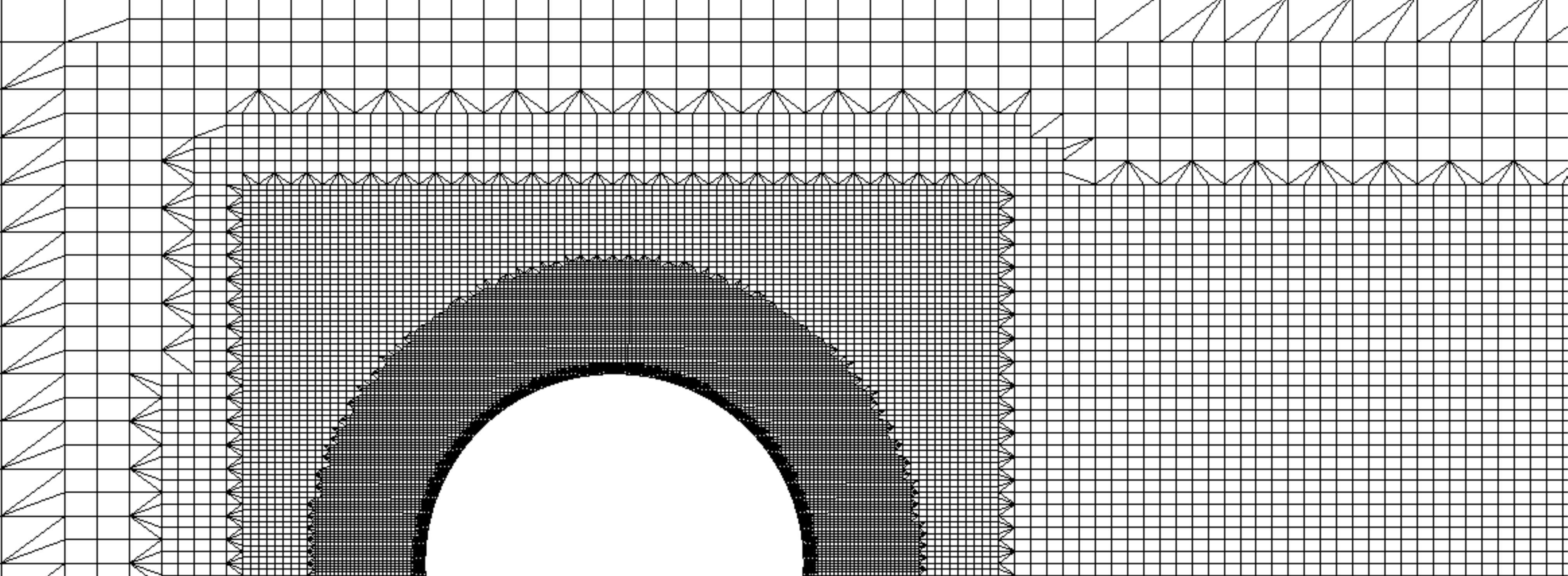}
\end{minipage}
\caption{The figure show from the top to the bottom: a general overview of the mesh with dimension and boundaries, a table with the imposed values at the boundaries, and a zoom of the mesh near to the cylinder.}\label{fig:mesh_cyl} 
\end{figure*}

\begin{figure*}
\begin{minipage}{1\textwidth}
\centering
\includegraphics[width=0.24\textwidth]{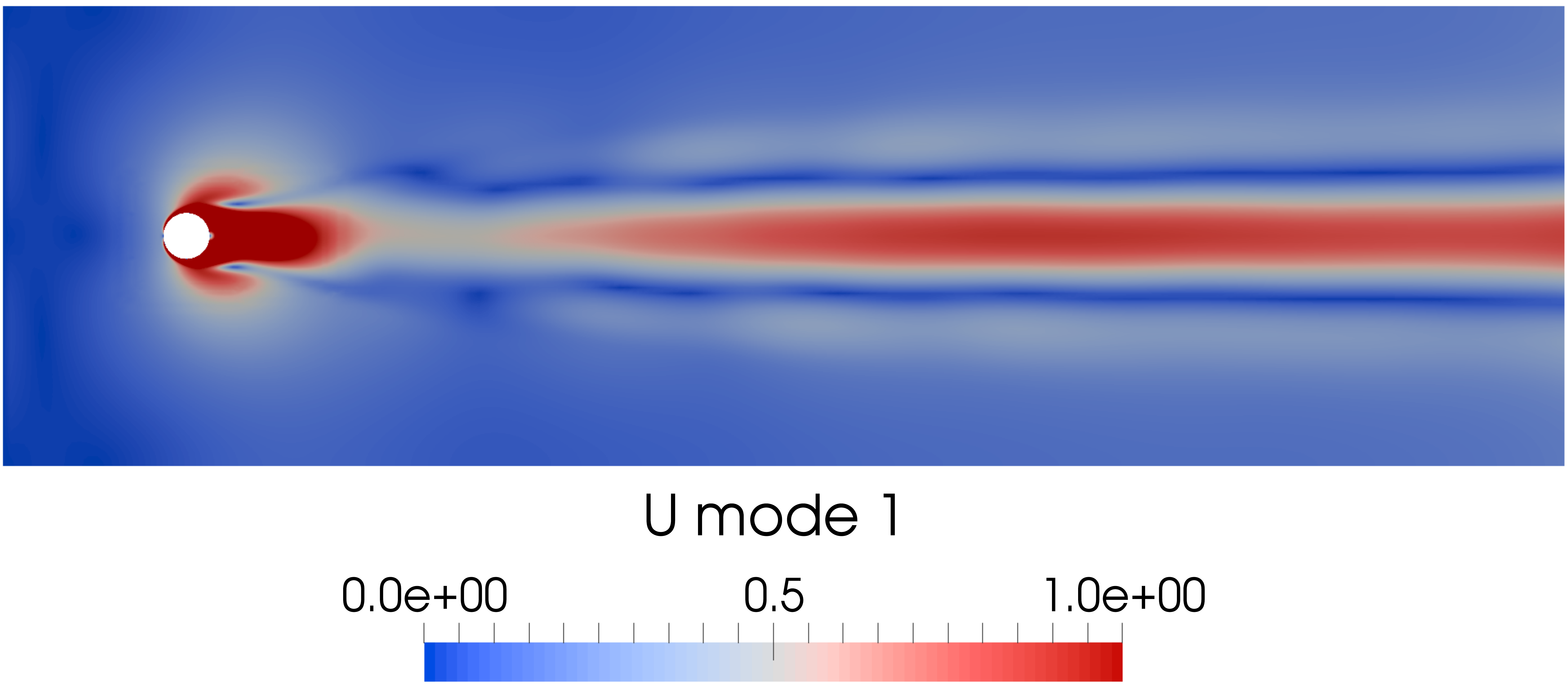}
\includegraphics[width=0.24\textwidth]{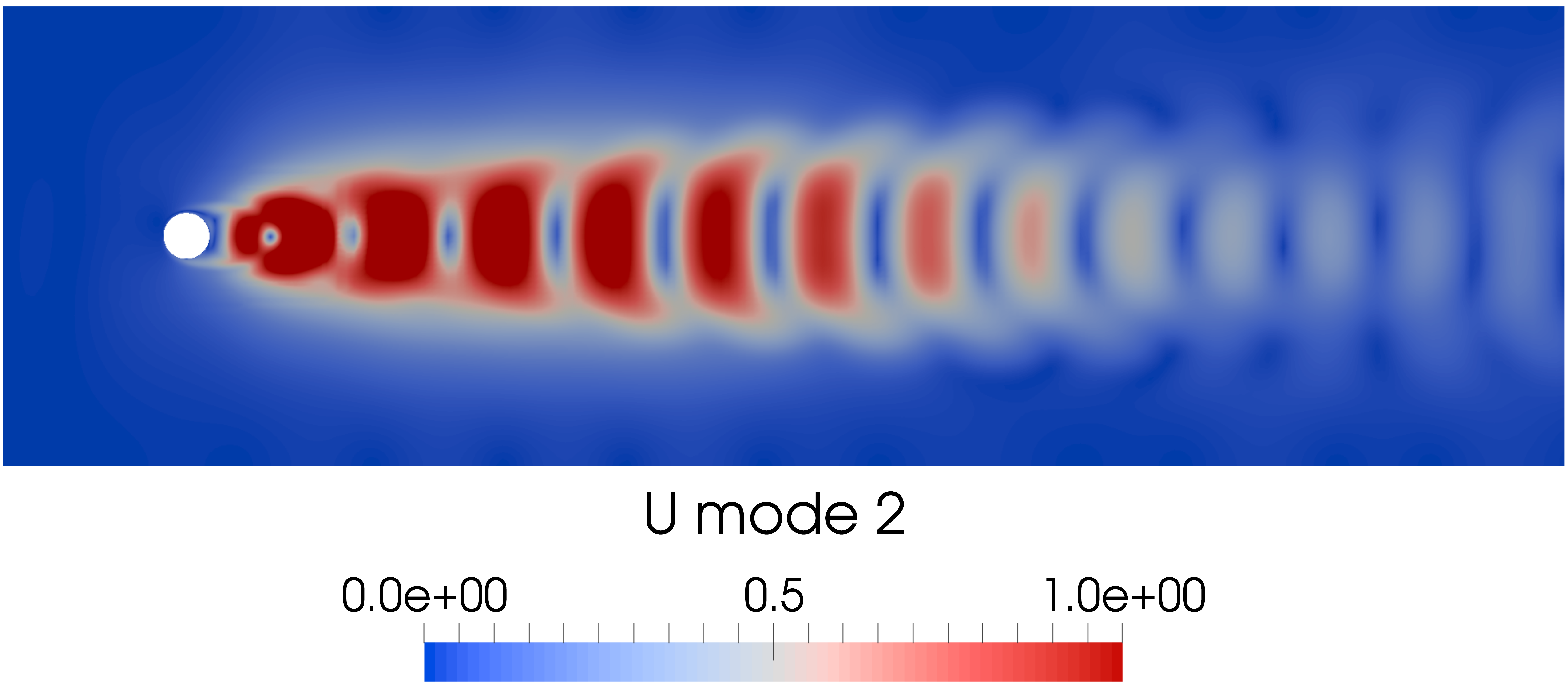}
\includegraphics[width=0.24\textwidth]{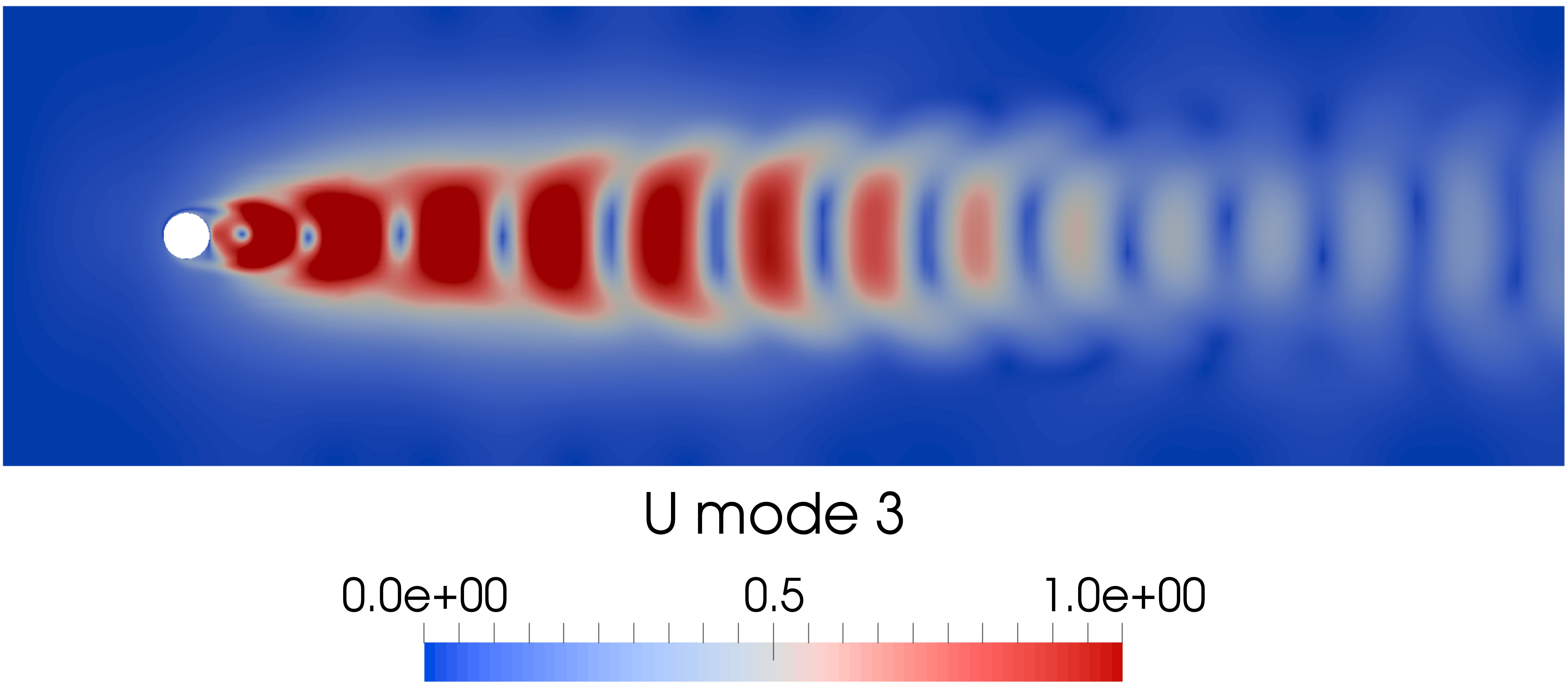}
\includegraphics[width=0.24\textwidth]{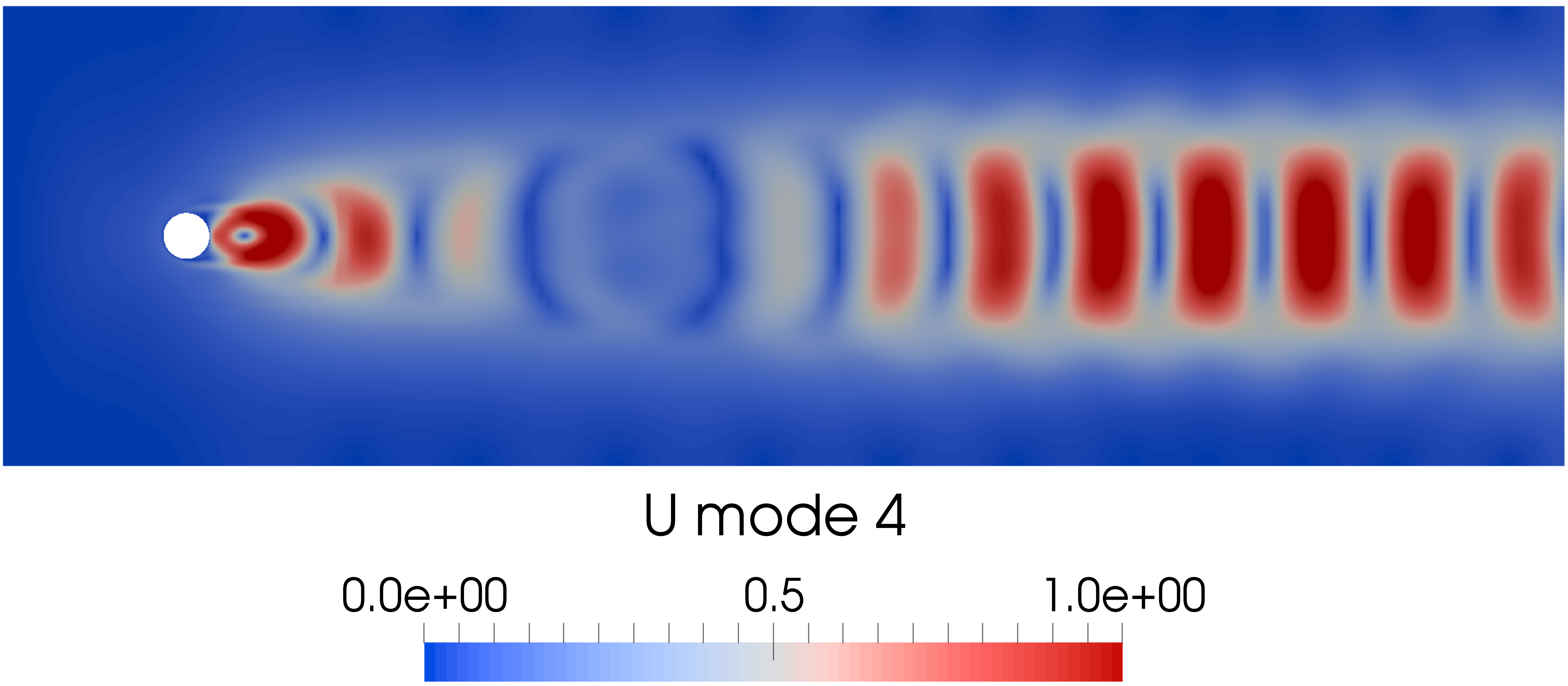}
\includegraphics[width=0.24\textwidth]{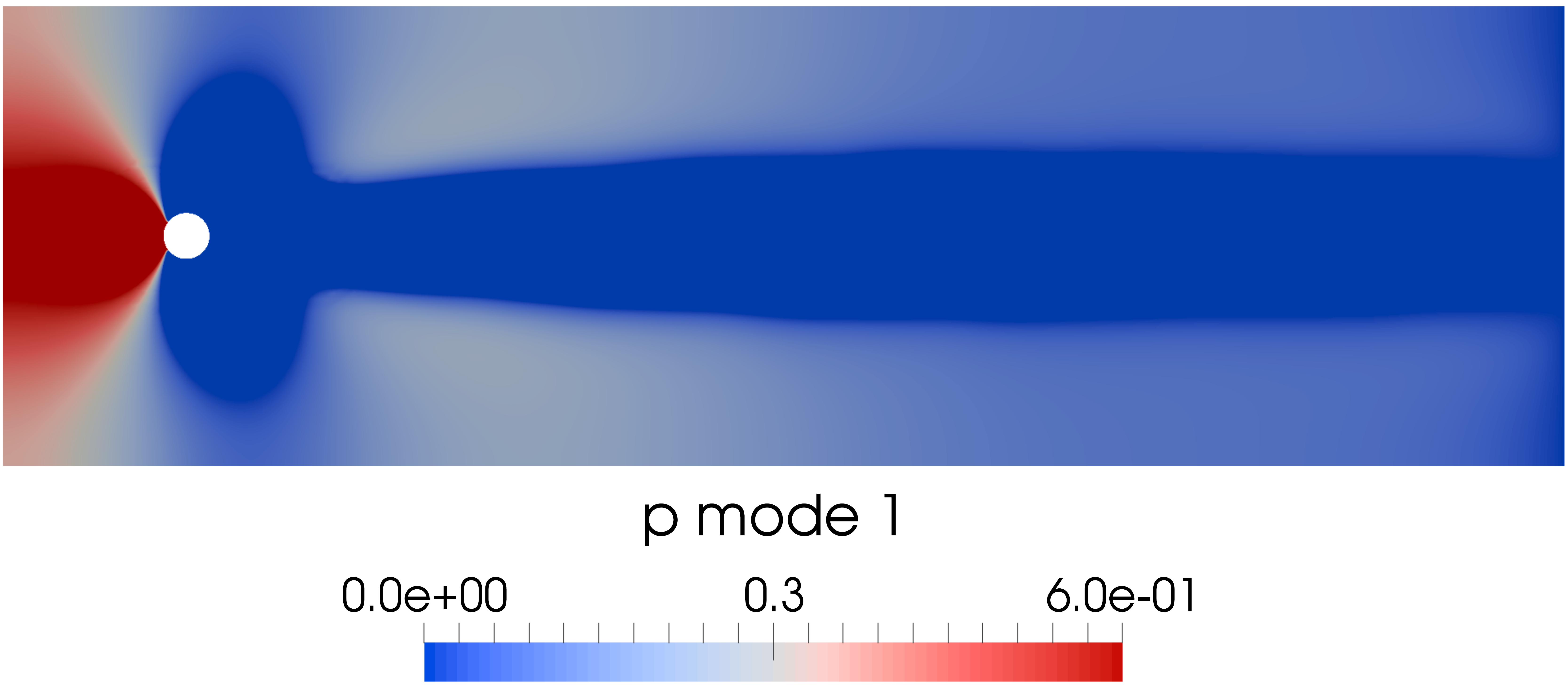}
\includegraphics[width=0.24\textwidth]{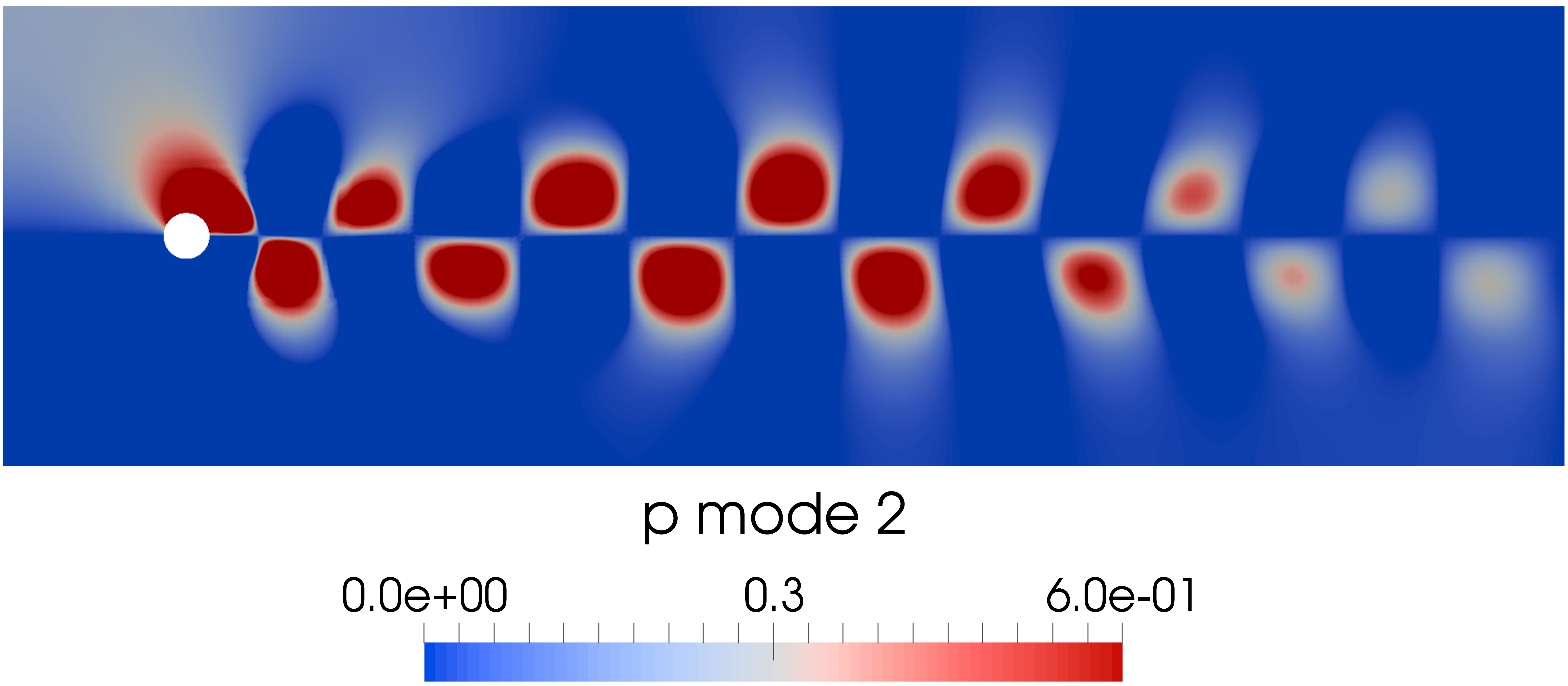}
\includegraphics[width=0.24\textwidth]{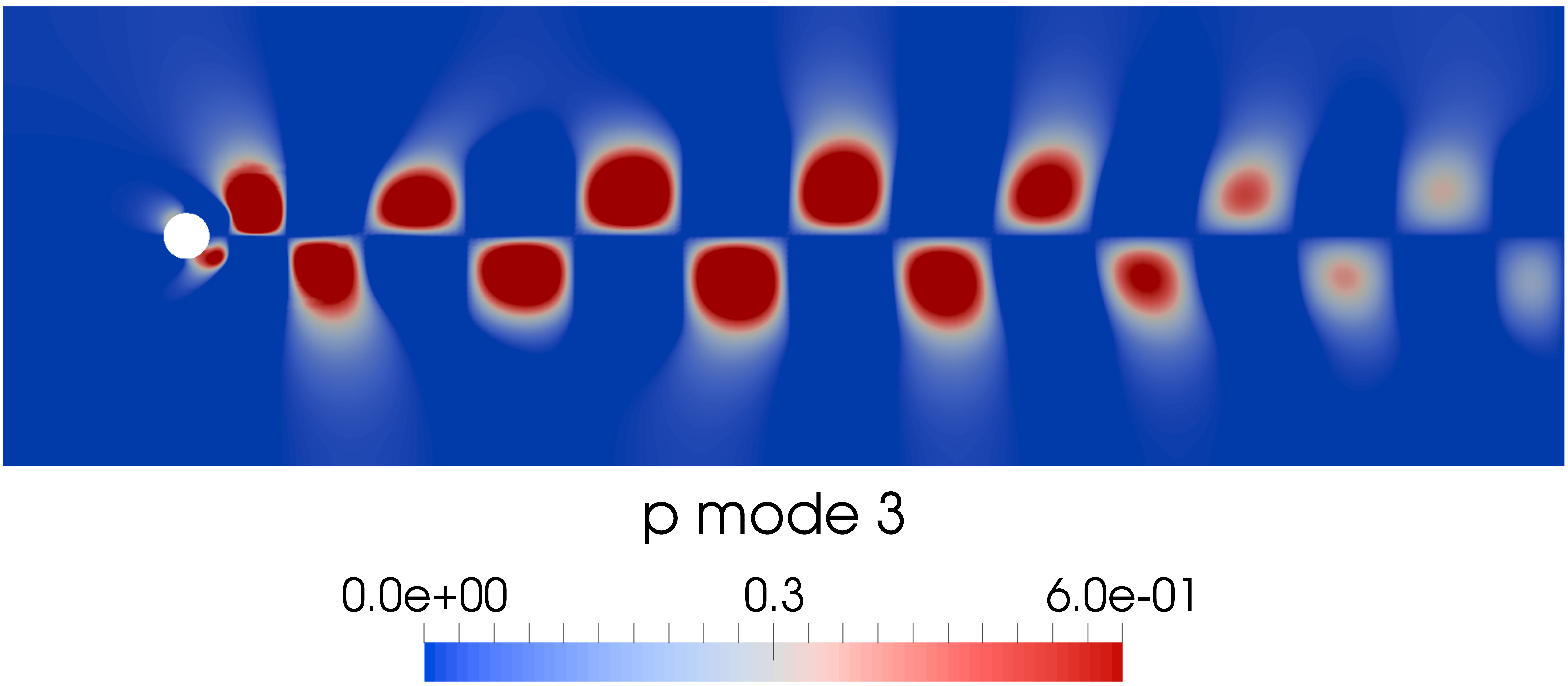}
\includegraphics[width=0.24\textwidth]{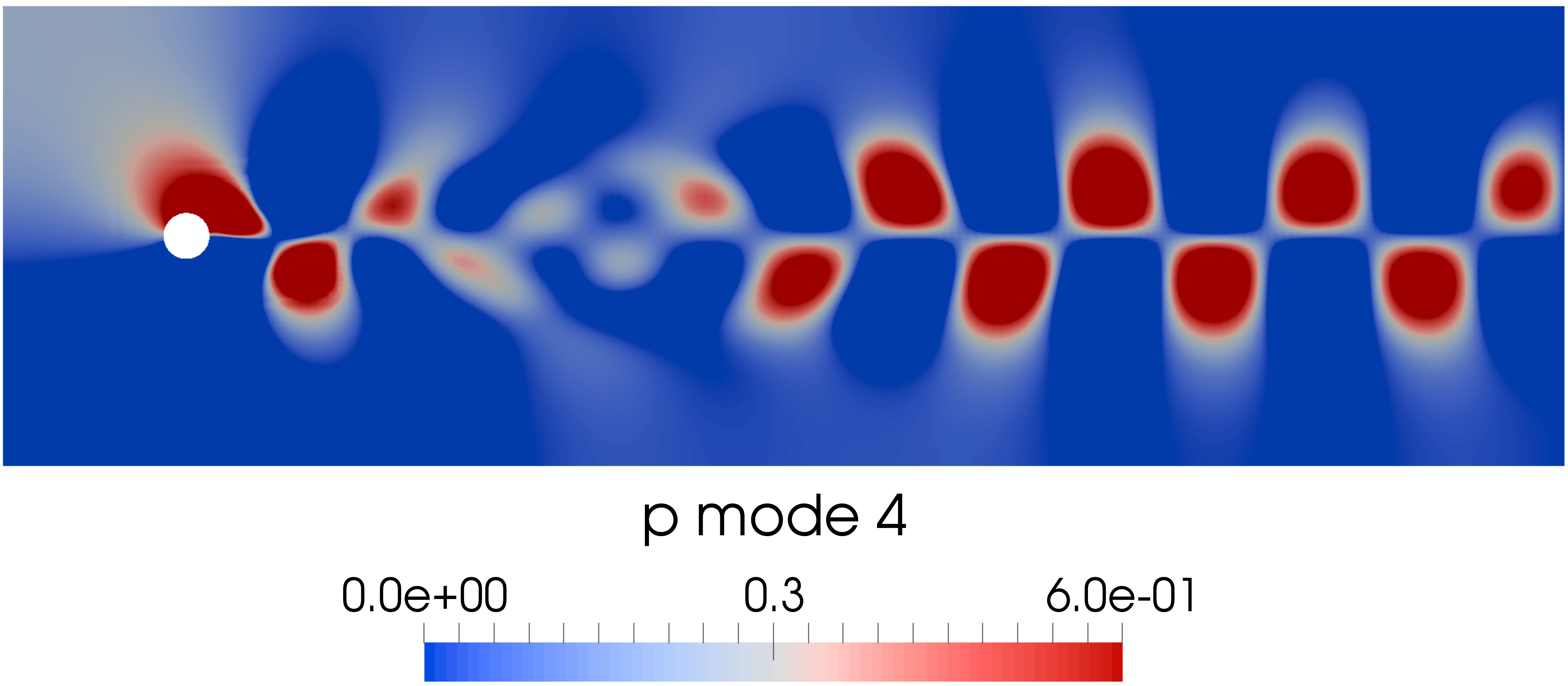} 
\includegraphics[width=0.24\textwidth]{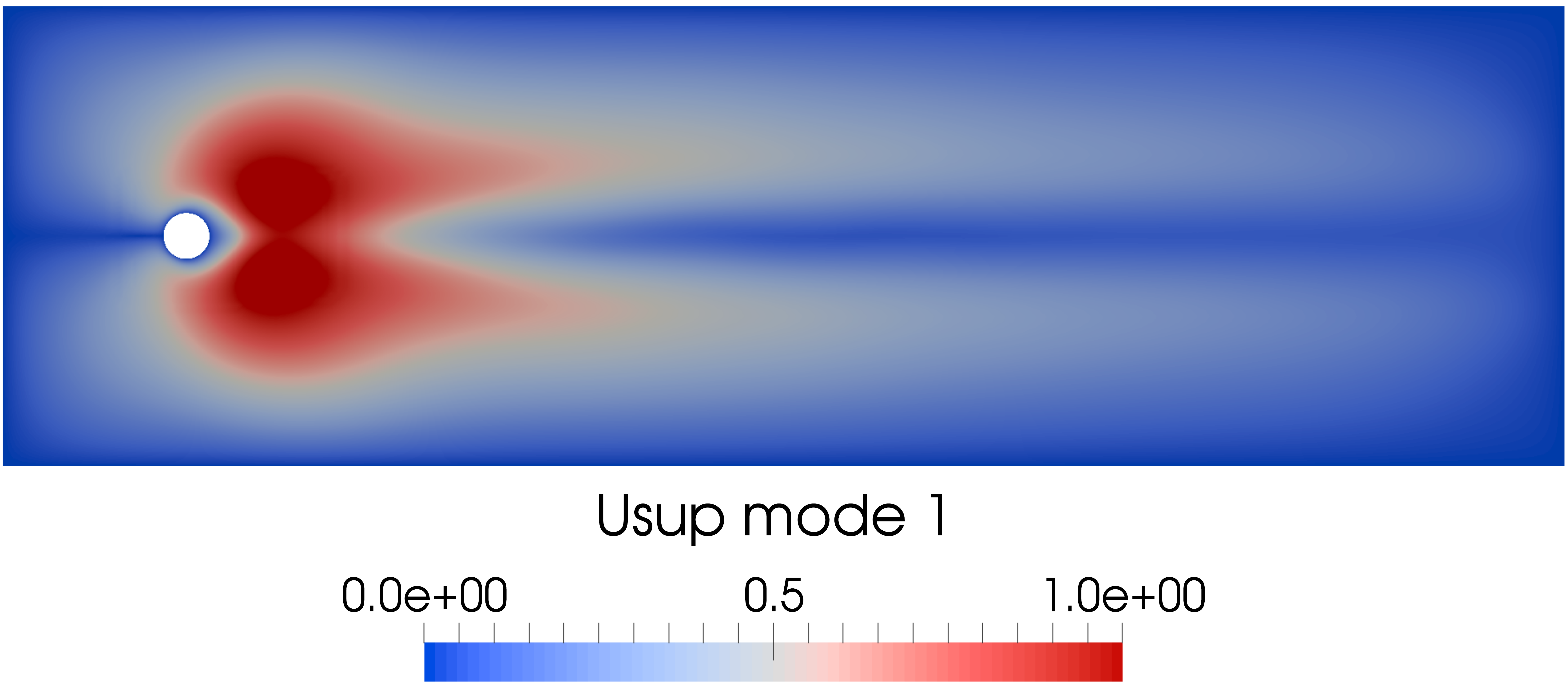}
\includegraphics[width=0.24\textwidth]{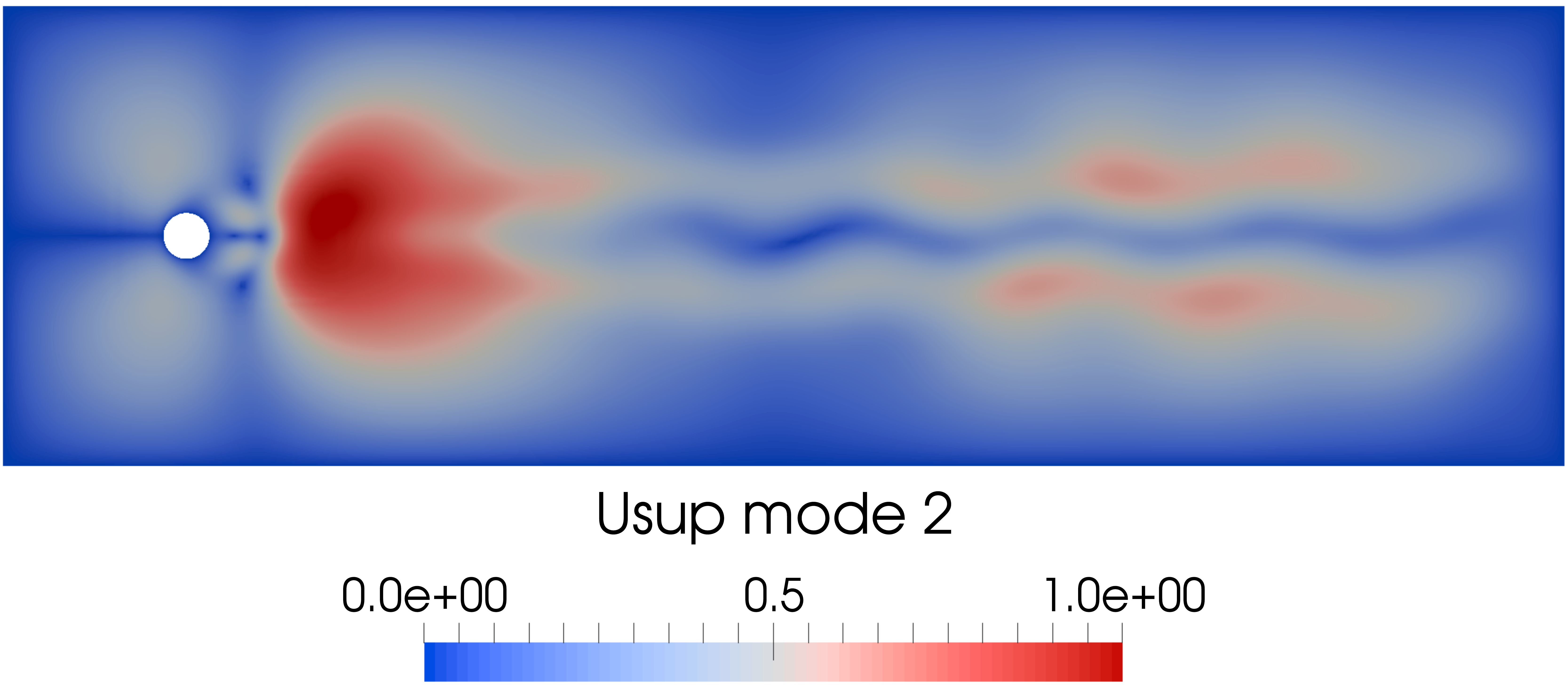}
\includegraphics[width=0.24\textwidth]{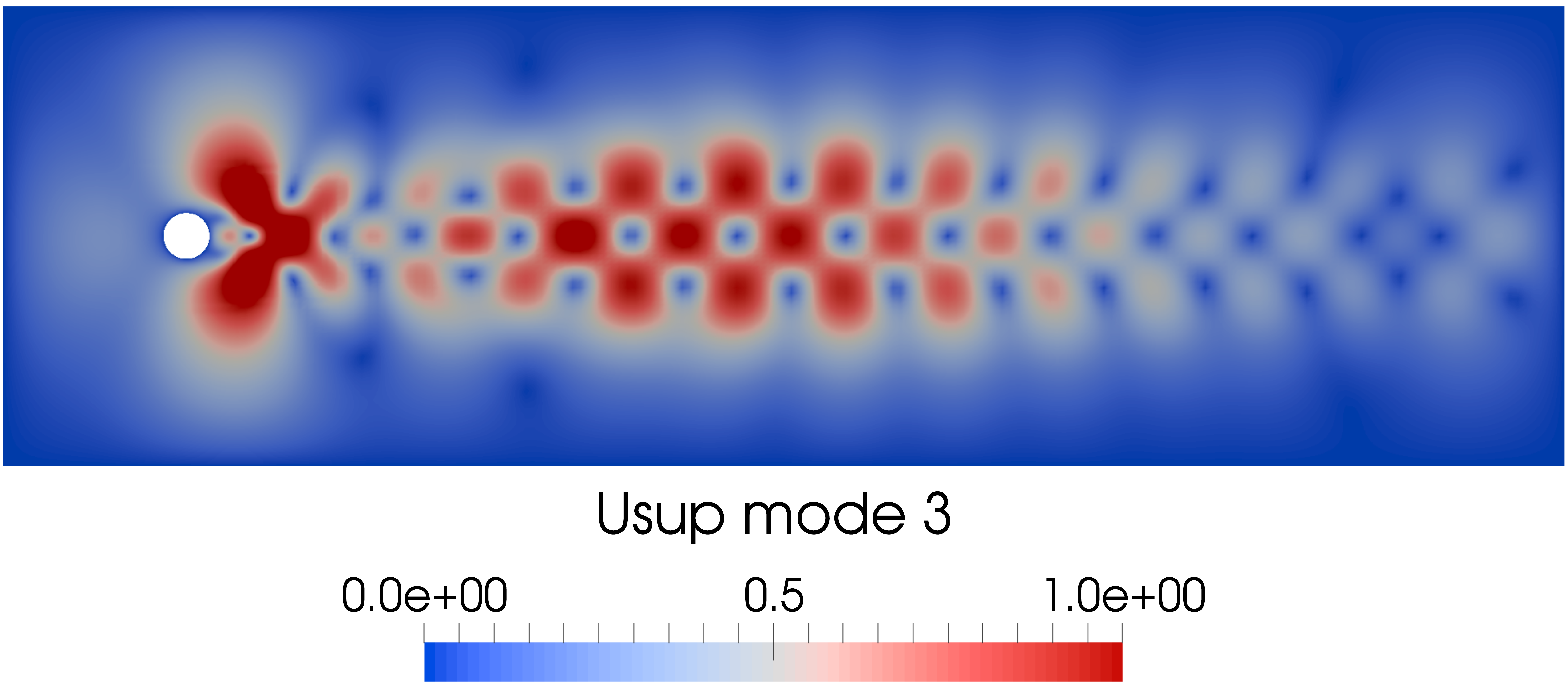}
\includegraphics[width=0.24\textwidth]{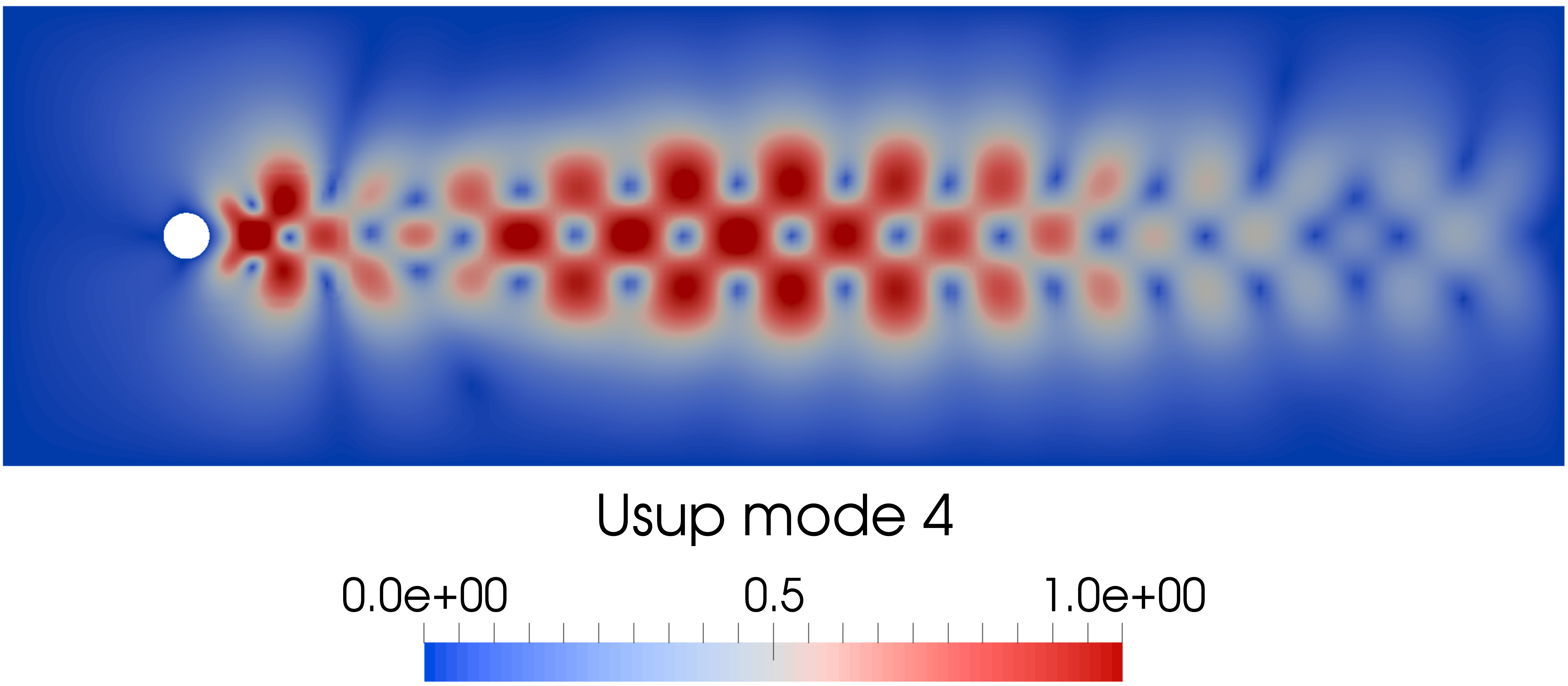}
\end{minipage} 
\caption{First four basis functions for velocity (first row), pressure (second row) and supremizers (third row).}\label{fig:basis_cyl} 
\end{figure*}

The second example, which aimed to test the methodologies on a more complex flow field and mesh structure, consists into the benchmark of the flow around a circular cylinder. In this numerical example also the physical parametrisation due to parametrised physical viscosity is introduced . Furthermore, this numerical experiment has been also used to test the performances of both stabilisation methods on periodic systems for long time integrations, wider respect to time window used to create the POD bases. The mesh, which is depicted in figure~\ref{fig:mesh_cyl} together with the boundary conditions, is mainly composed by quadrilateral cells. It is refined in the proximity of the cylinder and counts a total number of $43762$ cells. The mesh is generated starting from a structured grid with a base resolution of $200$ cells along the $x$ direction and $80$ cells along the $y$ direction and it is successively refined around the cylinder with $5$ progressive layers of refinement. The time step is set equal to $\Delta t = 0.005 \si{s}$, which is sufficiently small to meet the CFL condition in every part of the domain. Being the mesh sufficiently fine, and the time step sufficiently small, in the full-order simulation, for all the terms, including the convective term, a forth order spatial interpolation scheme is used. In this numerical experiment, the physical parametrisation given by the kinematic viscosity $\nu$ is introduced. To train the ROM $5$ different values of the kinematic viscosity are used. The values of the kinematic viscosity are chosen using an uniform distribution inside the range $\nu \in [0.005,0.01]$. These values of viscosity result into the values of the Reynolds number $\mbox{Re} \in [100,200]$. 
In this numerical experiment the simulation is run, for each value of the kinematic viscosity inside the training set, for $200\si{s}$. This time is long enough to achieve a completely evolved vortex shedding pattern. Since we were interested into the correct reproduction of the ROM during the periodic response regime, only the last $10$ seconds of simulation are used to collect the snapshots for the POD basis generation. Within this time window, the snapshots are collected every $0.05$ seconds, returning a total number of $200$ snapshots for each different value of the kinematic viscosity. The snapshots of the five different full-order simulations are then used to create the POD basis functions, which result in the cumulative eigenvalues of table~\ref{tab:cum_eig_cyl}. Figure~\ref{fig:basis_cyl} depicts the first $4$ basis functions for velocity, pressure and supremizers. The reduced order model counts, for the SUP-ROM, $15$ modes for velocity, $10$ modes for pressure and $12$ modes for supremizers, while for the PPE-ROM, $15$ modes for velocity, $10$ modes for pressure. 

It is worth remarking that in this example, in the SUP-ROM, the supremizer space counts more modes with respect to the pressure space. This choice is done in order to improve the accuracy and the stability of the results. For this particular case, in fact, we have experimentally observed, that an equal number of pressure and supremizer modes, leads to inaccurate results. This is justified by the fact that using the approximated approach described in \autoref{subsec:sup_enrich} for the supremizer enrichment, an equal dimension of the pressure and the supremizer spaces, does not automatically guaranty the fulfilment of the inf-sup condition. 

\begin{figure*}[t]
\centering
\includegraphics[width=0.48\textwidth]{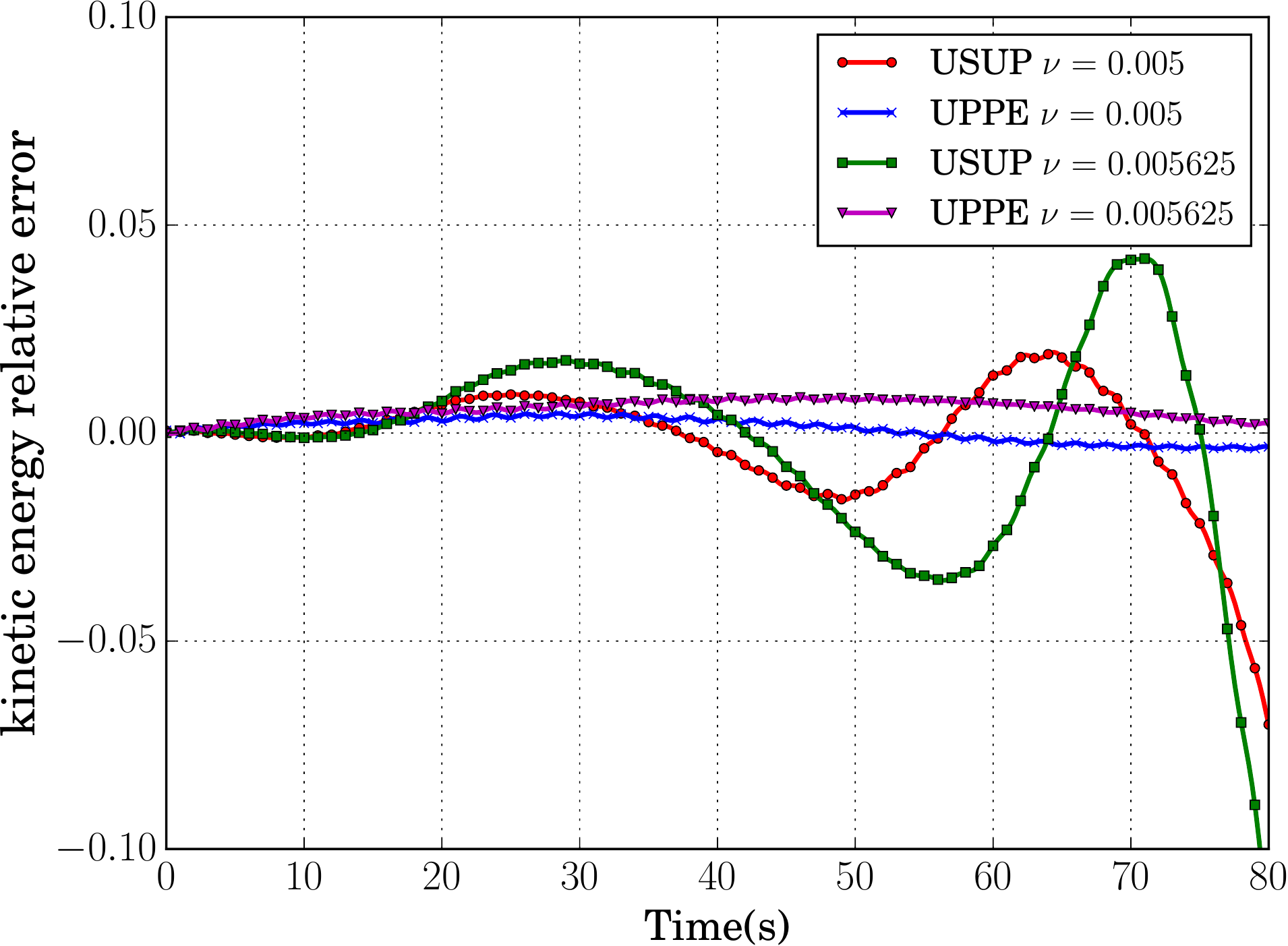}
\caption{Kinetic energy relative error in the cylinder example for $\nu = 0.005$ and $\nu = 0.005625$. The kinetic energy relative error is plotted over time for the two values of viscosity and for the two different models: with supremizer stabilisation (USUP and PSUP) and pressure Poisson equation stabilisation (UPPE and PPPE). The ROM solutions are obtained with 15 modes for velocity, 10 modes for pressure and 12 modes for supremizers. The time window in this case is wider respect to the one used for the generation of the reduced basis spaces ($\Delta T = 10$s)}\label{fig:kyn_cyl} 
\end{figure*}

\begin{figure*}
\begin{minipage}{1\textwidth}
\centering
\includegraphics[width=0.24\textwidth]{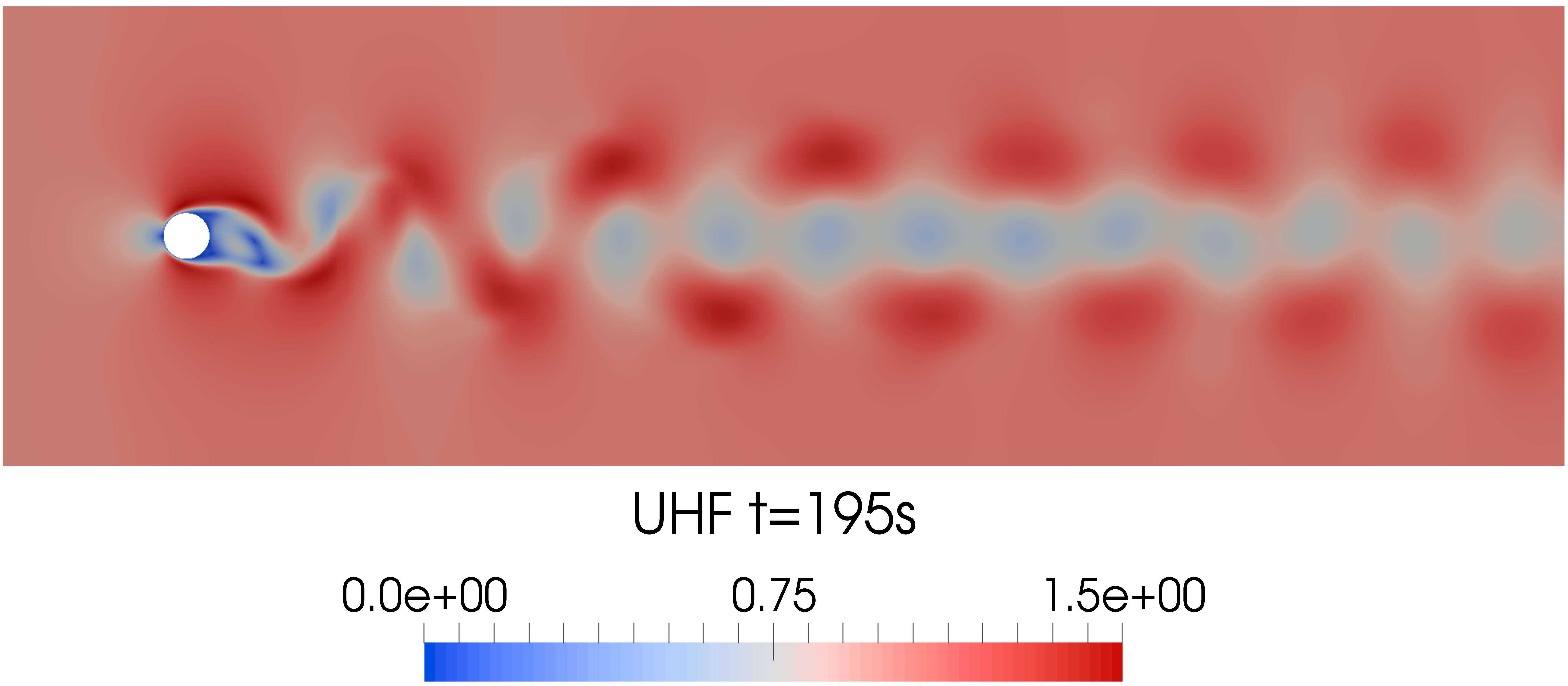}
\includegraphics[width=0.24\textwidth]{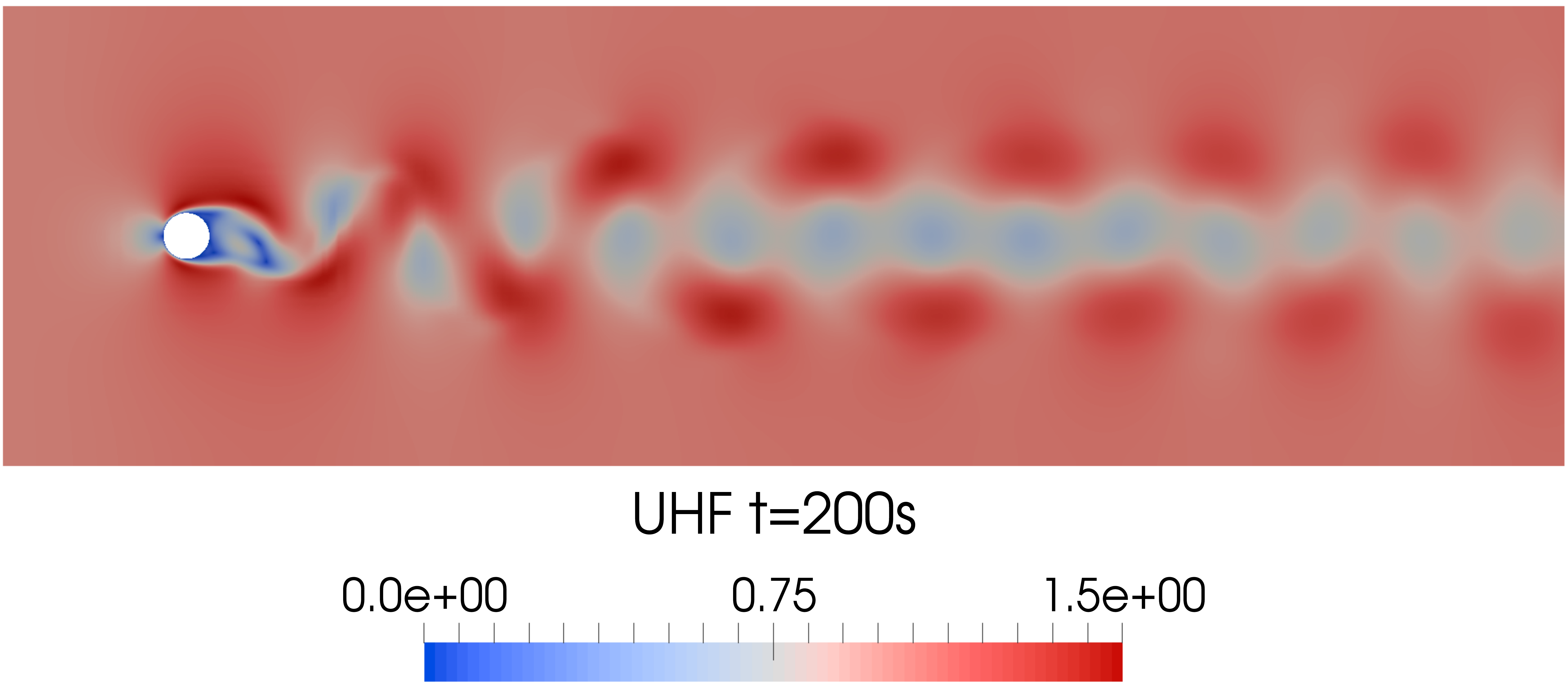}
\includegraphics[width=0.24\textwidth]{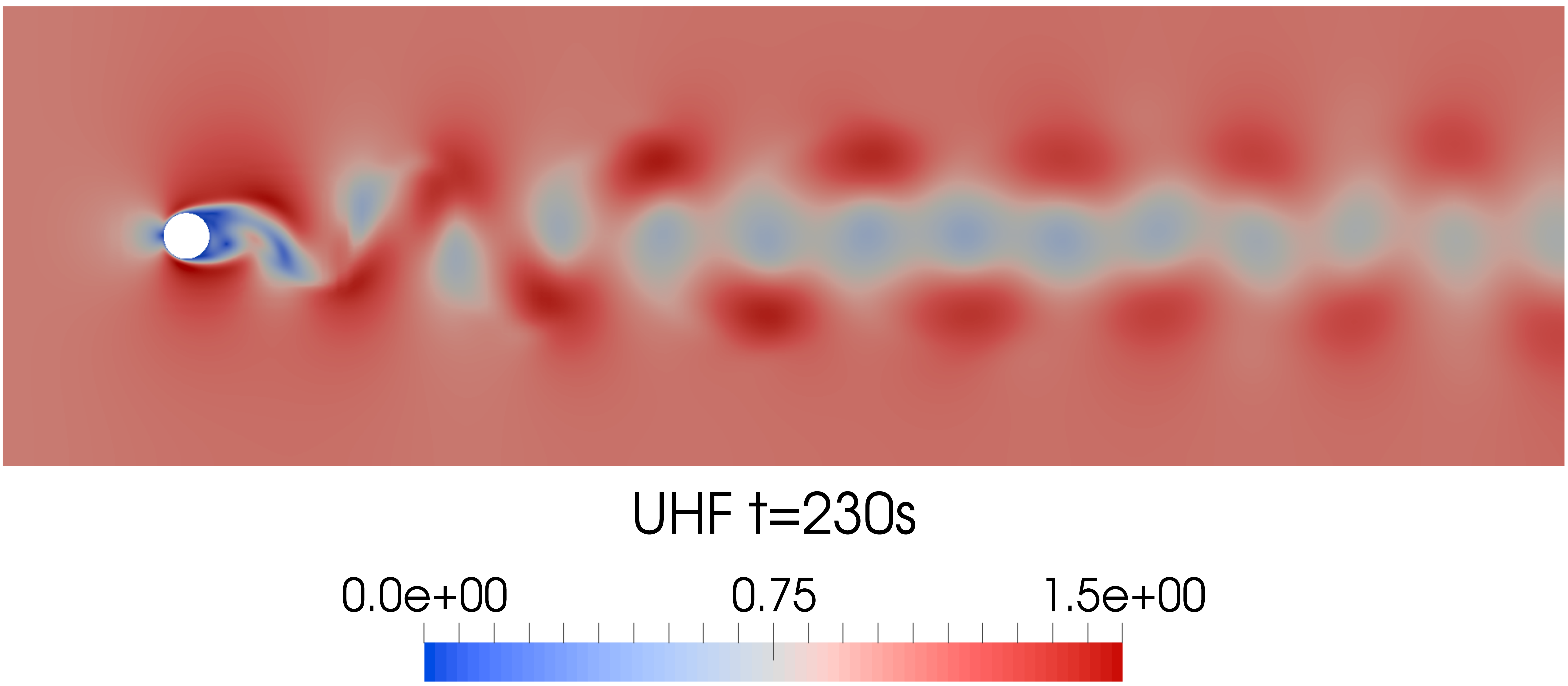}
\includegraphics[width=0.24\textwidth]{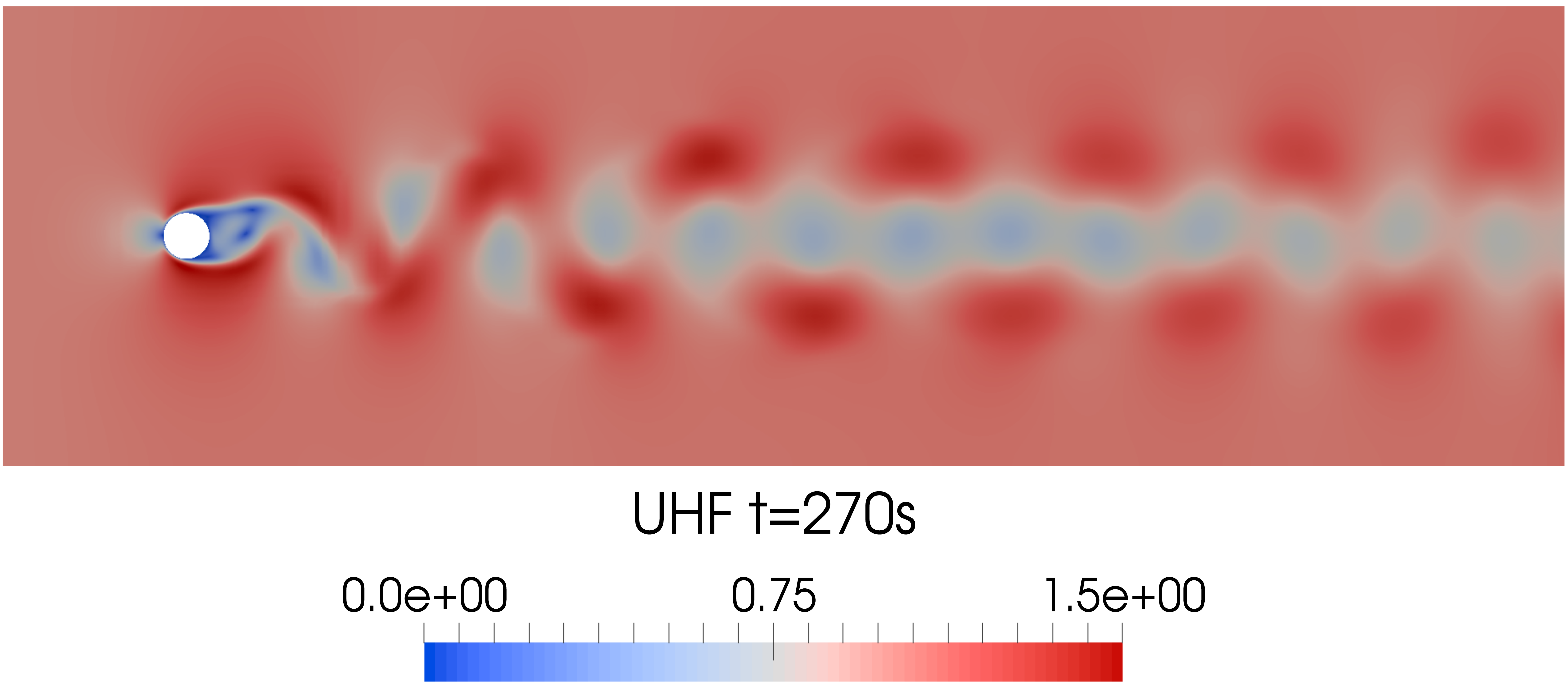}
\includegraphics[width=0.24\textwidth]{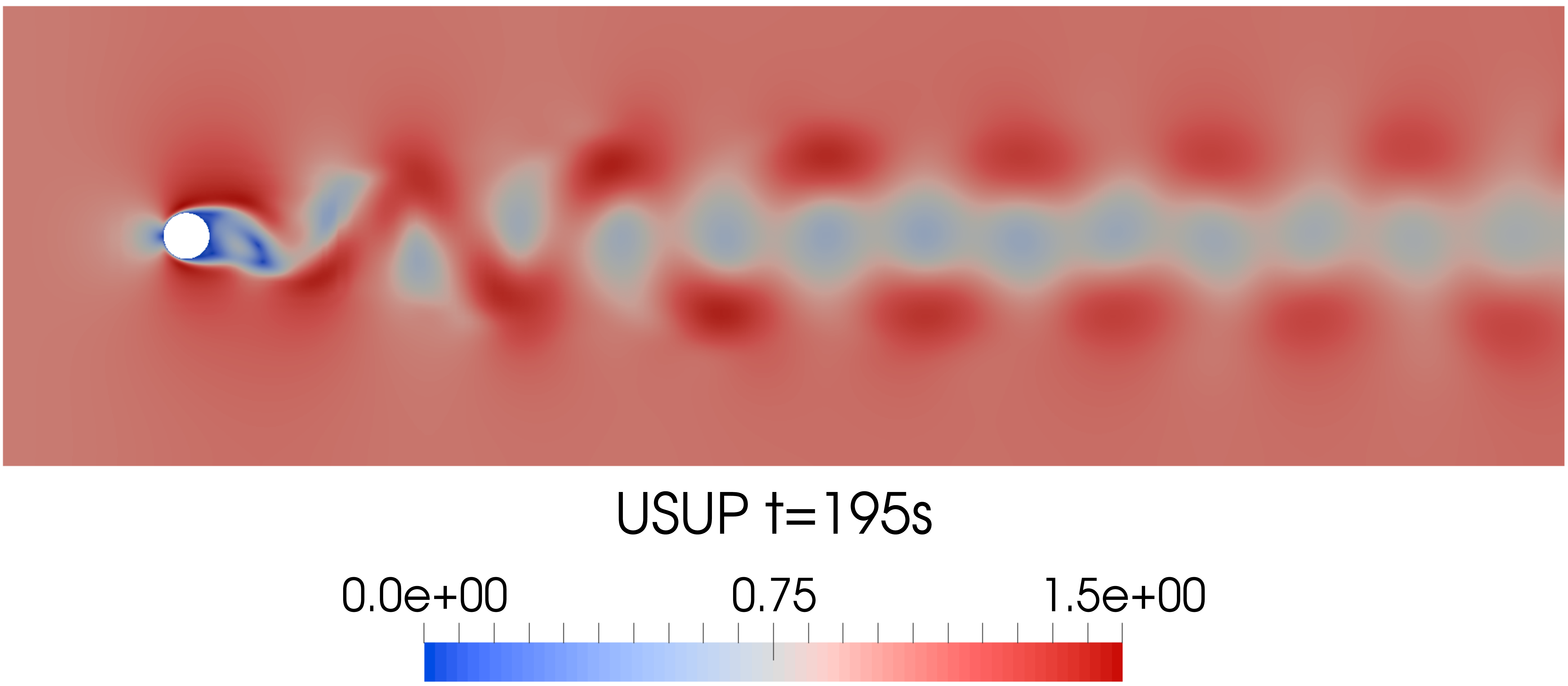}
\includegraphics[width=0.24\textwidth]{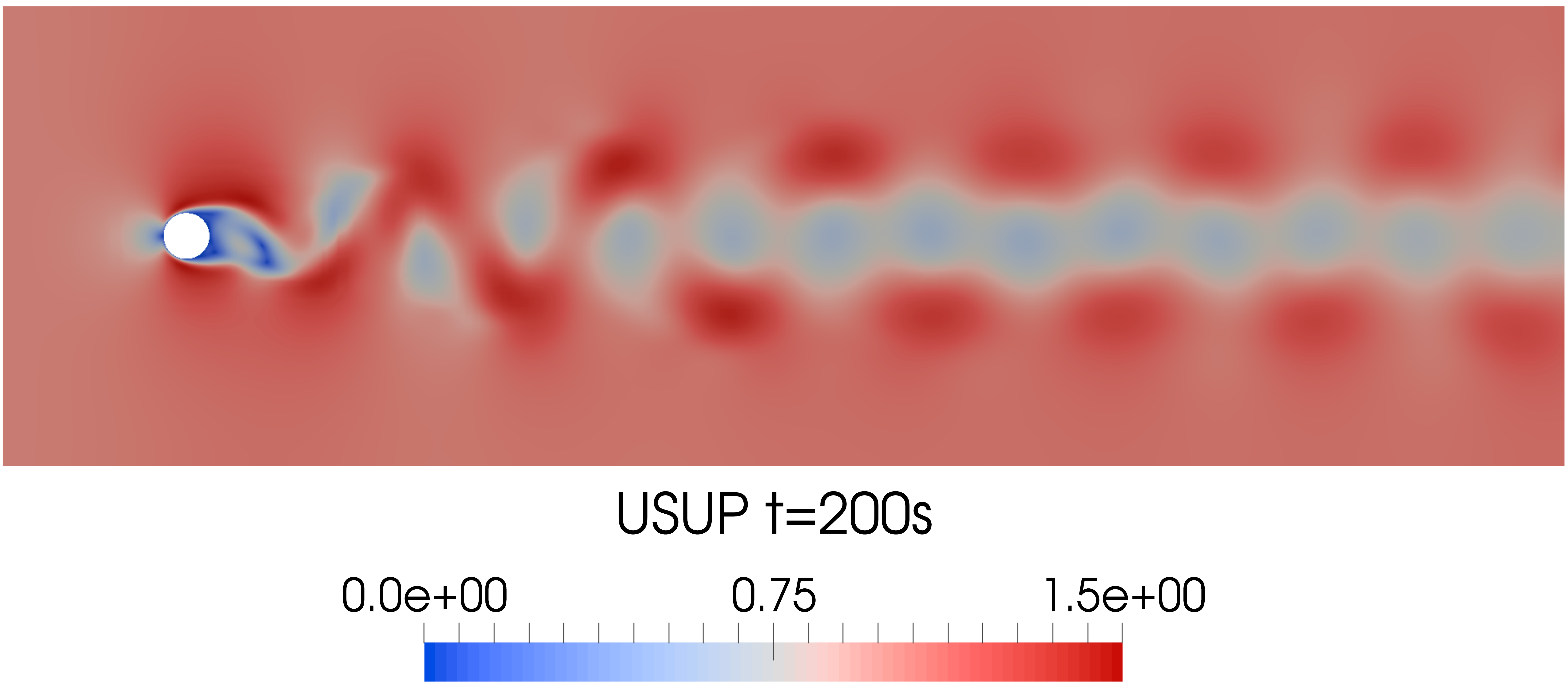}
\includegraphics[width=0.24\textwidth]{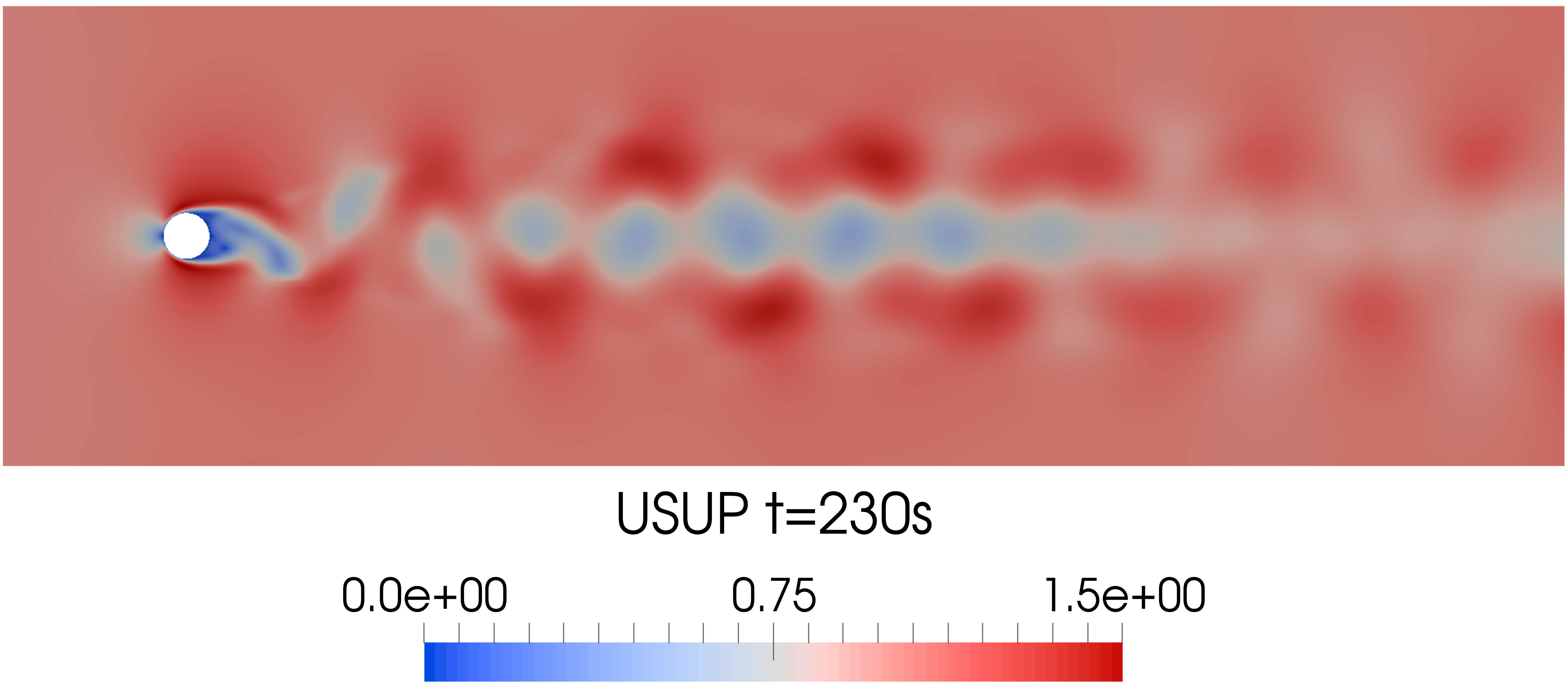}
\includegraphics[width=0.24\textwidth]{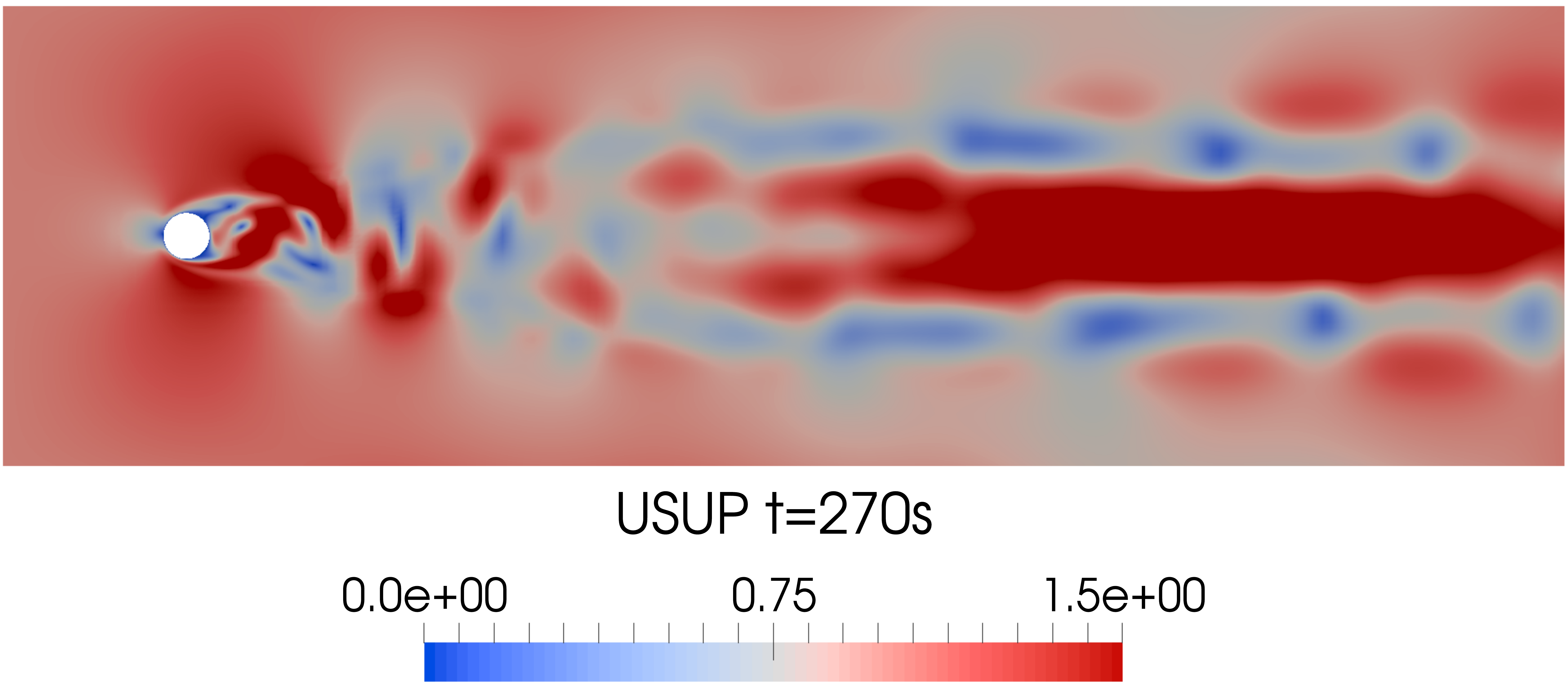} 
\includegraphics[width=0.24\textwidth]{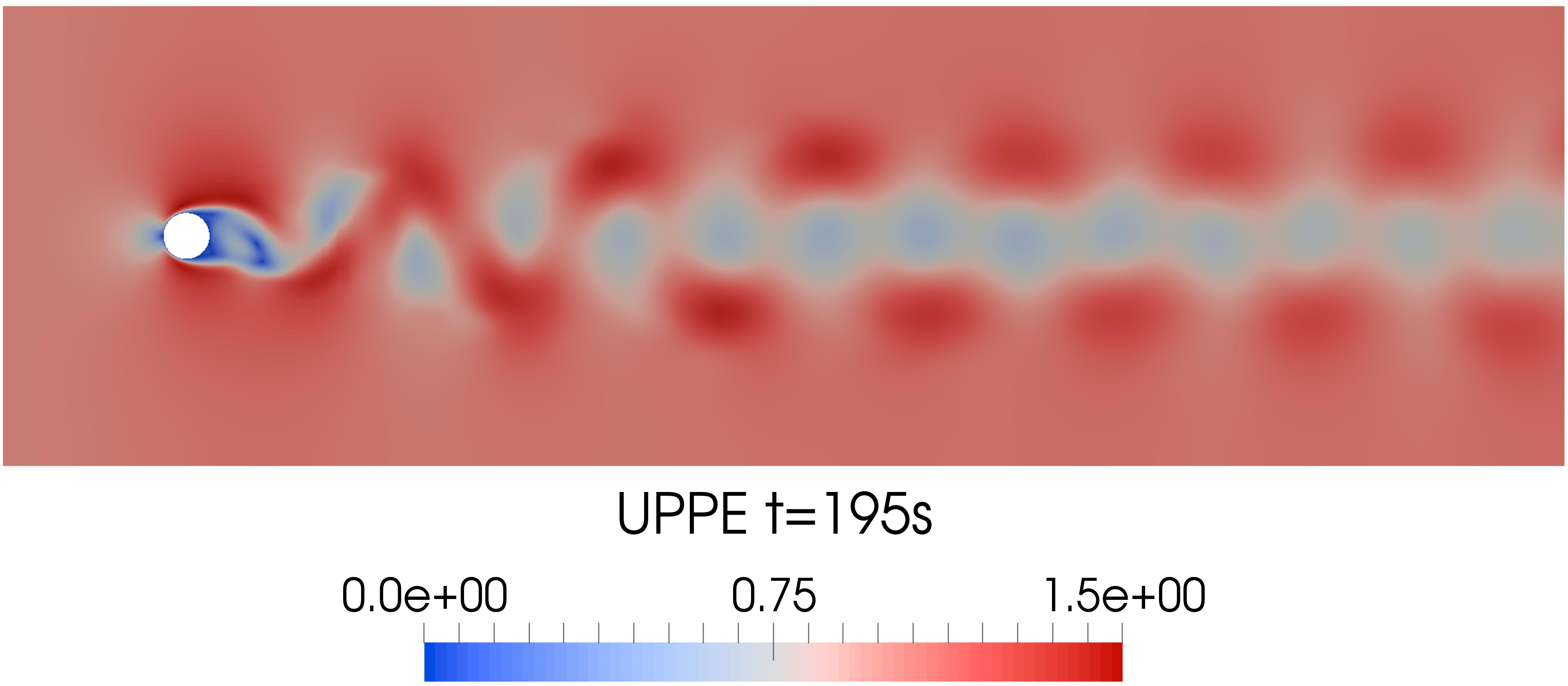}
\includegraphics[width=0.24\textwidth]{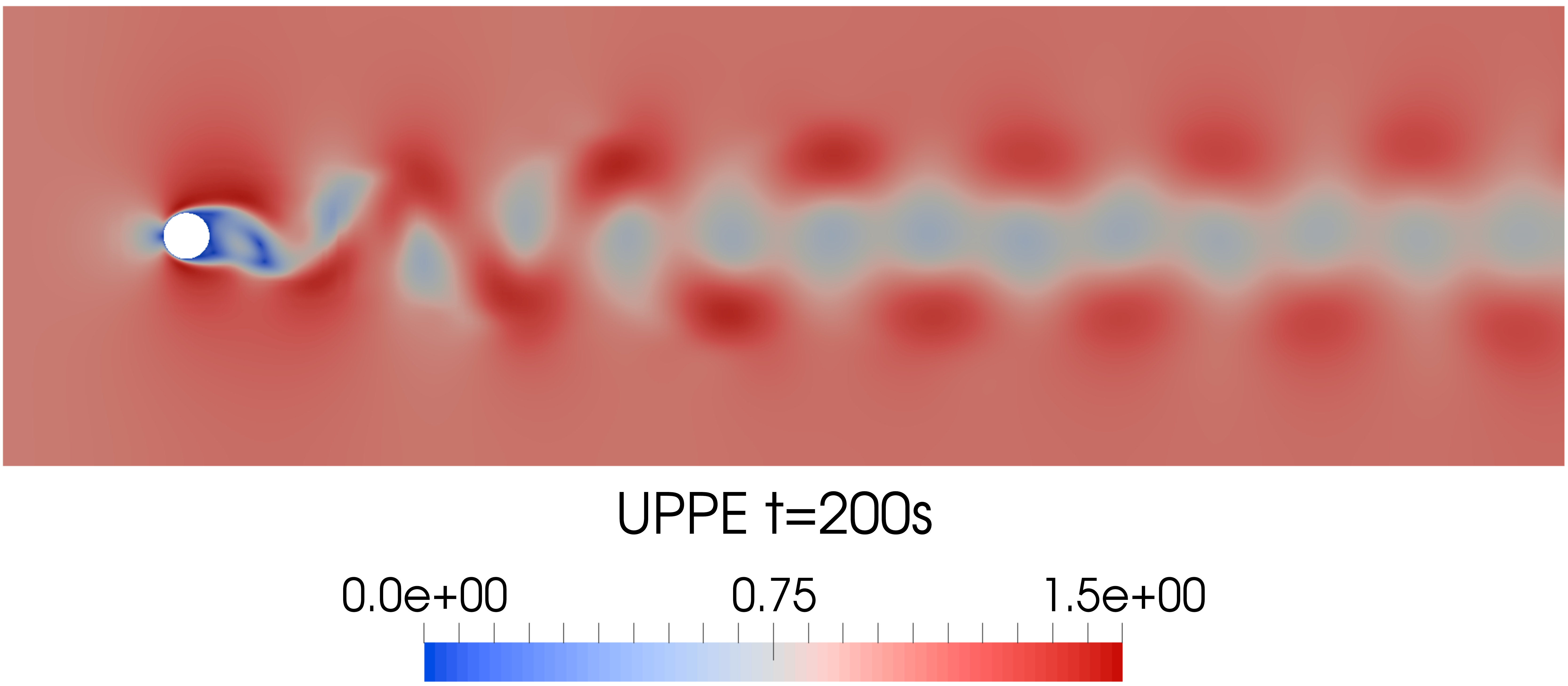}
\includegraphics[width=0.24\textwidth]{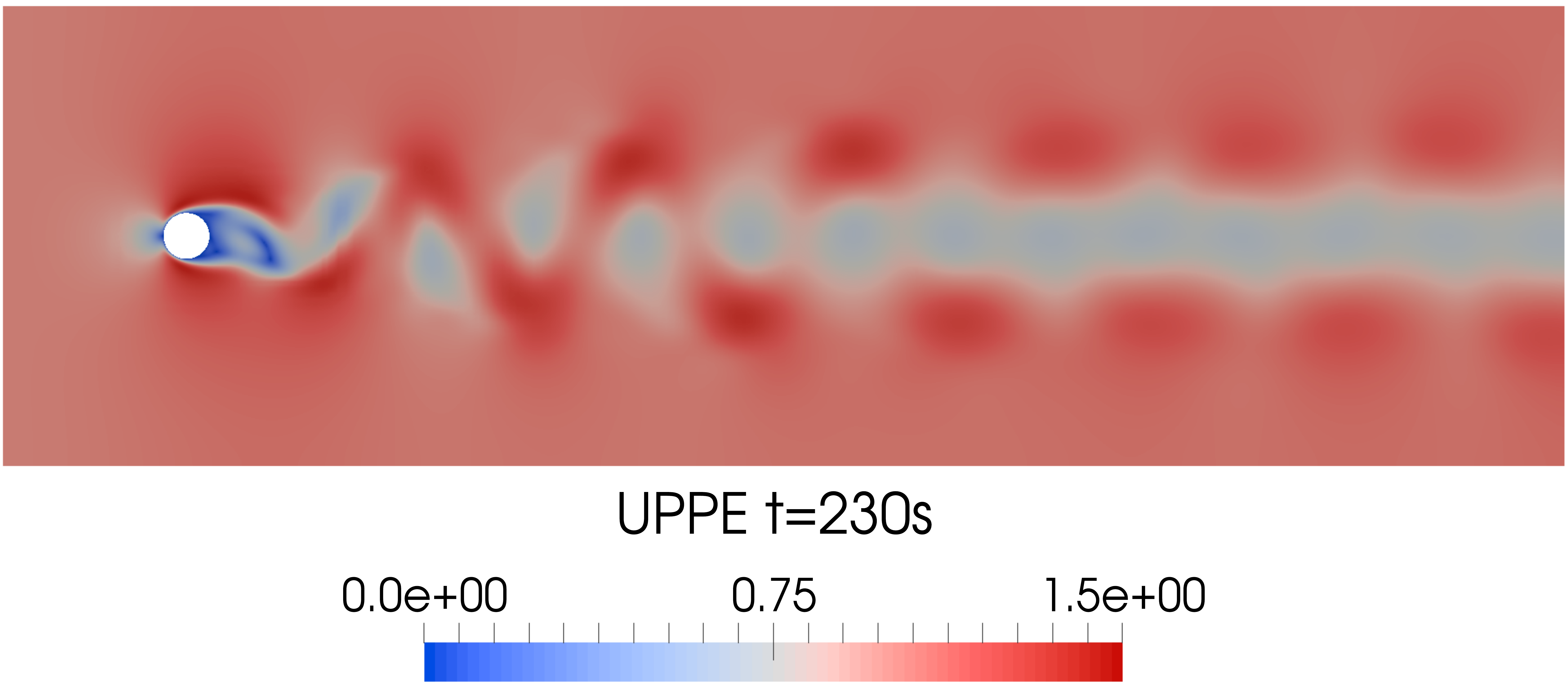}
\includegraphics[width=0.24\textwidth]{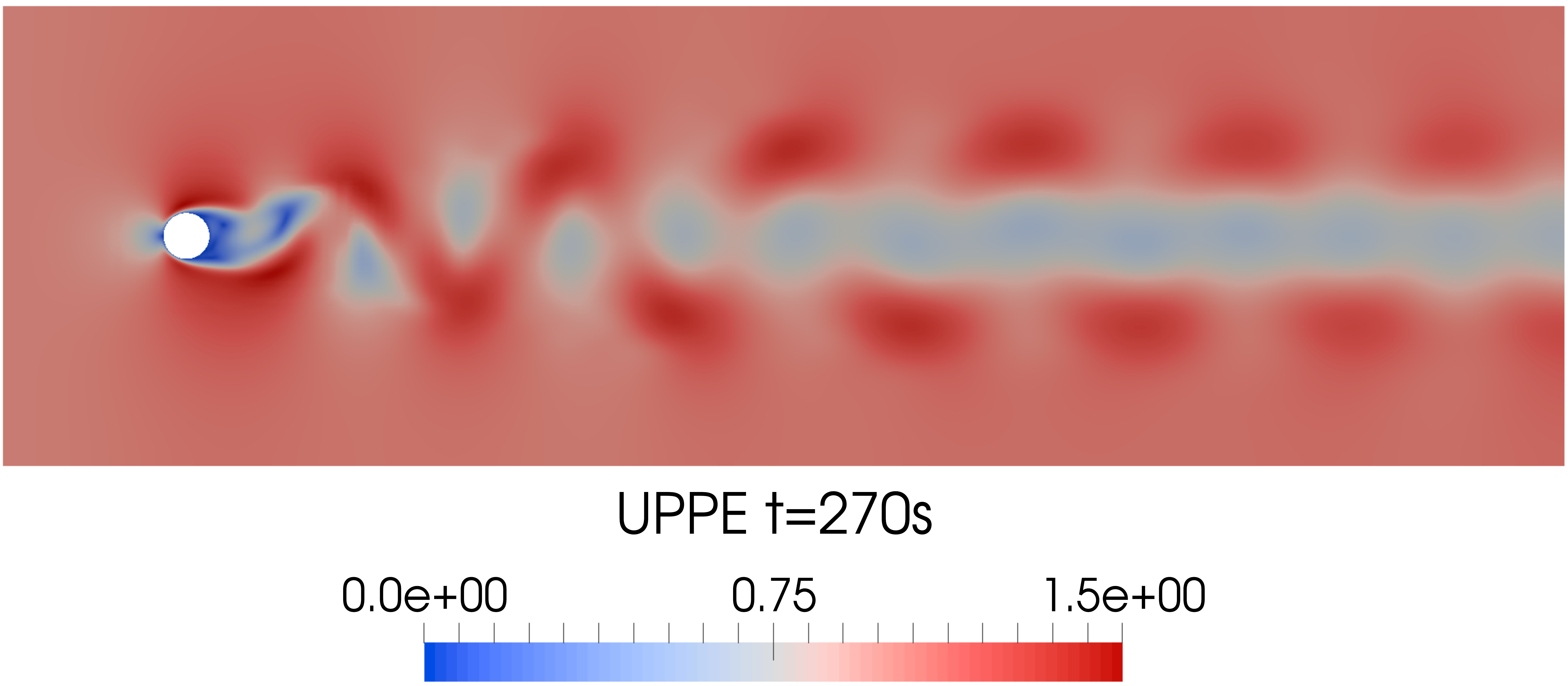}
\end{minipage} 
\caption{Comparison of the velocity field for high fidelity (UHF - first row), supremizer stabilised ROM (USUP - second row) and pressure Poisson equation stabilised ROM (UPPE third row). The fields are depicted for different time instant equal to $t=195 \si{s},200 \si{s},230 \si{s}$ and $270 \si{s}$ and increasing from left to right. The ROM solutions are obtained with $15$ modes for velocity and $10$ modes for pressure, and only for the SUP-ROM with $12$ additional supremizer modes. The velocity magnitude is shown in the images legends.}\label{fig:comparison_cavity}\label{fig:velocity_cyl} 
\end{figure*}

\begin{figure*}
\begin{minipage}{1\textwidth}
\centering
\includegraphics[width=0.24\textwidth]{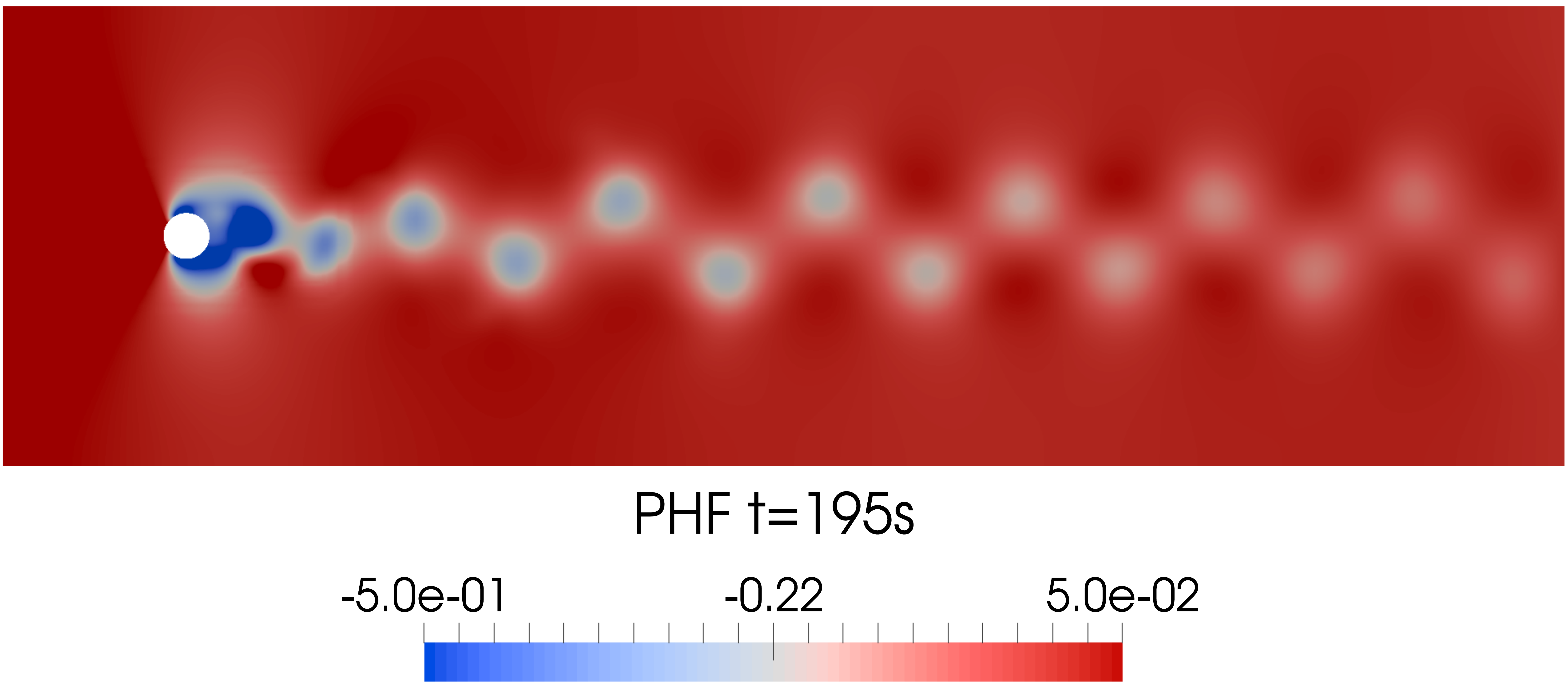}
\includegraphics[width=0.24\textwidth]{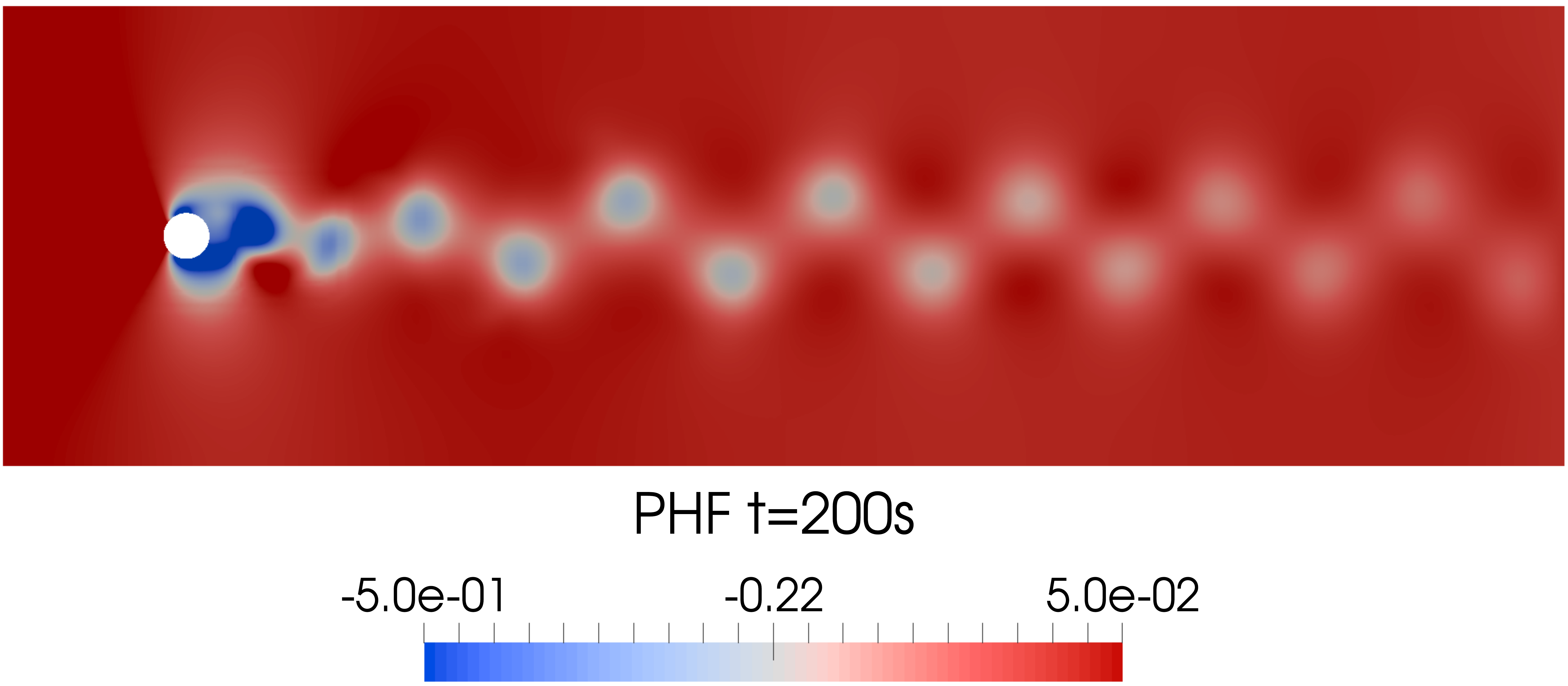}
\includegraphics[width=0.24\textwidth]{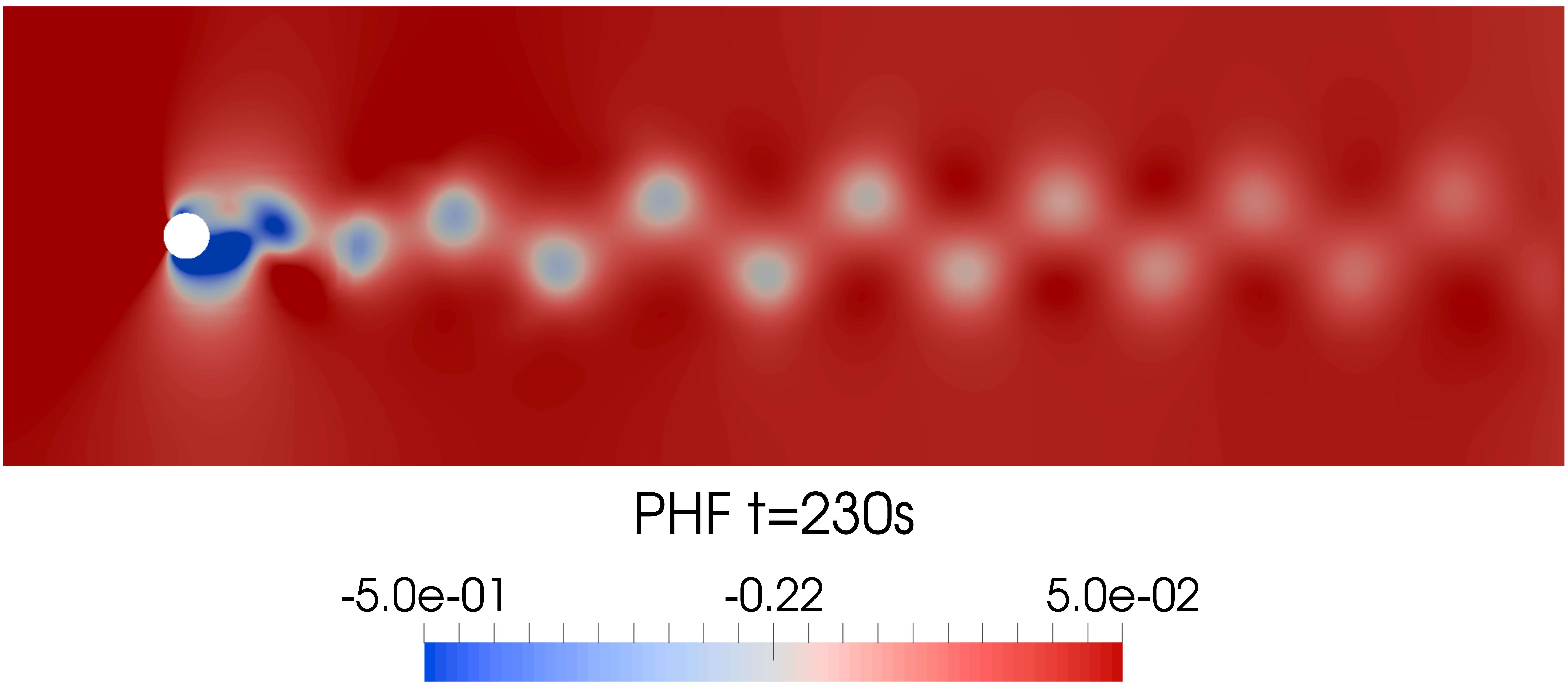}
\includegraphics[width=0.24\textwidth]{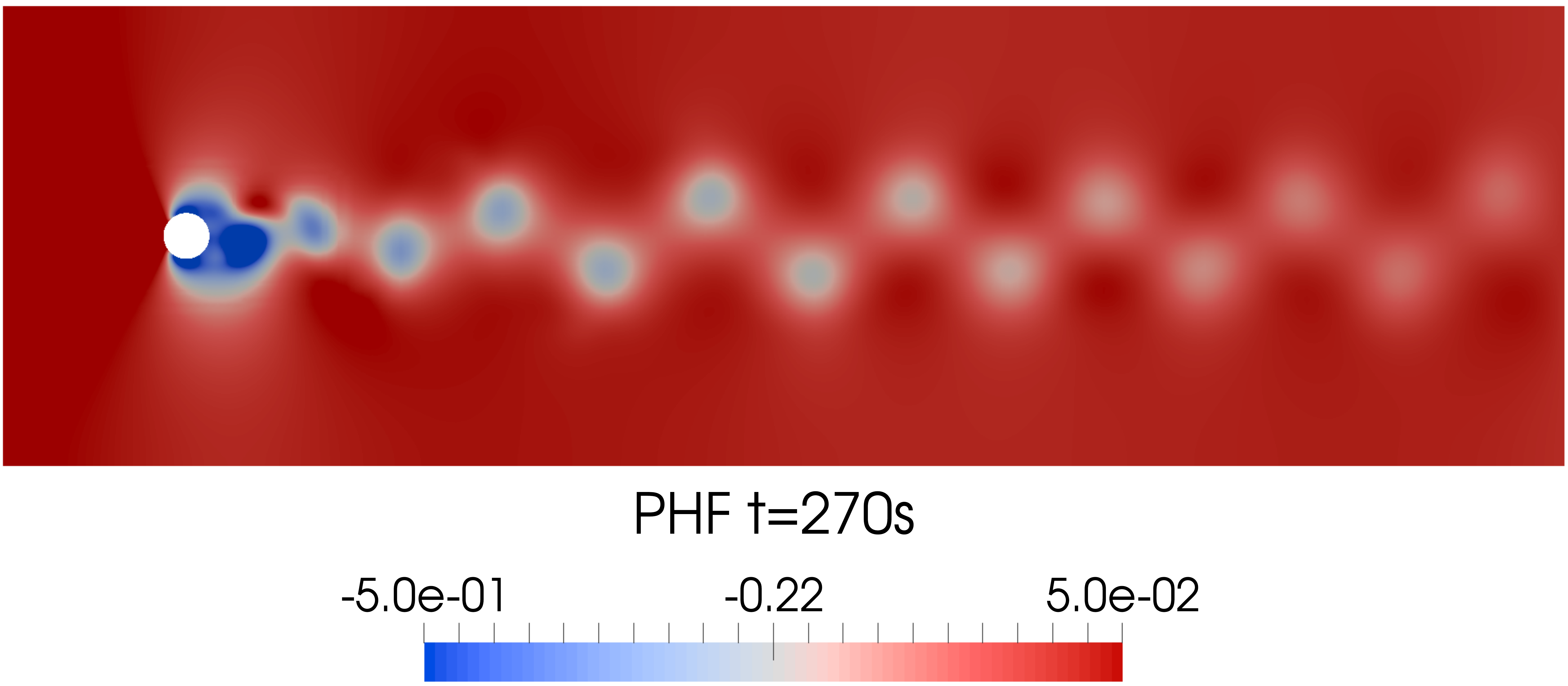}
\includegraphics[width=0.24\textwidth]{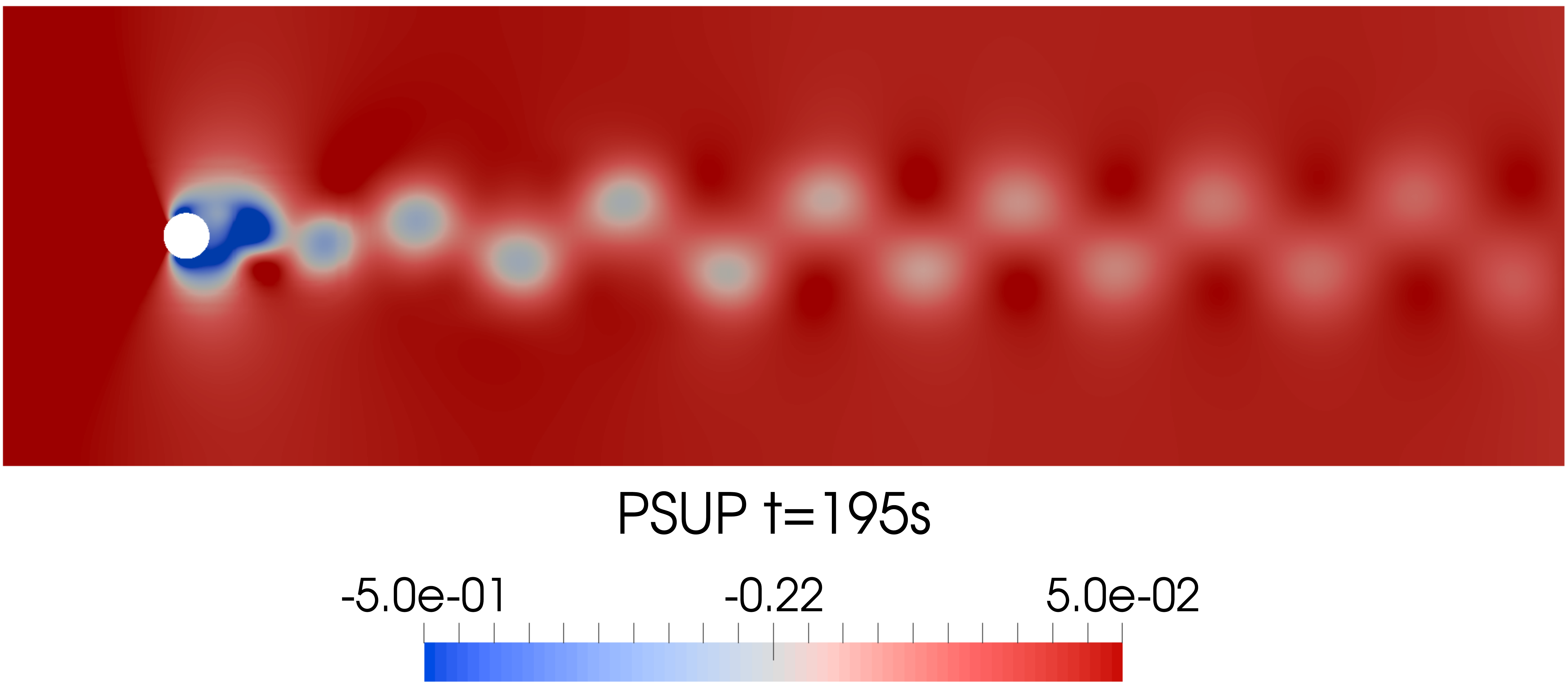}
\includegraphics[width=0.24\textwidth]{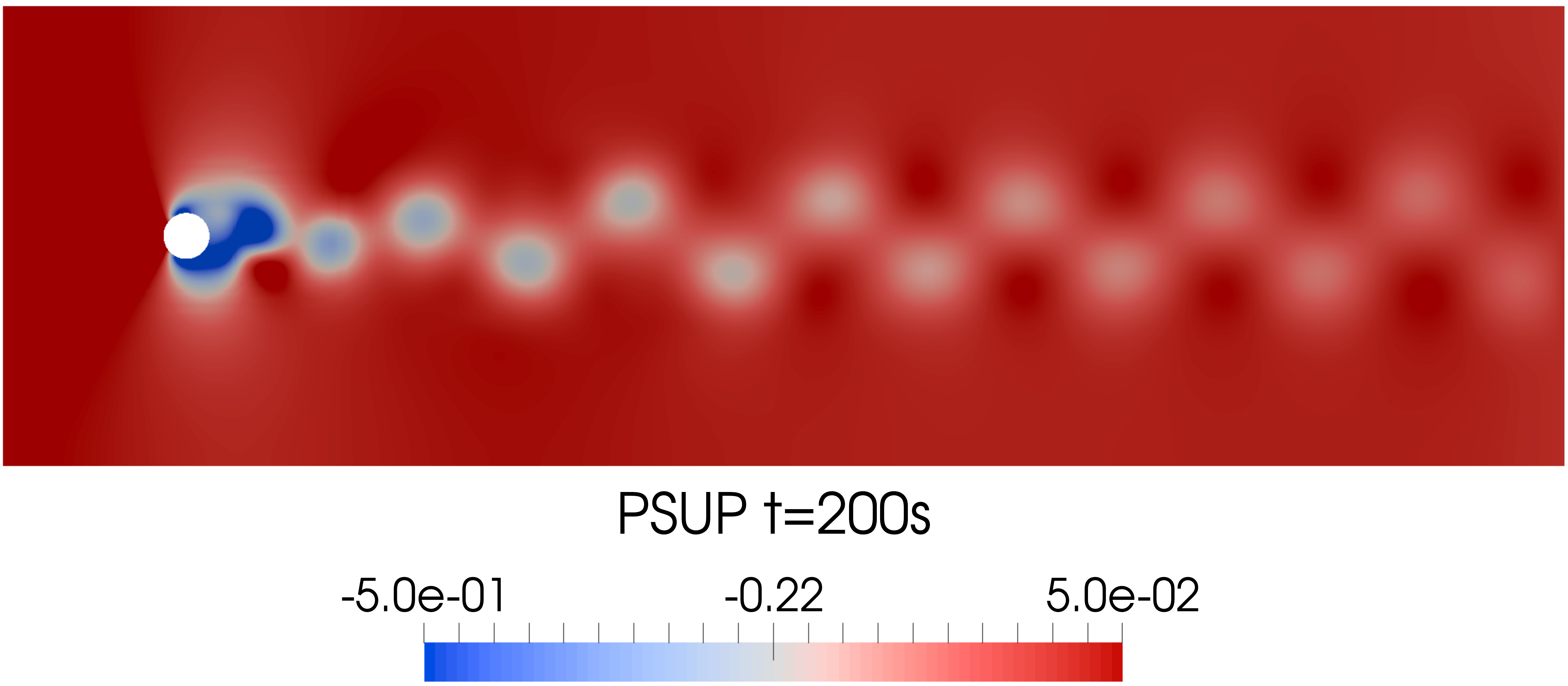}
\includegraphics[width=0.24\textwidth]{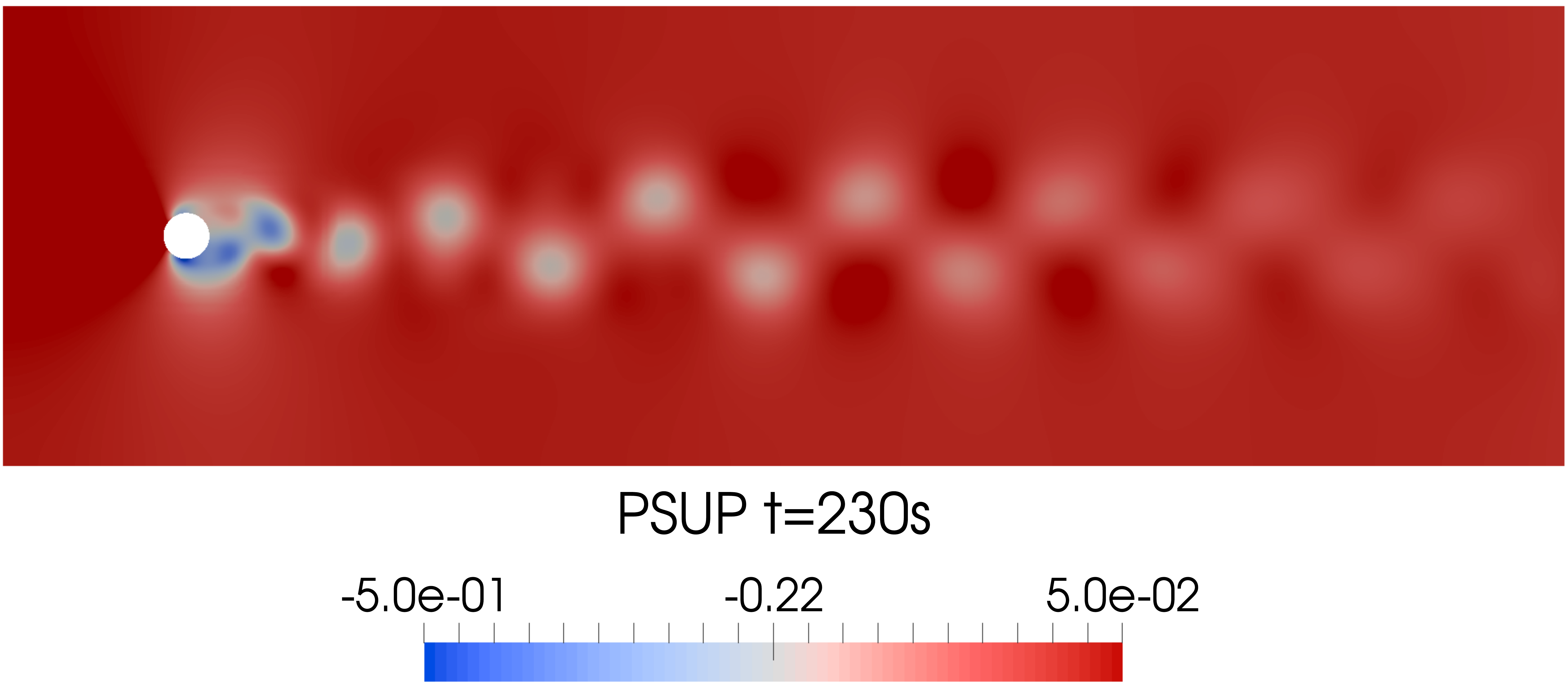}
\includegraphics[width=0.24\textwidth]{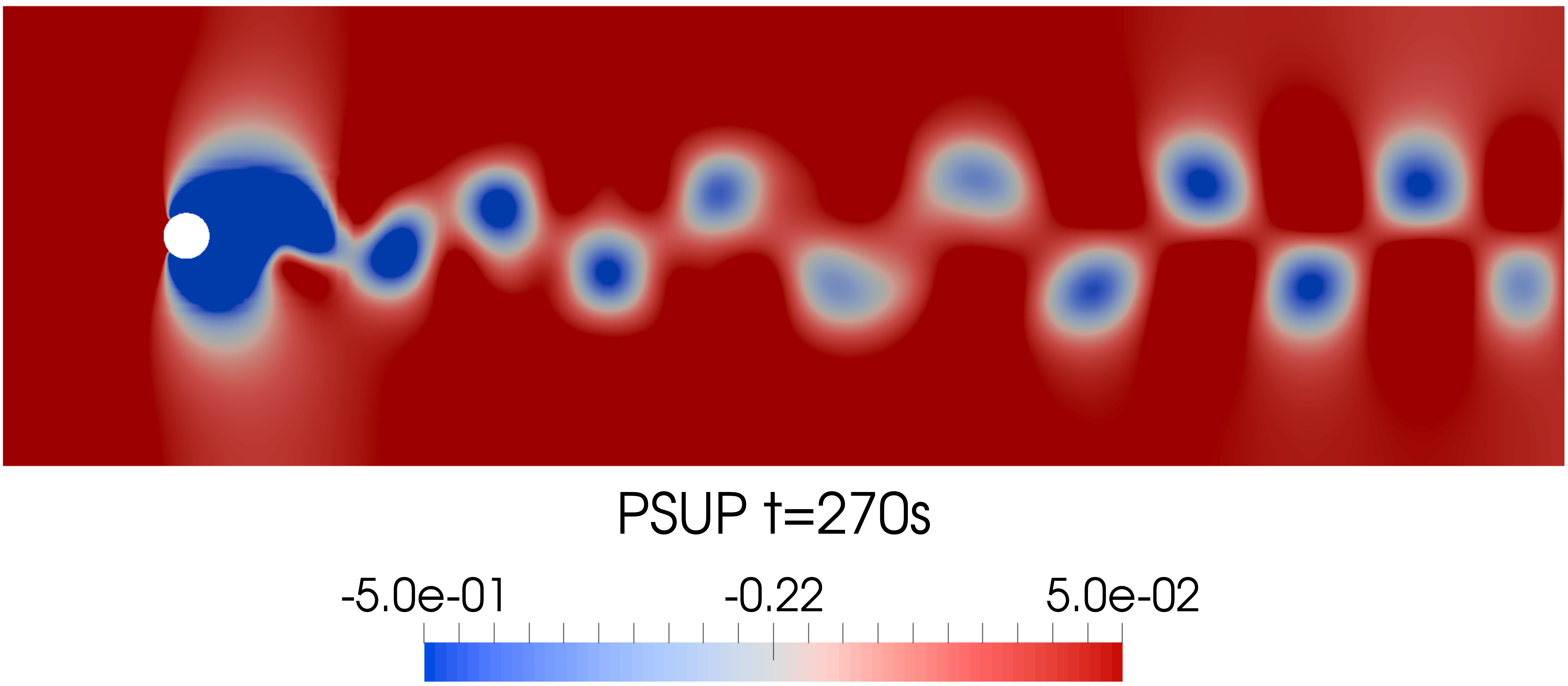} 
\includegraphics[width=0.24\textwidth]{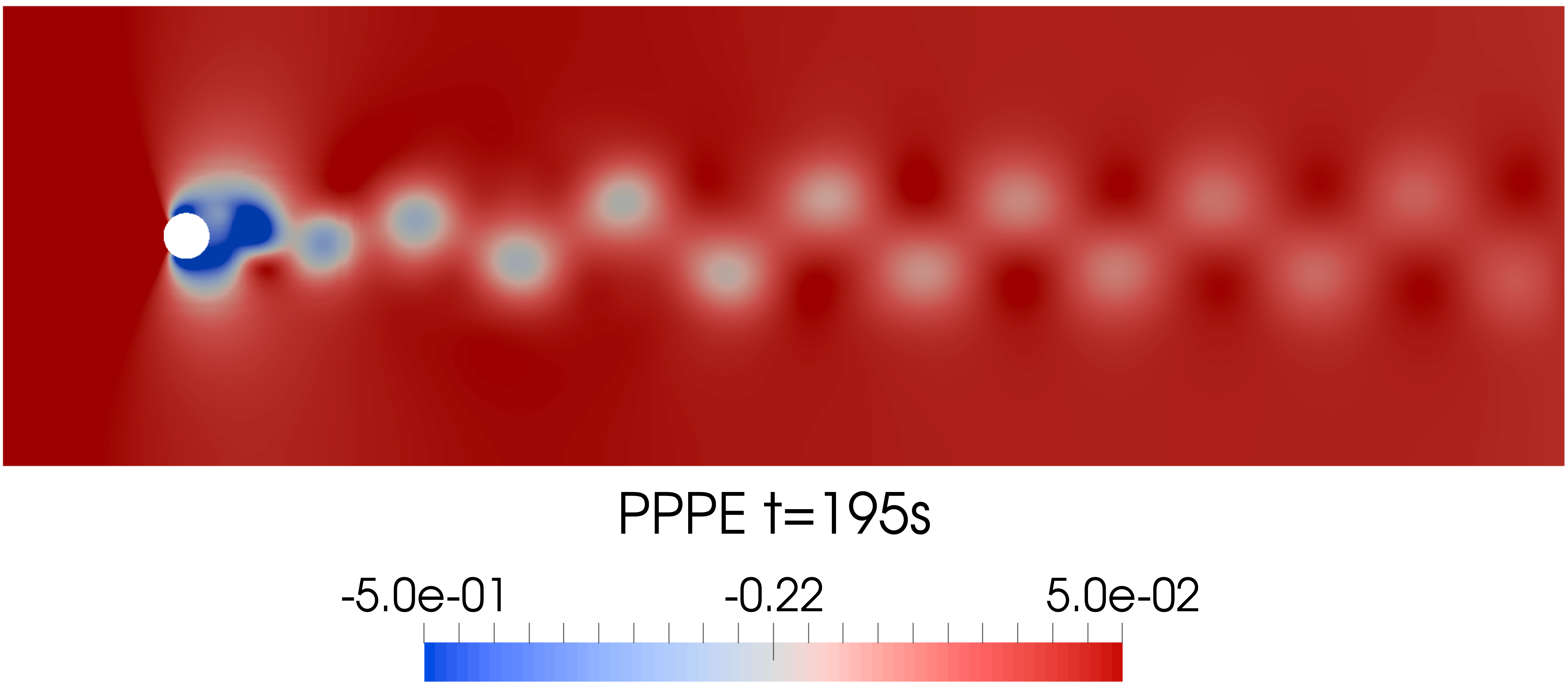}
\includegraphics[width=0.24\textwidth]{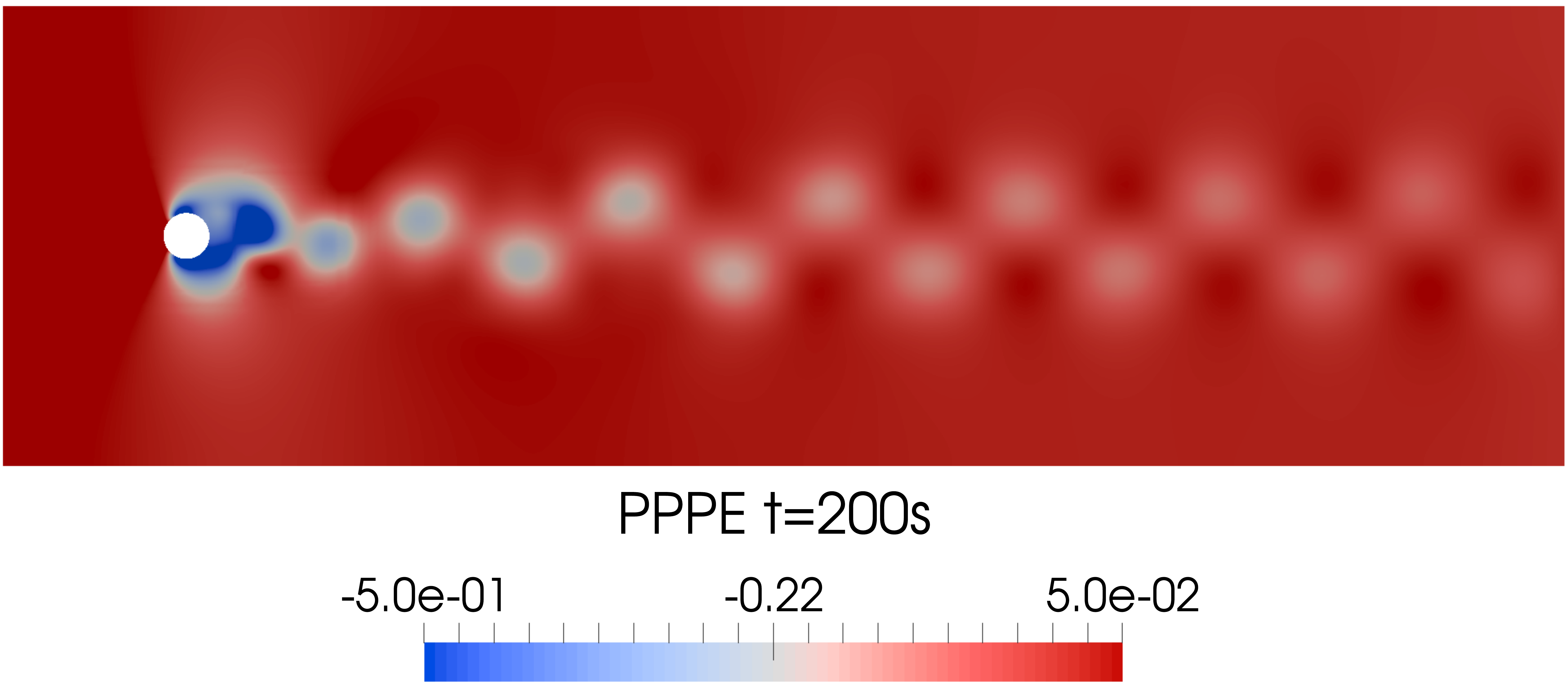}
\includegraphics[width=0.24\textwidth]{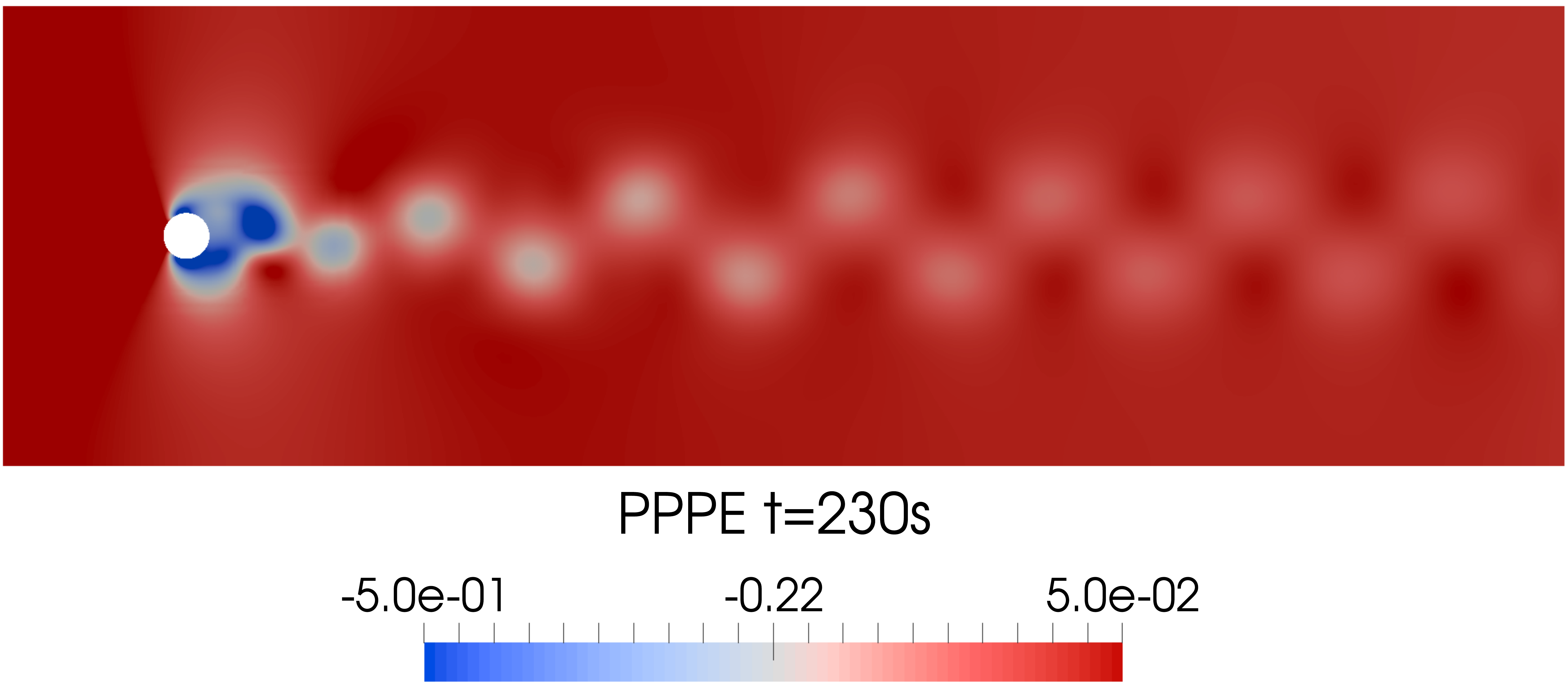}
\includegraphics[width=0.24\textwidth]{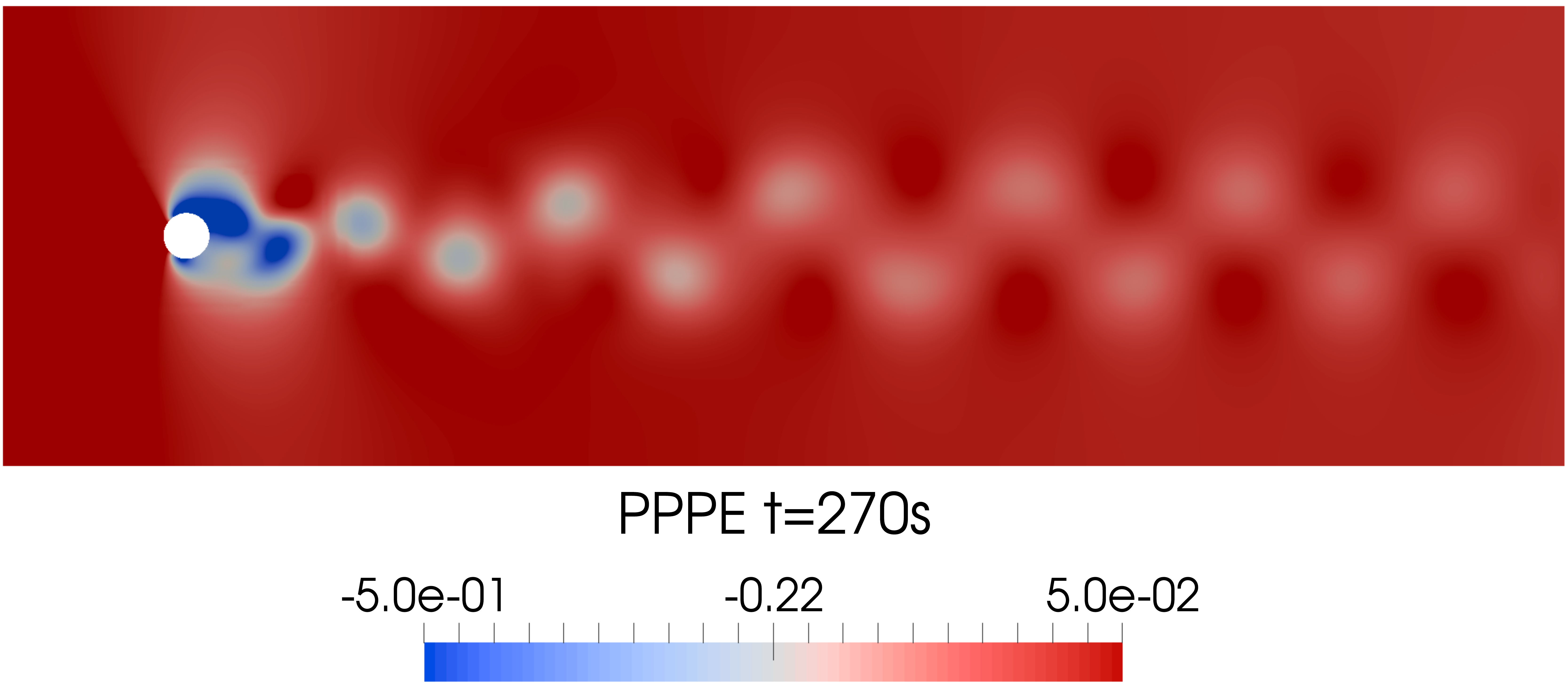}
\end{minipage} 
\caption{Comparison of the pressure field for high fidelity (PHF - first row), supremizer stabilised ROM (PSUP - second row) and pressure Poisson equation stabilised ROM (PPPE third row). The fields are depicted for different time instant equal to $t=195 \si{s},200 \si{s},230 \si{s}$ and $270 \si{s}$ and increasing from left to right. The ROM solutions are obtained with $15$ modes for velocity and $10$ modes for pressure, and only for the SUP-ROM with $12$ additional supremizer modes. The pressure magnitude is shown in the image legends.}\label{fig:pressure_cyl} 
\end{figure*}

\begin{figure*} 
\includegraphics[width=0.49\textwidth]{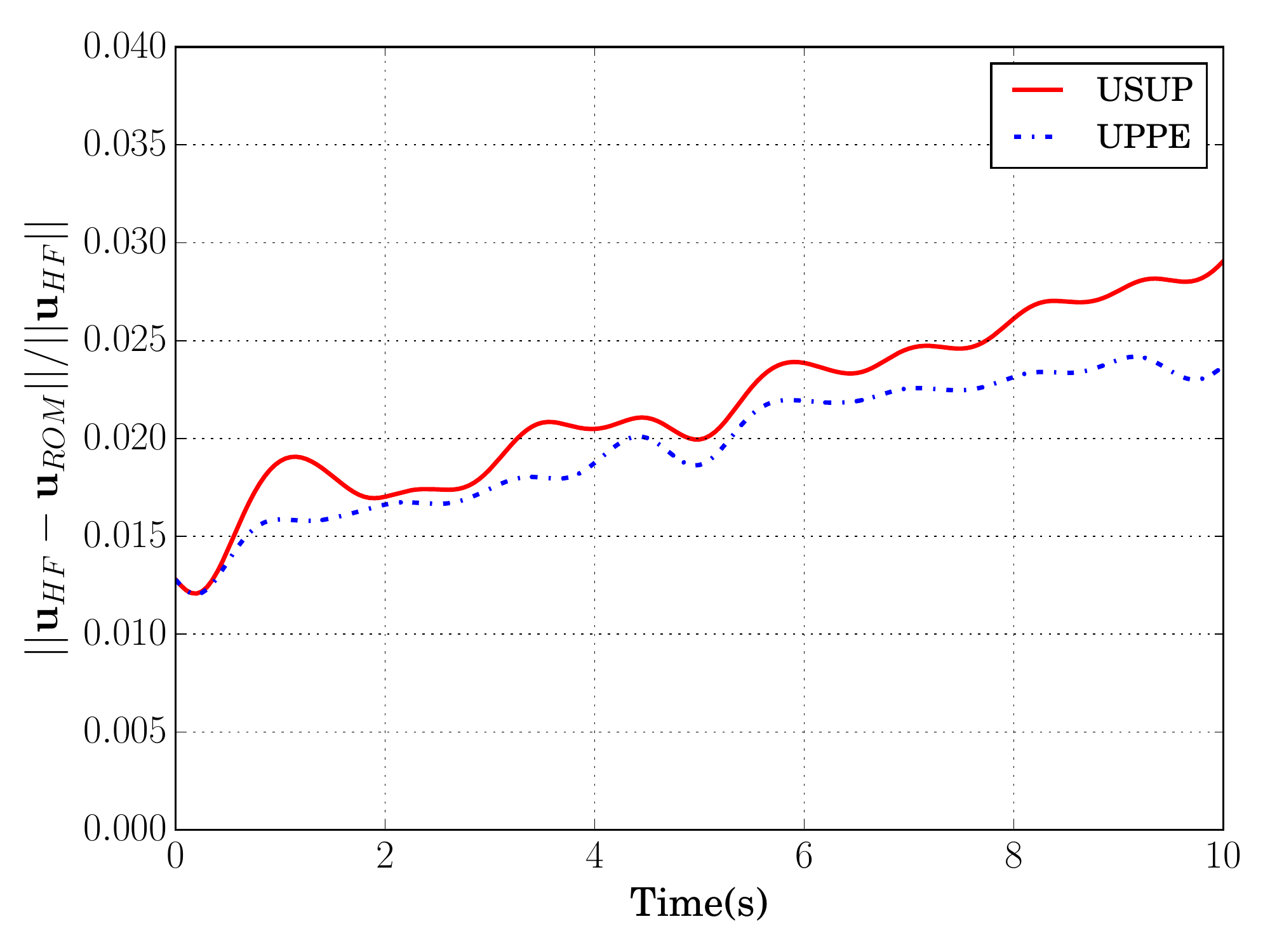}
\includegraphics[width=0.49\textwidth]{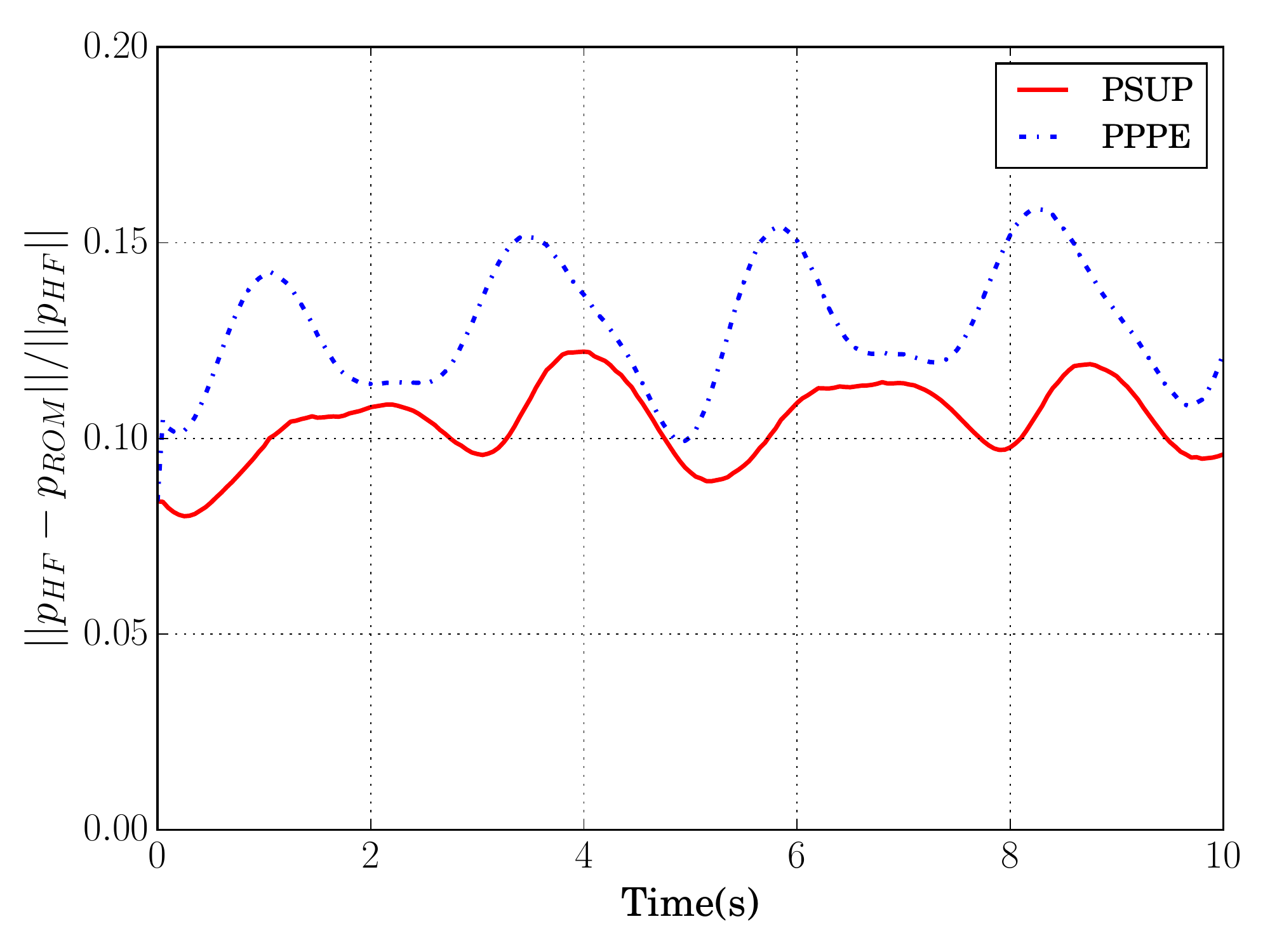}
\caption{Error analysis for the velocity (left plot) and pressure (right plot) fields in the cylinder example with $\nu = 0.005625$. The $L^2$ norm of the relative error is plotted over time for the two different models: with supremizer stabilisation (USUP and PSUP - continuous red line) and pressure Poisson equation stabilisation (UPPE and PPPE - dotted blue line). The ROM solutions are obtained with 15 modes for velocity, 10 modes for pressure and 12 modes for supremizers.}\label{fig:error_cilinder}
\end{figure*}

\begin{figure*}
\includegraphics[width=0.49\textwidth]{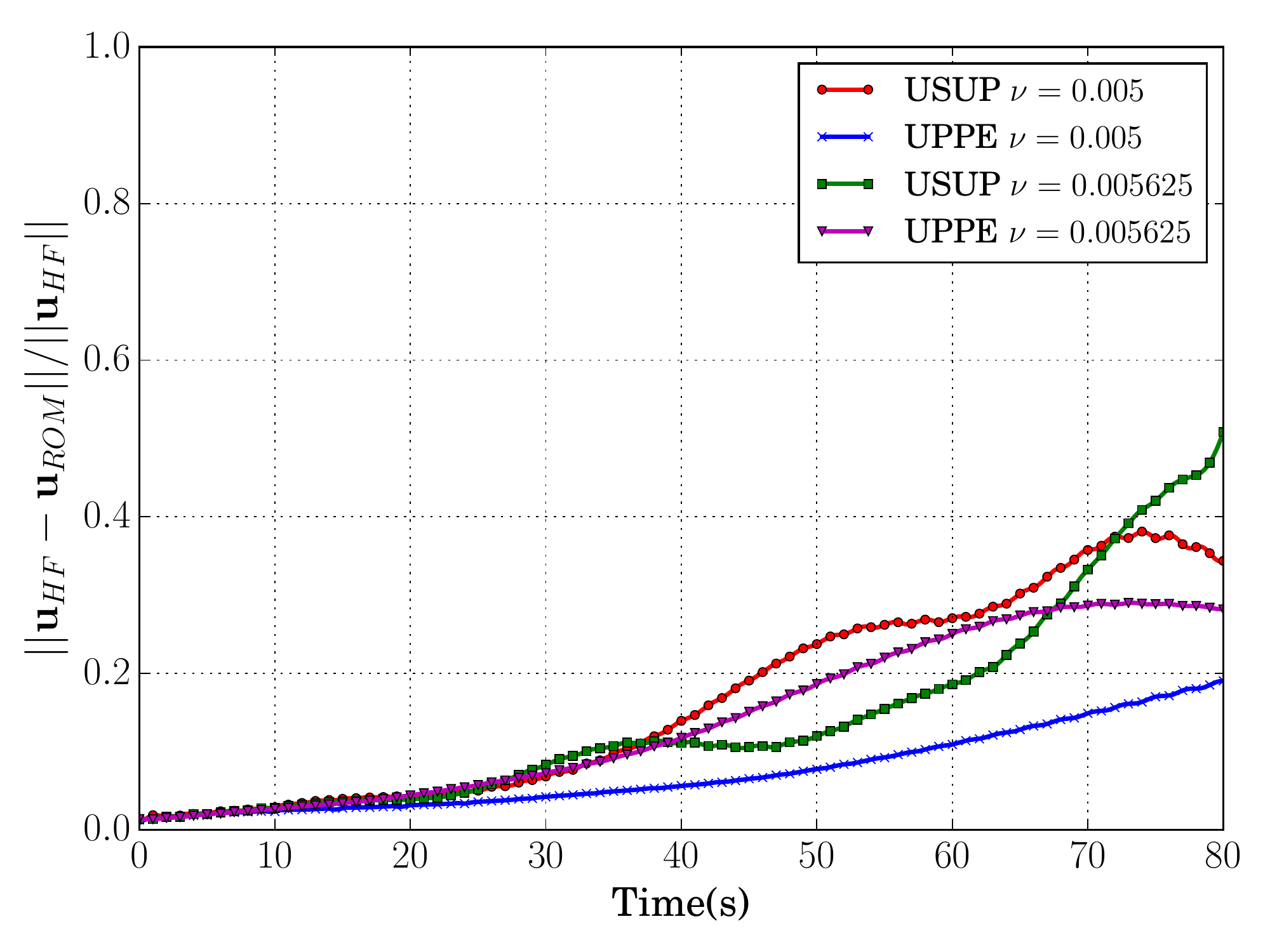}
\includegraphics[width=0.49\textwidth]{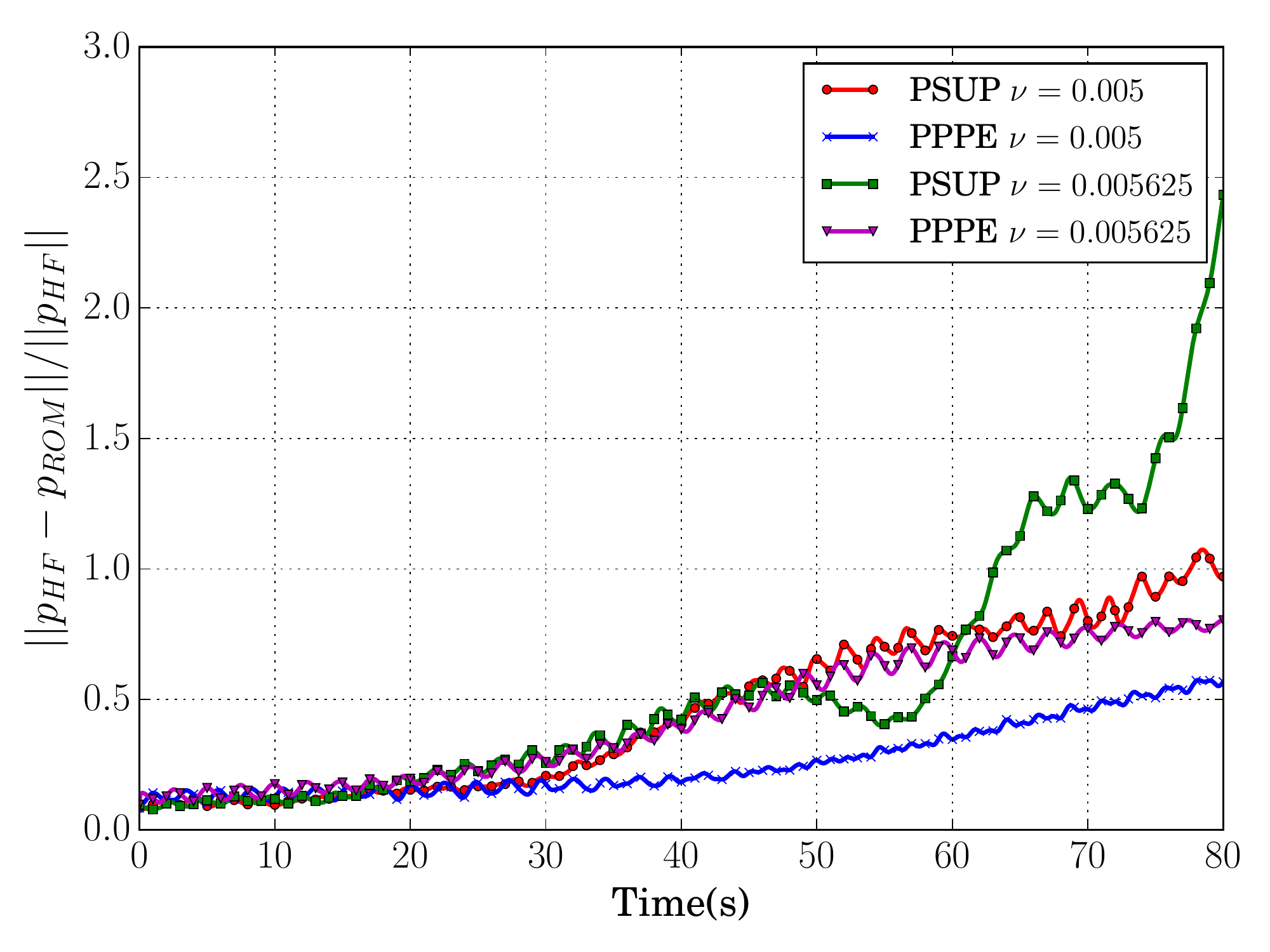}
\caption{Error analysis for the velocity (left plot) and pressure (right plot) fields in the cylinder example with $\nu = 0.005$ and $\nu = 0.005625$. The $L^2$ norm of the relative error is plotted over time for the two values of viscosity and for the two different models: with supremizer stabilisation (USUP and PSUP) and pressure Poisson equation stabilisation (UPPE and PPPE). The ROMs are obtained with 15 modes for velocity, 10 modes for pressure and 12 modes for supremizers. The time window is in this case wider respect to the one used for the generation of the reduced basis spaces ($\Delta T = 10$s)}\label{fig:error_cylinder_longer}  
\end{figure*}

To test the ROMs the results are compared against the full-order results for an intermediate value of the viscosity $\nu = 0.005625$, which is not included in the values of viscosity ($\nu = [0.005,0.00625,0.0075,0.00875,0.01]$) employed to generate the snapshots used to create the reduced basis spaces. The comparison has been performed on two different time windows. The first one covers 10 seconds of simulation and is coincident with the time window used for the generation of the snapshots,  the second one covers $80$ seconds of simulation and therefore is much wider with respect to the time window used for the generation of the snapshots. Figures~\ref{fig:velocity_cyl} and \ref{fig:pressure_cyl} show the comparison, for velocity and pressure respectively, between the results obtained with the full-order model, the SUP-ROM and the PPE-ROM. The fields are depicted at four different time instants equal to $t=195$s, $200$s, $230$s and $270$s. The first two time instants are respectively in the middle and at the end of the time window used to generate the snapshots while the two other time instants are outside of it. Figure~\ref{fig:error_cilinder} reports the $L^2$ norm of the relative error for velocity and pressure on the the $10$s wide time window for $\nu=0.005625$. The figure, also in this case, confirms the behaviour observed also in the cavity example: the SUP-ROM produces worse results for the velocity field but better results for the pressure field. Figure~\ref{fig:error_cylinder_longer} shows the same plots on a wider time window, and also for one of the value of viscosity ($\nu=0.005$) used to generate the full-order snapshots. For both ROMs, the relative error is increasing in time. Cross-referencing the data of figure~\ref{fig:error_cylinder_longer} with the plots of the figures~\ref{fig:velocity_cyl} and \ref{fig:pressure_cyl}, one can deduce that the increasing in the error is given, for the SUP-ROM, by the numerical instabilities that occur for long time integrations. In the last time step ($t=270$s), it is in fact possible to observe, for both velocity and pressure, a completely incorrect and non-physical flow pattern. For what concerns the PPE-ROM, instead, the flow pattern still looks regular and sufficiently similar to the high fidelity one but it is possible to observe a phase shift between the high fidelity and the ROM solution. The PPE-ROM, in fact, even though produces a still regular and physical pattern, has a period of vortex shedding which is slightly longer with respect the HF solution. To have a better idea about the behaviour of the different ROMs for long time integrations figure~\ref{fig:kyn_cyl} depicts the relative error of the total kinetic energy is plotted. It is well known that POD-Galerkin models are affected, in fact, by a blow-up energy issue \cite{couplet_sagaut_basdevant_2003,taddei2017}. It is not the objective of this work to deal with long time integration instabilities but it is worth checking which kind of pressure stabilisation method is likely prone to this issue. From the figure it is clear that the PPE-ROM accurately preserves the the total kinetic energy of the system. On the other hand, the SUP-ROM exhibits an oscillating behaviour with an increase of the total kinetic energy.      
\begin{table*}[] 
\centering
\caption{The table contains the cumulative eigenvalues for the cylinder problem. In the first, second and third columns are reported the cumulative eigenvalues for the velocity, pressure and supremizer fields respectively in function of the number of modes. In the last column is reported the value of the inf-sup constant, for the supremizer stabilisation case, for different different number of supremizer modes with a fixed number of velocity and pressure modes (15 modes for velocity and 10 modes for pressure)}\label{tab:cum_eig_cyl}
\label{my-label}
\begin{tabular}{ c | c | c | c | c}
N Modes & $\bm{u}$ & $p$ & $\bm{s}$ & $\beta$ \\ 
\hline 
1       & 0.390813 & 0.793239 & 0.921046 & 2.608e-04  \\
2       & 0.598176 & 0.85809  & 0.941746 & 4.492e-04  \\
3       & 0.802176 & 0.911636 & 0.961438 & 7.869e-03  \\
4       & 0.879096 & 0.934997 & 0.978072 & 1.662e-02  \\
5       & 0.949519 & 0.955578 & 0.98669  & 1.662e-02  \\
10      & 0.986025 & 0.992347 & 0.998307 & 1.098e-01  \\
15      & 0.995922 & 0.997994 & 0.999732 & 1.199e-01 
\end{tabular}
\end{table*}
\subsection{Comments on the results}
The two proposed numerical examples permit to draw some conclusions regarding the performances of the two different stabilisation methods. The SUP-ROM, demonstrates to produce better results for what concerns the pressure field and worse result for what concerns the velocity field. This fact may be justified by the additional and unnecessary (only in terms of correct representation of the velocity fields) supremizer modes that pollute the POD velocity space. The PPE-ROM, on the other hand, demonstrates do be more reliable for long time integrations. 

Table~\ref{tab:comp_time} shows the results in terms of computational costs. It is possible to deduce that both models, for both cases, permits to reach a considerable speed-up. For what concerns the cavity example both the offline and the online stages are computed in serial on one processor. On the other hand, in the cylinder example, the offline stage is performed in parallel with 6 processors while the online stage is still performed with a serial run on one processor.

The SUP-ROM demonstrates to be less efficient respect to the PPE-ROM. The SUP-ROM in fact, in comparison with the PPE-ROM, due to the additional supremizer modes, gives raise to a bigger reduced dynamical system. Both ROMs demonstrated to be able to capture with sufficient accuracy (especially from an engineering standpoint) the main features of the flow field for both velocity and pressure. The SUP-ROM demonstrates to be likely prone to instabilities issues for long time integrations.  

\section{Conclusions and perspectives}\label{sec:conclusions}
The main goal of this work was to compare and test the accuracy of two different pressure stabilisation strategies for POD-Galerkin ROMs based on a finite volume approximation. The ROMs are used to approximate the parametrised unsteady Navier-Stokes equations for moderate Reynolds numbers. 
The two analysed ROMs are based on the supremizer enrichment of the velocity space in order to meet the inf-sup condition and on to the exploitation of a pressure Poisson equation during the projection stage. The supremizer stabilisation is introduced here for the first time in a finite volume context and showed to effectively stabilise the resulting reduced system. It demonstrates moreover to be a valid alternative respect to the other stabilisation methods. The pressure Poisson equation is proposed here with an additional boundary condition for pressure that was neglected in previous works. 
Another goal of the article was also to test the behaviour of the two models for long time integrations. Concerning this aspect the PPE-ROM demonstrates to have better performances respect to the SUP-ROM.  
As future development, the interest is into higher Reynolds number and into turbulent flows. The attention will be in fact devoted to analyse the applicability of the proposed methods to turbulent flows. Moving to turbulent flows will be in fact essential to tackle real-world engineering problems. We will also further investigate the behaviour of the ROMs for long time integrations with the study of possible stabilisation techniques. The future interest is also into efficient methodologies for geometrical parametrisation. Reduced basis methods with FEM discretisation often employ domain decomposition and piecewise affine reference mappings. This decomposition is not trivial in a finite volume context because of the correlation among different parts of the domain introduced by the consistency requirement of the numerical fluxes.  
\begin{table*}[]
\centering
\caption{The table contains the computational time, for the supremizer (SUP) and the pressure Poisson equation (PPE) stabilisation techniques. In the cavity experiment the SUP-ROM is obtained with 10 modes for velocity, pressure and supremizers, while the PPE-ROM is obtained with 10 modes for pressure and supremizers. In the cylinder experiment the SUP-ROM is obtained with 15 modes for velocity, 10 for pressure and 12 for supremizers, while the PPE-ROM is obtained with 15 modes for velocity and 10 for pressure.}\label{tab:comp_time}
\label{my-label} 
\resizebox{0.48\textwidth}{!}{%
\begin{tabular}{ l | l | l | l }
& HF & SUP-ROM & PPE-ROM \\
\hline
Cavity Exp.  & $25 \mbox{min}$ & $7.64 \si{s}$ & $4.86 \si{s}$ \\
\hline
Cylinder Exp.& $18.5 \mbox{min} \times 6 \mbox{proc.}$ & $3.14 \si{s}$ & $0.971 \si{s}$
\end{tabular}}
\end{table*}
\section*{Acknowledgements} 
We acknowledge the support provided by the European Research Council Executive Agency by the Consolidator Grant project AROMA-CFD ``Advanced Reduced Order Methods with Applications in Computational Fluid Dynamics'' - GA 681447, H2020-ERC CoG 2015 AROMA-CFD and INdAM-GNCS projects. 
\section*{Appendix A. List of abbreviations and symbols}
\printnomenclature
\Urlmuskip=0mu plus 1mu\relax
\bibliographystyle{amsplain_gio}
\bibliography{bibfile} 
\end{multicols}
\end{document}